\documentclass[review,onefignum,onetabnum]{siamart190516}
\usepackage{amssymb,amsmath,version}
\usepackage{amssymb,version,graphicx,fancybox,mathrsfs,multirow}
\usepackage{threeparttable}
\usepackage{url,hyperref}
\usepackage{algpseudocode} 
\usepackage{verbatim}
\usepackage{color}
\usepackage{cases}
\usepackage{mathtools}
\usepackage{subfigure,float}
\usepackage{epstopdf}

\usepackage{booktabs}

\theoremstyle{plain}

\theoremstyle{remark}
\newtheorem{remark}{Remark}[section]

\theoremstyle{definition}
 
\newtheorem{example}{Example}[section]
\newtheorem{assum}{Assumption}[section]

\def\al{\alpha}
\def\fy{\varphi}
\def\Om{\Omega}

\def\Uad{\mathcal{A}}

\def\Ebb{\mathbb{E}}

\def\ul{u^{(\ell)}}

\def\wl{w^{(\ell)}}
\def\gl{g^{(\ell)}}
\def\fl{f^{(\ell)}}
\def\zl{z^{(\ell)}}

\def\twl{\widetilde{w}^{(\ell)}}
\def\fyl{\fy^{(\ell)}}
\def\II{{(\Omega)}}
\def\Omp{\Omega^{\prime}}
\def\II{(\Omega)}

\renewcommand{\d}{\mathrm{d}}

\def\II{(\Omega)}

\allowdisplaybreaks[1]
\begin{document}

%\subtitle{Subject Section}

\title{Finite element approximation for quantitative photoacoustic tomography in a diffusive regime}
\author{
Giovanni S.\ Alberti\thanks{Machine Learning Genoa Center (MaLGa), 
Department of Mathematics, Department of Excellence 2023-2027, University of Genoa, Via Dodecaneso 35, 16146 Genova, Italy. (\texttt{giovanni.alberti@unige.it}).}
\and
Siyu Cen\thanks{Department of Applied Mathematics, The Hong Kong Polytechnic University, Kowloon, Hong Kong, China. (\texttt{21037194r@connect.polyue.hk, zhizhou@polyu.edu.hk}).} \and Zhi Zhou$^\dag$}
\date{\today}

\maketitle
\begin{abstract}
In this paper, we focus on the numerical analysis of quantitative photoacoustic tomography. Our goal is to reconstruct the optical coefficients, i.e., the diffusion and absorption coefficients, using multiple internal observational data. The foundation of our numerical algorithm lies in solving an inverse diffusivity problem and a direct problem associated with {  an} elliptic equation. The stability of the inverse problem  depends critically on a non-zero condition in the internal observations, a condition that can be met using randomly chosen boundary excitation data. Utilizing these randomly generated boundary data, we implement an output least squares formulation combined with finite element discretization to solve the inverse problem. In this scenario, we provide a rigorous error estimate in   $L^2(\Omega)$ norm for the numerical reconstruction using a weighted energy estimate, inspired by the analysis of a newly proposed conditional stability result. The resulting error estimate provides valuable guidance for the \textsl{a priori} selection or adaptive adjustment of suitable regularization parameters and discretization mesh sizes. Several numerical experiments are presented to support our theoretical results and illustrate the effectiveness of our numerical scheme.
\vskip5pt
\textbf{Key words}: quantitative photoacoustic tomography,  inverse diffusivity problem, random boundary excitation data,
least-squares formulation, finite element approximation,  error estimate. 
\end{abstract} 

\maketitle

\vskip10pt

\section{Introduction}\label{sec:intro}
Photoacoustic tomography (PAT) is a biomedical imaging technique that combines the principles of optical imaging and ultrasound to produce high-resolution images of tissues within the body \cite{li2009photoacoustic,wang2009multiscale}. It offers unique advantages by capturing the functional and structural characteristics of tissues, making it particularly useful for medical diagnostics, including cancer detection, monitoring of vascular diseases, and studying brain functions. The first inverse problem in PAT concerns the reconstruction of the deposited optical energy from the time-dependent boundary measurement of the acoustic pressure. Explicit inversion formulas exist for a large class of geometries of interest, when the problem is in free space, with constant sound speed, and without accounting of acoustic attenuation. See some related discussion in \cite{Finch:2004,Kuchment:2008,Patch:2007,Kunyansky:2007} and the references therein. The free-space model of acoustic propagation ignores the boundary effects caused by the transducers. These can be taken into account by considering the problem in a bounded domain with suitable boundary conditions, see e.g. \cite{kunyansky-holman-cox-2013,Alberti2018,campodonico-2021}. Independently of the model, this first inverse problem is only moderately ill-posed: the model is based on the wave equation, which propagates { the} singularities, and allows for high-quality reconstructions of the optical energy inside the domain of interest.

In this paper, we assume that the aforementioned first step is done and that the deposited optical energy is known. Then we consider the second step of PAT, called quantitative photoacoustic tomography (QPAT), i.e., to recover simultaneously the diffusion coefficient and the absorption coefficient from the deposited optical energy.
We consider the case where radiation propagation is approximated by a second-order elliptic (diffusion) equation \cite{CaseZweife:1967,Arridge:1999}:
\begin{equation}\label{eqn:QPAT}
    \left\{
        \begin{aligned}
-\nabla \cdot (D(x) \nabla u)+\sigma(x) u&=0 && \text { in  $\Omega$,} \\
u&=g && \text { on $\partial \Omega$.}
        \end{aligned}
    \right.
\end{equation}
Here, $\Omega$ is a bounded Lipschitz domain in $\mathbb{R}^d$ ($d=2,3$) with boundary $\partial \Omega$.
%Here $\Omega$ is a convex polyhedral domain in $\mathbb{R}^d$ with boundary $\partial \Omega$.
The optical coefficients $(D(x), \sigma(x))$, with $D(x)$ being the diffusion coefficient and $\sigma(x)$ the absorption coefficient, are assumed to {be bounded and positive.}
The QPAT inverse problem consists of recovering $D(x)$ and $\sigma(x)$ from the internal observation of the optical energy
$$
H(x)=\sigma(x) u(x) \quad \text { for all } x \in \Omega.
$$

The problem of QPAT has been extensively studied in the literature. Since the inverse problem involves multiple parameters ($D$ and $\sigma$), a common method uses multiple illuminations $g$ to generate various optical energies $H$ and reconstruct the unknown parameters. In \cite{Bal:2010,Bal:2011}, the authors propose a decoupled procedure and prove the uniqueness and {the} H{\"o}lder stability for the inverse problem. The decoupled scheme relies on the following observation: if $u_1, u_2$ are two solutions to equation \eqref{eqn:QPAT} corresponding to illuminations $g_1, g_2$ respectively, then the quotient $u = u_2/u_1 = H_2 / H_1$ satisfies the following elliptic equation with one parameter:
\begin{equation}\label{eqn:IDP0}
    \left\{
        \begin{aligned}
            -\nabla \cdot (q \nabla u) & = 0  && \text{in } \Omega, \\
            u & = g  && \text{on } \partial\Omega,
        \end{aligned}
    \right.
\end{equation}
where $q = D u_1^2$ and $g = g_2/g_1$. Thus, the problem of QPAT is solved by a two-step procedure. The first step is to solve an inverse diffusivity problem (IDP) of recovering $q$ given $u$ and the boundary value $q|_{\partial\Omega}$.
After obtaining $q = D u_1^2$, the second step is solving a direct problem:
\begin{equation}\label{eqn:step2}
    \left\{
        \begin{aligned}
            -\nabla \cdot (D u_1^2 \nabla (1/u_1)) & = H_1  && \text{in } \Omega, \\
            1/u_1 & = 1/g_1  && \text{on } \partial\Omega,
        \end{aligned}
    \right.
\end{equation}
to find $u_1$ and hence determine $D$ and $\sigma$.

It is important to highlight that the following non-zero condition is crucial for the IDP:
\begin{align}\label{eqn:nonzero-0}
|\nabla u(x)| \ge C_0> 0 \quad \text{for all} \quad x \in \Omega.
\end{align}
There are several approaches for constructing a boundary illumination $g$ such that this condition holds. When $d=2$, the works \cite{Alessandrini:1986, Alessandrini:1994} provide a simple criterion for choosing a special boundary illumination $g$ that guarantees the non-zero condition. Roughly speaking, the graph of $g$ should have a single maximum point, a single minimum point, and be monotone in between.
For dimensions $d \ge 3$, ensuring the non-zero condition becomes more challenging \cite{alberti-bal-dicristo-2017}{: the suitable boundary values cannot be chosen independently of the unknown parameters of the PDE}. In \cite{Bal:2010}, the author uses the method of complex geometrical optics to construct boundary data $g$ satisfying the non-zero condition. However, this construction is not very explicit and depends on the interior values of the unknown coefficient $q$.
{To address this issue, bypassing the non-zero condition, \cite{Alessandrini:2017} derives a $\alpha$-H{\"o}lder stability result when the boundary data $g$ is quantitatively unimodal, which can be viewed as an extension of monotonicity in higher dimensions. However, the parameter $\alpha$ is not explicit and the construction of the boundary values is not easily implementable numerically.} 
%{We note that it is possible to obtain {the} $\alpha$-H{\"o}lder stability for the inverse problem even without requiring \eqref{eqn:nonzero-0}, provided the illuminations are suitably chosen \cite{Alessandrini:2017}. However, the parameter $\alpha$ is not explicit and the construction of the boundary values is not easily implementable numerically.}
Recently, \cite{Alberti:2022, Alberti:2022Runge} considered using random boundary illuminations and proved that the corresponding solutions will satisfy the non-zero condition with overwhelming probability. {This approach overcomes the drawbacks of the previous deterministic methods, as it avoids restrictive constraints on the boundary illuminations and provides an explicit construction that is more compatible with practical situations}.

In this project, we aim to develop a reconstruction scheme and to establish the approximation error for the diffusion coefficient $D$ and the absorption coefficient $\sigma$ in equation \eqref{eqn:QPAT} from multiple internal observations corresponding to carefully designed random boundary illuminations. 
The numerical analysis of parameter identification problems for elliptic equations, particularly the inverse diffusivity problem (IDP) using internal data, has been extensively studied in the literature \cite{Richter:1981inverse,Richter:1981numerical,Falk:1983,Wang:2010,Jin:2021error}. In \cite{Falk:1983}, one of the earliest works, Falk proposes a discrete least squares scheme to solve the IDP from a single observation and analyzes the reconstruction error under the \textit{a priori} non-zero condition \eqref{eqn:nonzero-0}. He derives {an}  approximation rate $O(h^r+h^{-2}\delta)$ in $L^2(\Omega)$ norm, where $r$ is the polynomial degree of the finite element space and $h$ is the mesh size. However, the constant appearing in the reconstruction rate increases exponentially  with respect to the $W^{2,\infty}(\Omega)$ norm of the solution. In \cite{Wang:2010}, the authors provide an improved convergence result under a stronger non-zero condition: $a_0 |\nabla u(x)|^2 > \max(f(x),0)$ almost everywhere in $\Omega$, where $a_0$ is a chosen constant and $f$ is the source term. Inspired by the stability analysis in \cite{Bonito:2017}, the work in \cite{Jin:2021error} develops a new error bound using a weighted energy estimate with a special test function. The estimate utilizes a much weaker non-zero condition, $q|\nabla u|^2+fu\ge C_0>0 $, and does not have the exponentially dependence on $\|u\|_{W^{2,\infty}(\Omega)}$. Recently, this approach was extended in \cite{Cen:2023} to reconstruct two parameters, $D$ and $\sigma$, with two internal measurements, $u_1$ and $u_2$, generated by different source terms $f_1$ and $f_2$, respectively. Nevertheless, the positivity condition used in these works requires a strictly positive or negative source term. Note that such a requirement is also needed in the numerical analysis of machine learning approaches \cite{CenJinQuanZhou:IMA2024, JinLiQuanZhou:SIIM2024}, as well as for inverse potential problems \cite{JinLuQuanZhou:2023,ZhangZhangZhou:2022,jin2024stochastic}.

In this work, we investigate the problem of QPAT raising in {a} practical scenario, where the source term vanishes and the measurement $H=\sigma u$ is generated by a boundary illumination. Compared with the work in \cite{Cen:2023}, the vanishing source term makes the required positivity condition fail in general,
and the measurement, which is the product between the function $u$ and the absorption coefficient $\sigma$, is more involved. In order to have a H{\"o}lder type stability, we employ specially designed random boundary illuminations \cite{Alberti:2022, Alberti:2022Runge}, and apply the weighted energy estimate with special test functions \cite{Bonito:2017, Jin:2021error}. We then discuss the numerical inversion formula and analyze the approximation error for the reconstruction. One popular reconstruction approach is to reformulate the IDP \eqref{eqn:IDP0} as a transport equation with variable $q$ \cite{Bal:2010,Bal:2011}. This approach is non-iterative and hence efficient for computation. However, it requires the non-zero condition to hold on the whole domain $\Omega$, while in our approach (see Proposition~\ref{prop:nonzero_condition}) the non-zero condition holds only locally for a specific boundary illumination. On the other hand, the least square formulation allows one to naturally incorporate  the local non-zero property into the error analysis.  Therefore, in this paper, we consider the least square fitting approach with a regularization term for the QPAT reconstruction. Motivated by the stability estimate, we employ weighted energy estimate with a special test function to analyze the approximation error in terms of the discretization mesh size $h$, the noise level $\delta$, and the regularization parameter $\alpha$. Our approach employs several technical tools, including the decoupled procedure for QPAT, the weighted energy estimate, the non-vanishing gradient property, and \textit{a priori} estimates for the finite element approximation.

The rest of the paper is organized as follows. In Section~\ref{sec:IDP}, we discuss the choice of random boundary illuminations and show the H{\"o}lder type stability of the inverse diffusivity problem under the non-zero condition. We also propose an iterative reconstruction algorithm and study the finite element approximation error. In Section~\ref{sec:QPAT}, we establish the numerical inversion scheme for QPAT and analyze the discrete approximation error. Numerical experiments are presented in Section~\ref{sec:numer} to validate the theoretical results. Throughout, we denote the standard Sobolev spaces of order $s$ by $W^{s,p}(\Omega )$ for any real $s \geq  0$ and  $p \geq  1$, equipped with the norm $\|  \cdot  \|_{  W^{s,p}(\Omega )}$. When $p=2$, we use the notation $H^{s}(\Omega)=W^{s,p}(\Omega)$. Moreover, we write $L^p(\Omega )$ with the norm $\|  \cdot  \|_{L^p(\Omega )}$ if $s = 0$. The spaces on the boundary $\partial \Omega$  are defined similarly. The notation $(\cdot , \cdot )$ denotes the standard $L^2(\Omega )$ inner product. The notation $a\lesssim b$ indicates that $a\le C b$  holds for some constant $C>0$, where $C$ is independent of the relevant parameters under consideration. We denote by $c$ and $C$ generic constants that are not necessarily the same at each occurrence, but are always independent of the noise level, the discretization parameter, and the penalty parameter.

\section{Inverse diffusivity problem}\label{sec:IDP}
In this section, we consider the inverse diffusivity problem associated to the second-order elliptic equation 
\begin{equation}\label{eqn:IDP}
    \left\{
        \begin{aligned}
            -\nabla \cdot (q \nabla w) & =0 && \mbox{in }\Om, \\
            w & =g && \mbox{on } \partial\Om.
        \end{aligned}
    \right.
\end{equation}
Let $\Omp\Subset\Om$ be a given Lipschitz subdomain and suppose that the exact diffusion coefficient $q^\dagger(x)$ is known for all $x \in \Om\setminus \Omp$.
The diffusion coefficient is assumed to be in  the following admissible set:
\begin{equation}\label{eq:Aq}
    \Uad_q=\{q\in H^1\II \mbox{ : }0<\Lambda_q^{-1} \le q\le \Lambda_q\mbox{ a.e. in }\Om   \mbox{, }  q=q^\dag  \mbox{ in }\Om\setminus\Omp  \},
\end{equation}
with an a priori known positive constant $\Lambda_q$. Moreover, we assume that the coefficient and {the} boundary data satisfy the following assumption.

\begin{assum}\label{assum:parameter_regularity}
Let $\Omega$ be a bounded Lipschitz domain in $\mathbb{R}^d$ and $\Omp\Subset\Omega$ be a given Lipschitz subdomain. We assume that the exact diffusivity coefficient $q^\dag \in C^{0,1}(\overline{\Omega}) \cap \Uad_q$. Further, we let $\gl$ (with $\ell=1,\dots,L$) denote {the} boundary data, which are taken as independent and identically distributed random variables in $H^{\frac{1}{2}}(\partial\Omega)$ satisfying the expansion
 \begin{equation}\label{eqn:reps_g}
        g^{(\ell)}=\sum_{k=1}^{M} a_k^{(\ell)} e_k,\qquad \ell=1,\dots,L,
\end{equation}
where $M$ is a given positive integer, $\{e_k\}_{k=1}^{\infty}$ is a fixed orthonormal basis of $H^{\frac{1}{2}}(\partial\Omega)$  and $a^{(\ell)}_{k}\sim N(0,\theta_k^2)$ are independent real Gaussian variables, with $\theta_k>0$ for every $k$ and $\sum_{k\ge 1} \theta_k<\infty$.
\end{assum}

\begin{remark}\label{rmk:reg-0}
Let $\wl(q^\dag)$ denote the solution to the elliptic problem \eqref{eqn:IDP} associated with the diffusion coefficient $q^\dag$ and the boundary excitation $\gl$. 
Under the regularity assumption, classical elliptic regularity theory (\cite[Theorem 5.20]{Giaquinta:2013} and  \cite[Theorem 8.8]{Gilbarg:1977}) implies that the corresponding solution to the elliptic equation \eqref{eqn:IDP} satisfies $\wl(q^\dag)\in C_{loc}^{1,\kappa}\II\cap H^1\II$ for all $\kappa\in (0,1)$.
\end{remark}

The inverse diffusivity problem (IDP) consists of recovering the diffusion coefficient in $\Omega'$ from the multiple internal observations $\wl(x;q^\dag)$ for all $x\in \Omp$, where $\ell=1,2,\ldots,L$. 
With the above choice of $\gl$, by using the result of \cite{Alberti:2022} we have the following non-zero condition, which is crucial for stability and error estimates. {Note that the non-zero condition holds with overwhelming probability, and the constants $C_0, C_1,$ and $C_2$ require only minimal \textit{a priori} information regarding $q$.}
\begin{proposition}\label{prop:nonzero_condition}
Suppose that Assumption~\ref{assum:parameter_regularity} holds. Take $\nu\in \mathbb{R}^d$ with $|\nu|=1$. Then, with a probability greater than 
\begin{equation}\label{eqn:probability_nonzero}
1-  L^d \exp\left(-C_1 L\right)- L\exp\left(-C_2  M\right),
\end{equation}
the following non-zero condition holds
\begin{equation}\label{eqn:nonzero_condition}
    \max_{\ell=1,\dots,L}|\nabla \wl(x)\cdot \nu |\ge  C_0,\quad x\in \Omega^{\prime},
\end{equation}
and the random boundary data {can be bounded by}
\begin{equation}\label{eqn:bounded_bdry}
    \max_{\ell=1,\dots,L} \|g^{(\ell)}\|_{H^{\frac{1}{2}}(\partial\Omega)}\le L^{\frac{1}{2}}.
\end{equation}
Here  
$\wl$ (with $\ell=1,\dots,L$) is the solution to \eqref{eqn:IDP} corresponding to the boundary illumination $\gl$. The positive constants $C_0$, $C_1$ and  $C_2$  depend only on $\Om$, $\Omp$, 
$\{\theta_k\}$, $\{e_k\}$, $\Lambda_q$ and  $\|q\|_{C^{0,1}(\overline{\Omega})}$.
\end{proposition}
\begin{proof}
All the constants appearing in the proof will depend only on $\Om$, $\Omp$, 
$\{\theta_k\}$, $\{e_k\}$, $\Lambda_q$ and  $\|q\|_{C^{0,1}(\overline{\Omega})}$.    Let $\overline{w}^{(\ell)}$ be the solution to \eqref{eqn:IDP} with boundary data
    \begin{equation*}
        \overline{g}^{(\ell)}=\sum_{k=1}^{\infty} a_k^{(\ell)} e_k, \quad \ell=1,\dots,L,
    \end{equation*} 
    where  {$\{e_k\}_{k=1}^{\infty}$ and $a^{(\ell)}_{k}$ are as in Assumption~\ref{assum:parameter_regularity}.}  
    By \cite[Theorem~1]{Alberti:2022} {(with the choice $\zeta(u) = \nabla w\cdot \nu$, $n=1$ and $N=L$, as a minor variation of \cite[Example~2]{Alberti:2022}) and \cite[Lemma~5]{Alberti:2022})},  with probability greater than $ 1- L^d \exp\left(-C_1 L\right)$, we have the following non-zero condition
    \begin{equation*} 
        \max_{\ell=1,\dots,L}|\nabla \overline{w}^{(\ell)}(x)\cdot \nu |\ge 2C_0,\quad x\in \Omega^{\prime}
    \end{equation*}
and 
    \begin{equation*}
        \max_{\ell=1,\dots,L} \|\overline{g}^{(\ell)}\|_{H^{\frac{1}{2}}(\partial\Omega)}\le  L^{\frac{1}{2}}/2.
    \end{equation*}
    
   Now we estimate the difference between $\overline{g}^{(\ell)}$ and the truncated boundary values $g^{(\ell)}$. 
     We view $\| \overline{g}^{(\ell)}-g^{(\ell)} \|_{H^{\frac{1}{2}}(\partial\Om)} $ as a random variable.  Since $a_k^{(\ell)}\sim N(0,\theta_k^2)$ and $ e_k$ are orthonormal in $H^{\frac{1}{2}}(\partial\Omega)$, the moment generating function satisfies for all $\lambda\in \mathbb{R}$:
    \begin{align*}
       \Ebb\exp\left(\lambda^2\| \overline{g}^{(\ell)}-g^{(\ell)} \|_{H^{\frac{1}{2}}(\partial\Om)}^2\right)=\Ebb\exp\left(\lambda^2\sum_{k=M+1}^{\infty}(a_k^{(\ell)})^2\right)=\exp\left(\lambda^2\sum_{k=M+1}^{\infty}\theta_k^2\right).
    \end{align*}
    The condition $\sum_{k=1}^{\infty}\theta_k<\infty$ implies that $ \sum_{k=M+1}^{\infty}\theta_k^2\le CM^{-1}$. By \cite[Proposition~2.5.2]{Vershynin:2018}, we have
    \begin{equation*}
        \mathbb{P}\left(\|\overline{g}^{(\ell)}-g^{(\ell)} \|_{H^{\frac{1}{2}}(\partial\Om)} \ge t \right)\le 2\exp\left(-C^2 t^2M\right),\quad \forall t\ge 0,\, \ell=1,\dots,L.
    \end{equation*}
    Thus, with probability greater than $ 1-2L\exp\left(-C^2 t^2 M\right)$, we have 
    \begin{equation*}
        \| \overline{g}^{(\ell)}-\gl \|_{H^{\frac{1}{2}}(\partial\Om)}\le t,\qquad    \ell=1,\dots,L.
    \end{equation*}
    Hence,  elliptic regularity yields
    \begin{equation*}
        \|\overline{w}^{(\ell)}-\wl\|_{C^1(\overline{\Omega'})}\le \widetilde{C}t,\qquad \ell=1,\dots,L.
    \end{equation*}
With the choice  $t= \min\{C_0/\widetilde{C}, L^{\frac{1}{2}}/2\}$, we have $\|\overline{w}^{(\ell)} -\wl \|_{C^1(\overline{\Omega'})}\le C_0$ and $\|\gl\|_{H^{\frac{1}{2}}(\partial\Omega)}\le L^{\frac{1}{2}}$. Let $C_2=C^2 t^2$, with a probability greater than
 \begin{equation*}
    1- L^d \exp\left(-C_1 L\right)-2L\exp\left(-C_2  M\right),
\end{equation*}
the  non-zero condition \eqref{eqn:nonzero_condition} and the upper bound on the boundary values \eqref{eqn:bounded_bdry} hold. 
\end{proof}

\subsection{Conditional  Stability}\label{subsec:stability}
In this part, we derive a useful conditional stability estimate in Sobolev spaces for the inverse diffusivity problem.
According to the non-zero condition \eqref{eqn:nonzero_condition} and the smoothness of the solutions $\wl \in C_{loc}^{1,\kappa}\II\cap H^1\II$,
there exist  
open sets $\Omega_\ell$, $\ell=1,\ldots,L$, covering $\Omp$ such that
\begin{equation}\label{eqn:cover_omega}
    \Omp\subset\bigcup_{\ell=1}^{L}\overline\Om_{\ell}\quad\text{where} \quad | \nabla \wl \cdot \nu| > C_0/2 \text{ for all} ~ x\in \Om_{\ell}.
    %=\{x\in \Omp: |(\ul_1\nabla\ul_2-\ul_2\nabla \ul_1)\cdot \nu|\ge C_0\}.
\end{equation}

\begin{theorem}\label{thm:stability}
Suppose the diffusion coefficient $q$ and the boundary terms $\gl$  (with $\ell=1,\dots,L$) satisfy  Assumption~\ref{assum:parameter_regularity}, and let $\widetilde{q}\in \mathcal{A}_q$  be a perturbation.  
Let $\wl$ and $\twl$ be the corresponding solutions to \eqref{eqn:IDP} with parameters $q$  and $\widetilde{q} $, respectively. Then, with a probability greater than
\eqref{eqn:probability_nonzero}, the following stability estimate holds:  
    \begin{equation}\label{eqn:stability}
        \|q-\widetilde{q}\|_{L^2\II} \le C C_0^{-1} L^{\frac{1}{4}} \Big(\sum_{\ell=1}^{L}\|\wl-\twl\|_{H^1(\Omp)}\Big)^\frac12.
    \end{equation}
     Here  $C> 0$ is a  constant depending only on $\Om$, $\Omp$, $ \Lambda_q$ and $\|q\|_{C^{0,1}(\overline{\Omega})}$, and $C_0$ is the lower bound of the non-zero condition given in \eqref{eqn:nonzero_condition}.  
\end{theorem}
\begin{proof}
{With an abuse of notation, several positive constants depending only on $\Om$, $\Omp$, $ \Lambda_q$ and $\|q\|_{C^{0,1}(\overline{\Omega})}$ will be denoted by the same letter $C$.}
By Proposition~\ref{prop:nonzero_condition}, with overwhelming probability \eqref{eqn:probability_nonzero}, both the non-zero condition \eqref{eqn:nonzero_condition} and the uniform bound \eqref{eqn:bounded_bdry} are satisfied.
Then for a given $\ell\in\{1,\dots,L\}$,  for any test function $\fyl \in H_0^1(\Omega)$,  integration by parts in \eqref{eqn:IDP} yields
\begin{equation}\label{eqn:stab_test}
    \big( (q-\widetilde{q})\nabla \wl,\nabla \fyl \big)= \big(  \widetilde{q} \nabla (\twl-\wl),\nabla \fyl \big).
\end{equation}
    Furthermore,  multiplying    both sides of \eqref{eqn:IDP} by$\frac{q-\widetilde{q}}{q}\fyl$  and applying integration by parts, we obtain
\begin{equation*}
    0=\big( q\nabla \wl,\nabla \frac{(q-\widetilde{q})\fyl}{q} \big)=  \big( q\fyl\nabla \wl,\nabla \frac{(q-\widetilde{q})}{q} \big)+\big( q \frac{(q-\widetilde{q})}{q} \nabla \wl,\nabla \fyl\big),
\end{equation*}
and hence 
\begin{equation}\label{eqn:stab_integrate}
    \big(  (q-\widetilde{q})  \nabla \wl,\nabla \fyl\big) =\frac{1}{2}\big(  (q-\widetilde{q})  \nabla \wl,\nabla \fyl\big)-\frac{1}{2}\big( q\fyl\nabla \wl,\nabla \frac{(q-\widetilde{q})}{q} \big).
\end{equation}
Now, we choose the test function $\fyl=(q-\widetilde{q})\wl/q$. Since $q=\widetilde{q}$ on $\Omega\setminus \Omp$, $\fyl$ vanishes on $\partial\Omega$. Noting that $q,\widetilde{q}\in \Uad_{q}$ and $\wl\in C^{1,\kappa}(\overline{\Omega'})$, we conclude that $\fyl\in H_0^1(\Omega)$, with 
\begin{align*}
    \|\fyl\|_{L^2(\Omega)}= \|(q-\widetilde{q})\wl/q\|_{L^2(\Omega)}\le 2\Lambda_q^2\|\wl \|_{L^2(\Omega)}\le CL^{\frac{1}{2}}
\end{align*}
and 
\begin{align*}
    &\|\nabla \fyl\|_{L^2(\Omega)}=\left\|\frac{q\nabla[(q-\widetilde{q})\wl]-(q-\tilde q)\wl \nabla q }{q^2}\right\|_{L^2(\Omp)}\\
    \le &  \Lambda_q^2 \Big( \Lambda_q \|\wl \|_{L^\infty(\Omp)}(\|\nabla q\|_{L^2(\Omp)}+\|\nabla \widetilde{q}\|_{L^2(\Omp)} )+2 \Lambda_q^2\| \nabla \wl \|_{L^2(\Omp)} \Big) \\
    &+2\Lambda_q^3 \| \wl \|_{L^\infty(\Omp)} \| \nabla q \|_{L^2(\Omp)}  \le CL^{\frac{1}{2}}.
\end{align*}
With the test function $\fyl$, the right hand side of \eqref{eqn:stab_integrate} equals to $\frac{1}{2}\int_{\Omega}\frac{(q-\widetilde{q})^2}{q}|\nabla \wl|^2\d x$. Therefore, by the relations \eqref{eqn:stab_test}, \eqref{eqn:stab_integrate} and the assumption $q=\tilde q$ in $\Omega\setminus\Omp$, we achieve 
\begin{align*}
    \frac{1}{2}\int_{\Omp}\frac{(q-\widetilde{q})^2}{q}|\nabla \wl|^2\d x= \int_{\Omp}  \widetilde{q} \nabla (\twl-\wl) \cdot \nabla \fyl \d x\le CL^{\frac{1}{2}} \| \twl -\wl\|_{H^1(\Omp)}.
\end{align*}
Taking summation with respect to $\ell$, we obtain
\begin{align*}
    \int_{\Omp}\frac{(q-\widetilde{q})^2}{q^2} \sum_{\ell=1}^{L}|\nabla \wl|^2\d x \le C L^{\frac{1}{2}}\sum_{\ell=1}^{L} \| \twl -\wl\|_{H^1(\Omp)} .
\end{align*}
The non-zero condition \eqref{eqn:nonzero_condition} indicates $ \sum_{\ell=1}^{L}|\nabla \wl(x)|^2\ge C_0^2 $, for all $x\in \Omp$. Hence, we conclude 
\begin{equation*}
    \|q-\tilde q\|_{L^2(\Omp)}^2 \le CC_0^{-2} L^{\frac{1}{2}}\sum_{\ell=1}^{L} \| \twl -\wl\|_{H^1(\Omp)}
\end{equation*}
Since $q = \widetilde{q}$ in $\Omega\setminus\Omp$, the proof is  completed.
\end{proof}

\begin{remark}\label{rmk:subdomain}
The proof of Theorem~\ref{thm:stability} depends on the non-zero condition \eqref{eqn:nonzero_condition} and the boundedness of   $\|\wl\|_{L^{\infty}(\Omp)}\le C\|\gl\|_{H^{\frac{1}{2}}(\partial\Omega)}\le CL^{\frac{1}{2}}$, which is satisfied under an overwhelming probability.  It is important to emphasize that the constant $C$ {and $C_0$} in \eqref{eqn:stability} are influenced by the distance between $\Omega'$ and $\partial \Omega$. { Consequently, our current stability estimates are restricted to a subdomain $\Omega' \Subset \Omega$, and we assume that the diffusivity coefficients are known in a neighborhood of $\partial \Omega$.} 
As the subdomain $\Omega'$ approaches {the} boundary of $\Omega$, controlling the regularity of the solutions and maintaining the stability of the inverse problem becomes increasingly challenging.
%In the limiting case, where $\Omega' = \Omega$ and $q = \widetilde{q}$ on $\partial \Omega$, the domain $\Omega$ and the boundary conditions $\gl$ must exhibit higher regularity to ensure that $\wl\in C^{1,\kappa}(\overline{\Om})$.
{ In the limiting case, where $\Omega' = \Omega$ and $q = \widetilde{q}$ on $\partial \Omega$, addressing the boundedness of $\|\wl \|_{L^\infty(\Omega)}$ requires higher regularity for both the domain $\Omega$ and the boundary data $\gl$. {Furthermore, establishing the non-zero condition \eqref{eqn:nonzero_condition} over the entire domain $\Omega$ requires an extension of the results in \cite{Alberti:2022}, which are based on a quantitative Runge approximation: this is expected to hold also when $\Omega'=\Omega$ (as shown in the numerical experiments in Section~\ref{sec:numer}), but no rigorous proof is available (see  \cite[Section~4.4]{bal-uhlmann-2013} for related results). }}
\end{remark}

\subsection{Finite element approximation and error estimate}\label{sec:fem}
{In this section, we introduce a numerical algorithm for the IDP and derive the reconstruction error estimation. First, we briefly state some standard results in Galerkin FEM approximation.} 
%In this section, we briefly state some standard results in Galerkin FEM approximation. 
We assume $\Omega\subset\mathbb{R}^d$ ($d=2,2$) is a bounded domain with sufficient smooth boundary $\partial\Omega$. Let $\mathcal{T}_h$ be a shape regular quasi-uniform partitions of $\Omega$ that fit the boundary exactly with a mesh size $h$. We assume that $\partial\Omp$ does not cross an element, that is,  $\Omp$ {is exactly represented as the union of a subset of elements in $\mathcal{T}_h$}. Let $V_h$ denote the conforming finite element space with piecewise polynomials of degree 1 and $\mathring{V}_{h}=V_h\cap H_0^1(\Omega)$. In particular the finite element space $V_h$ can characterized by curved element method \cite{Zlamal:1973,Zlamal:1974} when $d=2$ or isoparametric element method \cite{Ciarlet:1972,Lenoir:1986} when $d\ge 2$.

The following inverse inequality holds on the finite element space $\mathring{V}_h$ \cite[Lemma 4.5.3]{Brenner:2008}: for $0\le t\le s \le 1$ and $1\le p,q\le \infty$ we have
\begin{align}\label{eqn:inverse_ineq}
	\lVert\varphi_h\rVert_{W^{s,p}\II}\le Ch^{t-s+d/p-d/q}\lVert \varphi_h\rVert_{W^{t,q}\II},\quad \forall\varphi_h\in \mathring{V}_h.
\end{align}
Let $I_h\colon C(\overline{\Omega})\rightarrow V_h$ be the Lagrange nodal interpolation operator. Following interpolation error holds \cite[Corollary 4.4.20]{Brenner:2008}: for $s=1,2$ and $1\le p\le \infty $ (with $sp>d$ if $p>1$ and $sp\ge d $ if $p=1$)  
\begin{align}\label{eqn:error_I_h}
	\lVert v-I_h v \rVert_{L^p\II}+\lVert\nabla(v-I_h v)\rVert_{L^p\II}\le Ch^s\lVert v \rVert_{W^{s,p}\II},\quad\forall v\in W^{s,p}\II.
\end{align}
%{Let $I_h\colon L^1(\Omega)\rightarrow V_h$ be the Carstensen quasi-interpolation operator  \cite{Carstensen:1999}, defined by
%\begin{equation*}
 %I_hu:=\sum_{i=1}^{N_{\mathcal{T}_h}}\pi_i(u)\phi_i,\quad \pi_i(u):=\frac{\int_{\omega_i}u\phi_i\d x}{\int_{\omega_i}\phi_i\d x}\text{ with }\omega_i:=\mathrm{supp} \,\phi_i,
%\end{equation*}
%where $ \{\phi_i\}_{i=1}^{N_{\mathcal{T}_h}}\subset V_h$ denotes the canonical nodal basis of $V_h$. By definition, the interpolation ${I}_h$ preserves box constraint:
%\begin{equation*}
%    0<\Lambda_q^{-1}\le q\le \Lambda_q\implies %0<\Lambda_q^{-1}\le {I}_h q\le \Lambda_q.
%\end{equation*}
%Moreover, ${I}_h$ is stable in $H^1(\Omega)$ \cite{Carstensen:1999}
%\begin{equation}\label{eqn:stab_I_h}
%    \lVert   I_h v \rVert_{H^1\II} \le C \lVert v \rVert_{H^1\II},\quad\forall v\in H^1\II
%\end{equation}
%and satisfies the error estimate for $s=1,2$ and $1\le p\le \infty $  
%\begin{align}\label{eqn:error_I_h}
%	\lVert v-I_h v \rVert_{L^p\II}+h\lVert\nabla(v-I_h v)\rVert_{L^p\II}\le Ch^s\lVert v \rVert_{W^{s,p}\II},\quad\forall v\in W^{s,p}\II.
%\end{align} 
%    }
%{Let $I_h\colon H^1(\Omega)\rightarrow V_h$ be the Scott-Zhang interpolation operator \cite{ScottZhang:1990}. Following interpolation error holds \cite[Theorem 4.8.12]{Brenner:2008}: for $s=1,2$ and $1\le p\le \infty $  
%\begin{align}\label{eqn:error_I_h}
%	\lVert v-I_h v \rVert_{L^p\II}+h\lVert\nabla(v-I_h v)\rVert_{L^p\II}\le Ch^s\lVert v \rVert_{W^{s,p}\II},\quad\forall v\in W^{s,p}\II.
%\end{align} }
We use $I_h^{\partial}$ to denote the Lagrange interpolation operator on the boundary. We define the $L^2\II$-projection $P_h\colon L^2\II\rightarrow \mathring{V}_h$ by
\begin{align*}
	(P_h v,\varphi_h)=( v,\varphi_h),\quad\forall\varphi_h\in \mathring{V}_h.
\end{align*}
The operator $P_h$ satisfies the following error estimates \cite[p. 32]{Thomee:2007}: for every $s\in [1,2]$ we have
\begin{align}\label{eqn:error_P_h}
	\lVert v-P_h v\rVert_{L^2\II}+\lVert\nabla(v-P_h v)\rVert_{L^2\II}\le Ch^s\lVert v\rVert_{H^s\II},\quad\forall v \in H^s\II\cap H_0^1\II.
\end{align}

Now, we present the reconstruction algorithm. Slightly differently from the stability analysis, we aim to reconstruct the diffusion coefficient in the whole domain $\Omega$ using the measurement in the entire domain. Throughout this section, we let $z_{\delta}^{(\ell)}$ denote the practical noisy 
observations corresponding to $\wl(q^{\dag})$ with noise level $\delta$, i.e.  
\begin{equation}\label{eqn:noise}
     \| \wl(q^{\dag})-z_{\delta}^{(\ell)}\|_{L^2(\Om)}\le \delta, \quad\forall \ell=1,\dots,L.
\end{equation}  
The reconstruction is based on standard regularized least-squares  with further discretization using finite element methods. 
More precisely, the minimization problem is
\begin{equation}\label{eqn:cts_functional}
		\min_{q\in \Uad_q} J_{\al}(q)=\frac{1}{2}\sum_{\ell=1}^{L}\| \wl(q)-z_{\delta}^{(\ell)} \|_{L^2(\Om)}^2+\frac{\al L}{2}\|\nabla q\|_{L^2\II}^2,
\end{equation}
where $\al>0$ is the regularization parameter,
and $\wl(q)\in H^1\II$ is the weak solution of
\begin{equation}\label{eqn:cts_constraint}
 	\left\{
 		\begin{aligned}	
 			-\nabla\cdot(q \nabla \wl)& =0, && \mbox{in }\Om, \\
 			\wl & =\gl, && \mbox{on } \partial\Om.
 		\end{aligned}
 	\right.
\end{equation}
We formulate the finite element approximation of {the} problem \eqref{eqn:cts_functional}-\eqref{eqn:cts_constraint}:
\begin{equation}\label{eqn:dis_functional}
	\min_{q_h\in \Uad_{q,h}} J_{\al,h}(q_h)=\frac{1}{2}\sum_{\ell=1}^{L}\| \wl_{h}(q_h)-z_{\delta}^{(\ell)} \|_{L^2(\Om)}^2+\frac{\al L}{2}\|\nabla q_h \|_{L^2\II}^2,
\end{equation}
where $\wl_{h}(q_h)\in V_h$ is the weak solution of
\begin{equation}\label{eqn:dis_constraint}
	\left\{
	\begin{aligned}	
		(q_h \nabla \wl_{h},\nabla v_h)& =0, && \forall v_h\in \mathring{V}_h, \\
		\wl_{h} & =I_h^{\partial} \gl, && \mbox{on } \partial\Omega.
	\end{aligned}
	\right.
\end{equation}
Here, the admissible set is defined as
\begin{equation}\label{eq:Aqh}
      \Uad_{q, h}=\{q_h\in V_h \mbox{ : }0<\Lambda_q^{-1}\le q_h\le \Lambda_q\mbox{ a.e. in }\Om   \mbox{, }  q_h=I_h q^{\dag}\mbox{  on } \partial\Omega \}.
\end{equation}
The discrete problem \eqref{eqn:dis_functional}-\eqref{eqn:dis_constraint} is well-posed: for every $\alpha>0$, there exists at least one global minimizer $q_h^*$ and  it depends continuously on the data perturbation. {Moreover, as the noise level $\delta\to0^+$, the sequence of minimizers converges to the exact coefficient $q^\dagger$ in $H^1(\Omega)$, provided the regularization parameter $\alpha$ is chosen properly with respect to $\delta$ \cite{EnglKunischNeubauer:1989,ItoJin:2015}. Although the non-convex nature of the functional means that a global minimizer is not necessarily unique, our analysis holds valid for any such global minimizer.} 
The main objective in this section is to bound the approximation error $\|q^{\dag}-q_h^*\|_{L^2\II} $. The strategy is based upon the stability analysis in the preceding section. 
Furthermore, we need the following higher regularity assumption on the exact diffusivity coefficient and boundary data.
 
\begin{assum}\label{assum:high_regularity}
Let $\Omega\subset \mathbb{R}^d$ ($d=2,3$) be a bounded domain with  $C^{1,1}$ boundary $\partial\Omega$. Assume that the exact diffusivity coefficient $q^{\dag}\in W^{2,p}\II\cap \mathcal{A}_q$ with $p>d$.  Assume the boundary data $\gl$ (with $\ell=1,\dots,L$)   are taken as independent and identically distributed satisfying the expansion \eqref{eqn:reps_g}, where $\{e_k\}_{k=1}^{\infty}$  is the orthonormal basis of $H^{\frac12}(\partial\Omega)$ consisting of the eigenfunctions of the Laplace-Beltrami operator on $\partial\Omega$  and $a_k^{(\ell)}\sim N(0,  \theta_k^2)$, with $\theta_k^2\lesssim \frac{1}{k^\beta}$ with $\beta>\frac{3}{d-1}+1$.
\end{assum}

\begin{remark}\label{rmk:high_regularity}
Assumption~\ref{assum:high_regularity} requires higher regularity for the domain $\Omega$ as well as the parameter $q^\dag$ and $\gl$ to ensure that the finite element approximation achieves an optimal convergence rate. Indeed, under the regularity assumption, Sobolev embedding theory and elliptic regularity theory (\cite[Theorem 7.2]{Giaquinta:2013} and  \cite[Theorem 8.12]{Gilbarg:1977}) implies that the  solution satisfies  $\wl(q^\dag)\in    H^2\II\cap W^{1,\infty}\II$ when $d=2$,  $\wl(q^\dag)\in    H^2\II\cap W^{1,p}\II$ for all $2<p<\infty$ when $d=3$.  
Under Assumption~\ref{assum:high_regularity},  with a probability greater than 
\begin{equation*} 
1-  L^d \exp\left(-C_1 L\right)- L\exp\left(-C_2  M\right),
\end{equation*}
the  non-zero condition \eqref{eqn:nonzero_condition} holds and the random boundary data has the upper bound  
\begin{equation}\label{eqn:bounded_bdry_H2}
    \max_{\ell=1,\dots,L} \|g^{(\ell)}\|_{H^2(\partial\Omega)}\le L^{\frac{1}{2}},
\end{equation} 
where the positive constants $C_0$, $C_1$ and  $C_2$  depend only on $s$, $\Om$, $\Omp$, $\{\theta_k\}$, $\Lambda_q$ and  $\|q^{\dagger}\|_{C^{0,1}(\overline{\Omega})}$. The nonzero condition is a direct consequence of Proposition \ref{prop:nonzero_condition}. 
It suffices to investigate the upper bound of $\|\gl\|_{H^2(\partial\Omega)}$. Note that the Laplace--Beltrami operator $-\Delta$ on $\partial\Omega$ admits a positive sequence $\{\lambda_k\}_{k=1}^{\infty}$ of eigenvalues and the corresponding eigenfunctions $\{\varphi_k\}_{k=1}^{\infty}$ form an orthonormal basis of $L^2(\partial\Omega)$. Here we use the equivalent norm in space $H^{s}(\partial\Omega) $, with $s > 0$, defined by  \cite[Chapter 1, Remark 7.6]{Lions:1972} 
\begin{equation}\label{eq:lions-magenes}
    \|g\|_{H^{s}(\partial\Omega)}^2=   \sum_{k=1}^{\infty}(1+\lambda_k)^{s}(g,\varphi_k)_{\partial\Omega}^2 .
\end{equation}
Therefore,   $e_k=(1+\lambda_k)^{-1/4}\varphi_k$.  Thus, recalling \eqref{eqn:reps_g}, we have $g^{(l)}=\sum_{k=1}^M a_k^{(l)} (1+\lambda_k)^{-1/4} \varphi_k$. Hence, by \eqref{eq:lions-magenes}
 we obtain
 \[ \|g^{(l)}\|_{H^2(\partial\Omega)}^2 = \sum_{k=1}^M (1+\lambda_k)^2 (g^{(l)},\varphi_k)_{\partial\Omega}^2
= \sum_{k=1}^M (1+\lambda_k)^{\frac{3}{2}} (a_k^{(l)})^2.
\]
 
By  the asymptotic behavior of eigenvalues $\lambda_k\sim k^{\frac{2}{d-1}}$ \cite[Theorem 1.1]{Safarov:1997}, the moment generating function of $\|\gl\|_{H^2(\partial\Omega)}$ satisfies for all $\lambda\in\mathbb{R}$: 
\begin{align*}
    \Ebb\exp\left(\lambda^2\| \gl\|_{H^2(\partial\Om)}^2\right) 
    =& \Ebb\exp\left(\lambda^2\sum_{k=1}^M (1+\lambda_k)^{\frac{3}{2}} (a_k^{(l)})^2\right) \\
    = & \exp\left(\lambda^2\sum_{k=1}^M (1+\lambda_k)^{\frac{3}{2}} \theta_k^2\right) \lesssim \exp\left(\lambda^2\sum_{k=1}^M k^{\frac{3}{d-1}-\beta} \right).
\end{align*}
Then, since $\frac{3}{d-1}-\beta<-1$,  by \cite[Proposition 2.5.2]{Vershynin:2018}, with probability greater than $1-L\exp(-C_1L)$, we have 
\begin{equation*}
\max_{\ell=1,\dots,L}\|\gl\|_{H^2(\partial\Omega)}\le L^{\frac{1}{2}}.
\end{equation*} 
\end{remark}

We have the following $L^2(\Omega)$ error estimate for $w_h(q^{\dag})-{w_h}(I_h q^{\dag})$.
\begin{lemma}\label{lem:error_coeff}
Let Assumption~\ref{assum:high_regularity} hold and the boundary data satisfy $\|g\|_{H^2(\partial\Omega)}\le L^{\frac{1}{2}}$.  We denote the solutions of equation \eqref{eqn:dis_constraint} with coefficients $q^{\dag}$ and $I_h q^{\dag}$ by $w_h(q^{\dag})$ and $w_h(I_h q^{\dag})$, respectively. Then 
\begin{equation*}
    \|w_h(q^{\dag})-{w_h}(I_h q^{\dag}) \|_{L^2(\Om)}\le C h^2 L^{\frac{1}{2}},
\end{equation*}
where C is a positive constant depending only on $\Omega$ and $q^\dagger$.
\end{lemma} 
\begin{proof}
With an abuse of notation, several positive constants depending only on $\Omega$ and $q^\dagger$ will be denoted by the same letter $C$.
We start with the estimate in energy norm. By subtracting the weak formulations of ${w_h}(q^{\dag})$ and ${w_h}(I_h q^{\dag})$, we derive
\begin{align*}
\big(I_h q^{\dag} (\nabla {w_h}(I_h q^{\dag})-\nabla{w_h}(q^{\dag}) ) ,\nabla {v_h}\big)=\big( (q^{\dag}-I_h q^{\dag})\nabla {w_h}(q^{\dag}),\nabla {v_h} \big),\quad \text{for all}  ~~{v_h}\in \mathring{V}_h.
\end{align*}
Select the test function ${v_h}= {w_h}(I_h q^{\dag})-{w_h}( q^{\dag})$. Note that it belongs to $\mathring{V}_h$ since ${u_h}(I_h q^{\dag})$ and ${u_h}(q^{\dag})$ share the same boundary value.
Using the box constraint on $q^{\dag}$ and the Cauchy--Schwarz inequality, we obtain
\begin{align*}
 {\|\nabla {w_h}(I_h q^{\dag})-\nabla{w_h}(q^{\dag})\|_{L^2(\Omega)}^2 
\le C\|q^{\dag}-I_h q^{\dag}\|_{L^\infty(\Omega)} \|\nabla {w_h}(q^{\dag})\|_{L^2(\Omega)}\|\nabla {w_h}(I_h q^{\dag})-\nabla{w_h}(q^{\dag})\|_{L^2(\Omega)}.}
\end{align*}
Then the approximation estimate \eqref{eqn:error_I_h}  implies 
\begin{equation}\label{eqn:grad_u_h(I_hq)-u(q)}
\|\nabla {w_h}(I_h q^{\dag})-\nabla{w_h}(q^{\dag})\|_{L^2(\Omega)}\le C h\|\nabla {w_h}(q^{\dag})\|_{L^2(\Omega)}\le ChL^{\frac{1}{2}}.
\end{equation}
Next, we apply the duality argument to get the estimate in $L^2(\Omega)$ norm. % $\|{u_h}(I_h q^{\dag})-{u_h}(q^{\dag})\|_{L^2\II}$.
Let $\psi$ satisfy
\begin{equation*}
\begin{aligned}
 -\nabla\cdot (q^{\dag}\nabla \psi)={w_h}(I_h q^{\dag})-{w_h}(q^{\dag})~~\mbox{in }\Omega, \quad\mbox{with}~~
  \psi =0~~ \mbox{on } \partial\Omega.
\end{aligned}
\end{equation*}
Then we have
\begin{equation*}
\begin{aligned}
\|{w_h}(I_h q^{\dag})-{w_h}(q^{\dag})\|_{L^2(\Omega)}^2
         =&\big(-\nabla\cdot (q^{\dag}\nabla \psi),{w_h}(I_h q^{\dag})-{w_h}(q^{\dag}) \big)\\
         =&\big( q^{\dag}\nabla \psi,\nabla ({w_h}(I_h q^{\dag})-{w_h}(q^{\dag})) \big)\\
         =&\big( (q^{\dag}-I_h q^{\dag})\nabla \psi,\nabla ({w_h}(I_h q^{\dag})-{w_h}(q^{\dag})) \big)\\
         &+ \big( I_h q^{\dag}\nabla (\psi- P_h \psi),\nabla ({w_h}(I_h q^{\dag})-{w_h}(q^{\dag})) \big)\\
         &+\big( (q^{\dag}-I_h q^{\dag} )\nabla P_h \psi,\nabla {w_h}(q^{\dag}) \big),
\end{aligned}
\end{equation*}
where we used the weak formulation of ${w_h}(q^{\dag})$ and  ${w_h}(I_h q^{\dag})$ in the last equality. Therefore, by H\"older inequality, error estimate  \eqref{eqn:error_I_h}, \eqref{eqn:error_P_h}  and \eqref{eqn:grad_u_h(I_hq)-u(q)} yield  that
    \begin{align*}
        \|{w_h}(I_h q^{\dag})-{w_h}(q^{\dag})\|_{L^2(\Omega)}^2
        \le & Ch^2 \| q^{\dag}\|_{W^{1,\infty}(\Omega)}\|\nabla\psi\|_{L^2(\Omega)}\|\nabla {w_h}(q^{\dag})\|_{L^2(\Omega)}\\
        &+Ch^2\|I_h q^{\dag}\|_{L^{\infty}(\Omega)}\|\psi\|_{H^2(\Omega)}\|\nabla {w_h}(q^{\dag})\|_{L^2(\Omega)}\\
        & +Ch^2\|q^{\dag}\|_{W^{2,p}(\Omega)}\|\nabla P_h \psi\|_{L^q(\Omega)}\|\nabla {w_h}(q^{\dag})\|_{L^{2}(\Omega)}.
    \end{align*}
	Here $\frac{1}{p}+\frac{1}{q}+\frac{1}{2}=1$ and, by Assumption~\ref{assum:high_regularity}, $q=\frac{2p}{p-2}<\frac{2d}{d-2}$. Thus the stability of the $L^2(\Omega) $ projection (see \cite[Theorem 4]{Crouzeix:1987} and \cite[Lemma 2.1]{Bakaev:2001}) and the Sobolev embedding imply $\|\nabla P_h\psi\|_{L^q(\Omega)}\le C \|\nabla  \psi\|_{L^q(\Omega)} \le C \|   \psi\|_{H^2(\Omega)}$. By using standard  elliptic regularity estimates, according to which $\|\psi \|_{H^2(\Omega)}\le C \|{w_h}(I_h q^{\dag})-{w_h}(q^{\dag})\|_{L^2(\Omega)}$, we obtain
    \begin{align*}
         \|{w_h}(I_h q^{\dag})-{w_h}(q^{\dag})\|_{L^2(\Omega)}\le  C h^2\|\nabla {w_h}(q^{\dag})\|_{L^2(\Omega)} \le Ch^2 L^{\frac{1}{2}}.
    \end{align*}
    This completes the proof of the lemma.
\end{proof}

\begin{corollary}\label{cor:error_discrete}
Let Assumption~\ref{assum:high_regularity} hold and the boundary data satisfy $\|g\|_{H^2(\partial\Omega)}\le L^{\frac{1}{2}}$.  Let $w( q^{\dag})$ be the solution of equation \eqref{eqn:cts_constraint} and $w_{h}(I_h q^{\dag})$ be the solution of equation \eqref{eqn:dis_constraint}. Then
\begin{equation*}
    \|w_{h}(I_h q^{\dag})-w( q^{\dag}) \|_{L^2(\Om)}\le Ch^2 L^{\frac{1}{2}},
\end{equation*}
where $C$ is a positive constant depending only on $\Omega$ and $q^\dagger$. 
\end{corollary}
\begin{proof}
We use the following splitting
    \begin{align*}
        \|w_{h}(I_h q^{\dag})-w( q^{\dag}) \|_{L^2(\Om)}\le& \|w_{h}(I_h q^{\dag})-w_{h}( q^{\dag}) \|_{L^2(\Om)}+\|w_{h}( q^{\dag})-w( q^{\dag}) \|_{L^2(\Om)}.
    \end{align*}
    For the first term, we apply Lemma~\ref{lem:error_coeff} and obtain
    $ \|w_{h}(I_h q^{\dag})-w_{h}( q^{\dag}) \|_{L^2\II}   \le Ch^2 L^{\frac{1}{2}}$. 
   The second term can be estimated by utilizing the standard duality argument with the interpolation estimate $\|g-I_h^\partial g\|_{L^2(\partial\Omega)}\le ch^2L^{\frac{1}{2}}$.%, which holds by the construction of the interpolation operator and the property of nonlinear map $\Phi_T$.
\end{proof}

The next lemma gives an \textit{a priori} estimate.  {It measures the accuracy of the discrete forward state and provides the necessary \textit{a priori} control over the coefficient $q_h^*$.}
\begin{lemma}\label{lem:apriori}
Let Assumption~\ref{assum:high_regularity} hold and boundary data satisfy $\|g^{(\ell)}\|_{H^2(\partial\Omega)}\le L^{\frac{1}{2}}$, $\ell=1,\dots,L$. Let $q_h^*\in \Uad_{q,h}$  be a minimizer of problem \eqref{eqn:dis_functional}-\eqref{eqn:dis_constraint}. Then we have
\begin{equation*}
\sum_{\ell=1}^{L}\| \wl_{h}(q_h^*)-\wl(q^{\dag}) \|_{L^2(\Om)}+L\al^{\frac{1}{2}} \|\nabla q_h^* \|_{L^2(\Om)}\le C  L  (h^2 L^{\frac{1}{2}}+\delta+\al^{\frac{1}{2}}),
\end{equation*}
where $C$ is a positive constant depending only on $\Omega$ and $q^\dagger$.
\end{lemma}
\begin{proof}
With an abuse of notation, several positive constants depending only on $\Omega$ and $q^\dagger$ will be denoted by the same letter $C$.
    Since $q_h^*$ is a minimizer of $J_{\alpha,h}$, we have $J_{\al,h}(q_h^*)\le J_{\al,h}(I_h q^{\dag})$. As a result,
    \begin{align*}
        &\frac{1}{2}\sum_{\ell=1}^{L}\| \wl_{h}(q_h^*)-z_{\delta}^{(\ell)} \|_{L^2(\Om)}^2+\frac{\al L}{2}\|\nabla q_h^* \|_{L^2(\Om)}^2\\
        &\le \frac{1}{2}\sum_{\ell=1}^{L}\| \wl_{h}(I_h q^{\dag})-z_{\delta}^{(\ell)} \|_{L^2(\Om)}^2+\frac{\al L}{2}\|\nabla I_h q^{\dag} \|_{L^2(\Om)}^2\\
        &\le \sum_{\ell=1}^{L} \Big(\| \wl_{h}(I_h q^{\dag})-\wl(q^{\dag}) \|_{L^2(\Om)}^2+\| \ul(q^{\dag})-z_{\delta}^{(\ell)} \|_{L^2(\Om)}^2\Big)+\frac{\al L}{2}\|\nabla I_h q^{\dag} \|_{L^2(\Om)}^2.
    \end{align*}
    %{By the stability of Scott-Zhang interpolation \cite[Corollary 4.8.15]{Brenner:2008}, we have $\|\nabla I_h q^{\dag} \|_{L^2(\Om)} \le c \|\nabla q^\dagger\|_{L^2(\Omega)}\le c$. } 
    By the interpolation property \eqref{eqn:error_I_h} and regularity of $q^\dag$, the term $\|\nabla I_h q^{\dag} \|_{L^2(\Om)} $ can be bounded by
    \begin{align*}
        \|\nabla I_h q^\dag\|_{L^2(\Omega)}\le & \|\nabla I_h q^\dag-\nabla q^\dag \|_{L^2(\Omega)}+\|\nabla   q^\dag\|_{L^2(\Omega)}\\
        \le & Ch\|q^\dag\|_{H^2(\Omega)}+\|q^\dag \|_{H^1(\Omega)}\le C.
    \end{align*} 
   This, together with  Corollary~\ref{cor:error_discrete} and the bound for the noise level  in \eqref{eqn:noise}, implies that
    \begin{align*}
        \frac{1}{2}\sum_{\ell=1}^{L}\| \wl_{h}(q_h^*)-z_{\delta}^{(\ell)} \|_{L^2(\Om )}^2+\frac{\al L}{2}\|\nabla q_h^* \|_{L^2(\Om)}^2\le CL(h^4L+\delta^2+\alpha).
    \end{align*}
Hence, we derive $ \al^{\frac{1}{2}} \|\nabla q_h^* \|_{L^2\II}\le  C(h^{2}L^{\frac{1}{2}}+\delta+\alpha^{\frac{1}{2}} )$.
Then the triangle inequality and the Cauchy-Schwarz inequality lead to
    \begin{align*}
\sum_{\ell=1}^{L}\| \wl_{h}(q_h^*)-\wl(q^{\dag}) \|_{L^2(\Om )}&\le \sum_{\ell=1}^{L}\| \wl_{h}(q_h^*)-z_{\delta}^{(\ell)}\|_{L^2(\Om )}+\sum_{\ell=1}^{L}\| z_{\delta}^{(\ell)}-\wl(q^{\dag}) \|_{L^2(\Om )}\\
        &\le L^{\frac{1}{2}}\Big(\sum_{\ell=1}^{L}\| \wl_{h}(q_h^*)-z_{\delta}^{(\ell)}\|_{L^2(\Om)}^2\Big)^{\frac{1}{2}}+ L\delta\le C L  (h^2 L^{\frac{1}{2}}+\delta+\al^{\frac{1}{2}}).
    \end{align*}
\end{proof}

Next, we state our main theorem, estimating the error between the exact diffusivity coefficient $q^\dag$ and the numerical reconstruction $q_h^*$. {The error analysis below employs the conditional stability argument implicitly within the discrete setting: rather than applying Theorem \ref{thm:stability} on the discrete minimizer $q_h^*$ directly, we reconstruct the weighted energy estimate with the special test function at the finite element level, using Lemma \ref{lem:apriori} to establish the required estimates.}

\begin{theorem}\label{thm:error_estimate}
Suppose the exact diffusivity coefficient $q^{\dag}$ and the random boundary illuminations $g^{(\ell)}$ (with $\ell=1,\dots,L$) satisfy Assumption~\ref{assum:high_regularity}. Let $q_h^*\in \Uad_{q,h}$ be a minimizer of problem \eqref{eqn:dis_functional}-\eqref{eqn:dis_constraint}. Set $\xi=h^2 L^{\frac{1}{2}}+\delta+\al^{\frac{1}{2}}$. Then, with probability greater than \eqref{eqn:probability_nonzero}, 
    we have 
    \begin{equation*} 
        \|q^{\dag}- q_h^*\|_{L^2(\Omp)}^2\le C C_0^{-2}  L^2(1+\alpha^{-\frac{1}{2}}\xi)\left(h+ h^{1-\epsilon}(1+\alpha^{-\frac{1}{2}}\xi)+\min\left(1,h + h^{-1}L^{-\frac{1}{2}}\xi \right)  \right),
    \end{equation*}
where $\epsilon = 0$ when $d = 2$, and $\epsilon > 0$ is arbitrarily small when $d = 3$. Here, $C$ is a positive constant depending only on $\epsilon$, $\Omega$, and $q^{\dagger}$, while $C_0$ is defined in \eqref{eqn:nonzero_condition}.
\end{theorem}
\begin{proof}
With an abuse of notation, several positive constants depending only on $\Omega$ and $q^\dagger$ will be denoted by the same letter $C$. Let $\wl=\wl(q^\dag)$ be the solution to \eqref{eqn:cts_constraint} with boundary value $\gl$. For a test function $\fyl\in H^1_0(\Omega)$, we multiply  both sides of \eqref{eqn:cts_constraint} by $ (I_h q^\dag-q_h^*)\fyl / q^\dag$, and apply integration by parts:
\begin{equation*}
       0=\big( q^\dag\nabla \wl,\nabla \frac{(I_h q^\dag-q_h^*)\fyl }{q^\dag} \big)=  \big( q^\dag \fyl\nabla \wl,\nabla \frac{(I_h q^\dag-q_h^*)}{q^\dag} \big)+\big( (I_h q^\dag-q_h^*) \nabla \wl,\nabla \fyl\big).
\end{equation*}
Thus, we obtain
\begin{equation}\label{eq:test-intermediate}
    \big( (I_h q^\dag-q_h^*) \nabla \wl,\nabla \fyl\big)=\frac{1}{2}\big( (I_h q^\dag-q_h^*) \nabla \wl,\nabla \fyl\big)-\frac{1}{2} \big( q^\dag \fyl\nabla \wl,\nabla \frac{(I_h q^\dag-q_h^*)}{q^\dag} \big).
\end{equation}
Set the test function $\fyl=(I_h q^\dag-q_h^*)\ul /q^\dag$. We first verify $\fyl\in H_0^1(\Omega)$. Since $q_h^*\in \mathcal{A}_{q,h}$, $\fyl$ vanishes on $ \partial\Omega$. Recall that, under the current assumptions, we have $\|g^{(\ell)}\|_{H^2(\partial\Omega)}\le L^{\frac{1}{2}}$ for every $\ell=1,\dots,L$, cf. Remark~\ref{rmk:high_regularity}.  By the regularity of $q^\dag$ and $\ul$, and  in view of Lemma~\ref{lem:apriori}, we conclude that $\fyl\in H_0^1(\Omega)$, with 
\begin{align*}
    \|\fyl\|_{L^2(\Omega)}= \|(I_h q^\dag-q_h^*)\wl/q^\dag\|_{L^2(\Omega)}\le 2\Lambda_q^2\|\wl \|_{L^2(\Omega)}\le CL^{\frac{1}{2}}
\end{align*}
and 
\begin{equation}\label{eq:nabla-phi}
\begin{aligned}
    &\|\nabla \fyl\|_{L^2(\Omega)}=\left\|\frac{q^\dag\nabla[(I_h q^\dag-q_h^*)\wl]-(I_h q^\dag-q_h^*)\wl \nabla q^\dag }{(q^\dag)^2}\right\|_{L^2(\Om)}\\
    &\le  \Lambda_q^3 \|\wl \|_{L^\infty(\Om)}(\|\nabla I_h q^\dag\|_{L^2(\Om)}+\|\nabla q_h^*\|_{L^2(\Om)} ) \\
    &\;+\Lambda_q^2\left(2 \Lambda_q^2\| \nabla \wl \|_{L^2(\Om)}+ 2\Lambda_q \| \wl \|_{L^\infty(\Om)} \| \nabla q^\dag \|_{L^2(\Om)} \right) \\
    &\le  CL^{\frac{1}{2}}(1+\|\nabla q_h^*\|_{L^2(\Om)})\le CL^{\frac{1}{2}}(1+\alpha^{-\frac{1}{2}}\xi).
\end{aligned}
\end{equation}
With this test function $\fyl$, by direct computation, we can further write the left hand side of \eqref{eq:test-intermediate} as
\begin{equation}\label{eqn:error_weak_form}
    \big( (I_h q^\dag-q_h^*) \nabla \wl,\nabla \fyl\big)=\frac{1}{2} \int_{\Omega}\frac{ (I_h q^\dag-q_h^*)^2}{q^\dag}|\nabla \wl|^2 \d x.
\end{equation}
On the other hand, by the weak formulation of \eqref{eqn:cts_constraint} and \eqref{eqn:dis_constraint}, we have 
\begin{align*}
    \big( (I_h q^\dag-q_h^*) \nabla \wl,\nabla \fyl\big) &=\big( (I_h q^\dag-q^\dag) \nabla \wl,\nabla \fyl\big)+\big( (  q^\dag-q_h^*) \nabla \wl,\nabla \fyl\big)\\
    & =\big( (I_h q^\dag-q^\dag) \nabla \wl,\nabla \fyl\big)
    +\big( (  q^\dag-q_h^*) \nabla \wl,\nabla (\fyl-P_h\fyl) \big)\\
    &\qquad\qquad\;\quad\qquad \qquad\qquad \quad +\big( q_h^* \nabla (w_h^{(\ell)}(q_h^*)-\wl),\nabla P_h\fyl\big)\\
     & =\mathrm{I}_1^{(\ell)}+\mathrm{I}_2^{(\ell)}+\mathrm{I}_3^{(\ell)}.
\end{align*}
For $\mathrm{I}_1^{(\ell)}$, the interpolation error \eqref{eqn:error_I_h} and the estimate \eqref{eq:nabla-phi} yield that
\begin{equation*}
    |\mathrm{I}_1^{(\ell)} |\le C \|I_h q^\dag-q^\dag\|_{L^{\infty}(\Omega)}\|\nabla \wl\|_{L^2(\Omega)}\|\nabla\fyl\|_{L^2(\Omega)}\le ChL (1+\alpha^{-\frac{1}{2}}\xi).
\end{equation*}
Now, we consider $\mathrm{I}_2^{(\ell)}$. %We first note that $P_h\fyl$ is supported in $\Omp$, since $\fyl=0$ on $\Omega\setminus\Omp$. 
Applying integration by parts, the regularity of $q^\dag$ and $\wl$, the inverse inequality \eqref{eqn:inverse_ineq}, the projection error \eqref{eqn:error_P_h}  and estimate \eqref{eq:nabla-phi} imply that
\begin{align*}
    |\mathrm{I}_2^{(\ell)} |
    &=|\big( \nabla\cdot((  q^\dag-q_h^*) \nabla \wl),\fyl-P_h\fyl \big)|\\
    &\le \bigl( \|\nabla(q^\dag-q_h^*)\|_{L^q\II}  \|\nabla \wl\|_{L^p\II} + \|q^\dag-q_h^*\|_{L^\infty\II}  \|\Delta\wl\|_{L^2\II} \bigr)\|\fyl-P_h\fyl\|_{L^2\II}\\
    &\le Ch\left(L^{\frac{1}{2}} + L^{\frac{1}{2}}h^{d/q-d/2} \|\nabla q_h^* \|_{L^2\II} ) \right)\|\fyl\|_{H^1\II}\\
    &\le Ch^{1+d/q-d/2}L(1+\alpha^{-\frac12}\xi)^2=Ch^{1-\epsilon}L(1+\alpha^{-\frac12}\xi)^2 .
\end{align*}
Here we use the regularity results for $\wl$, as stated in Remark \ref{rmk:high_regularity}, and take $\frac{1}{p}+\frac{1}{q}+\frac{1}{2}=1$ with $p=\infty$ when $d=2$,  $p=\frac{d}{\epsilon}$ when $d=3$.  To estimate $\mathrm{I}_3^{(\ell)} $,   by the inverse inequality \eqref{eqn:inverse_ineq} and the projection error  \eqref{eqn:error_P_h}, we first derive that
    \begin{align*}
        \| \nabla \wl -\nabla \wl_{h}(q_h^*) \|_{L^2(\Om)} 
        \le&   \| \nabla \wl -\nabla P_h \wl \|_{L^2(\Om)} +  \| \nabla P_h \wl -\nabla \wl_{h}(q_h^*) \|_{L^2(\Om)} \\
        \le& C \big( h\| \wl\|_{H^2(\Om)}+h^{-1}\|P_h \wl - \wl_{h}(q_h^*)\|_{L^2(\Om)}\big)\\
        \le& C \big( hL^{\frac{1}{2}}+h^{-1}\| \wl -  \wl_{h}(q_h^*)\|_{L^2(\Om)}\big).
    \end{align*}
    There obviously holds that $ \| \nabla \wl -\nabla \wl_{h}(q_h^*) \|_{L^2\II}\le CL^{\frac{1}{2}}$.
    Therefore, by using these  two inequalities, \eqref{eq:nabla-phi}  and Lemma~\ref{lem:apriori}, we obtain
    \begin{align*}
       \sum_{\ell=1}^{L} |\mathrm{I}_3^{(\ell)}|&\le  \sum_{\ell=1}^{L} \|q_h^* \|_{L^{\infty}\II}  \| \nabla \wl_{h}(q_h^*)-\nabla\wl   \|_{L^2\II}  \|\nabla P_h \fyl\|_{L^2\II}\\
        &\le CL^{\frac{1}{2}}(1+\alpha^{-\frac{1}{2}}\xi )\sum_{\ell=1}^{L}  \| \nabla \wl_{h}(q_h^*)-\nabla\wl   \|_{L^2\II}  \\
        &\le CL^{\frac{1}{2}}(1+\alpha^{-\frac{1}{2}}\xi ) \min\Big(L^{\frac{3}{2}},L^{\frac{3}{2}}h+h^{-1} \sum_{\ell=1}^{L}\| \wl_h({q_h^*}) -  \wl \|_{L^2\II} \Big)\\
        &\le  CL^2 (1+\alpha^{-\frac{1}{2}}\xi )\min\Big(1,h + h^{-1}L^{-\frac{1}{2}}\xi \Big).
    \end{align*}
Taking summation with respect to $\ell=1,\dots,L$ in \eqref{eqn:error_weak_form},  the estimates of $\mathrm{I}_1^{(\ell)},\mathrm{I}_2^{(\ell)},\mathrm{I}_3^{(\ell)}$ yield that
\begin{equation*}
\begin{split}
    \frac{1}{2} \int_{\Omega} & \frac{ (I_h q^\dag-q_h^*)^2}{q^\dag} \sum_{\ell=1}^{L}|\nabla \wl|^2 \d x\\
    &\qquad \le CL^2(1+\alpha^{-\frac{1}{2}}\xi)\left(h+ h^{1-\epsilon}(1+\alpha^{-\frac{1}{2}}\xi)+\min\left(1,h + h^{-1}L^{-\frac{1}{2}}\xi \right)  \right).
\end{split}
\end{equation*}
Applying the interpolation error bound $\|q^\dag-I_h q^\dag\|_{L^2(\Omega)}\le Ch^2\|q^\dag\|_{H^2(\Omega)}$ (see \eqref{eqn:error_I_h}), we arrive at the weighted estimate
\begin{equation*}
    \begin{split}
        \frac{1}{2} \int_{\Omega} & \frac{ (  q^\dag-q_h^*)^2}{q^\dag} \sum_{\ell=1}^{L}|\nabla \wl|^2 \d x\\
        &\qquad \le CL^2h^4+ CL^2(1+\alpha^{-\frac{1}{2}}\xi)\left(h+ h^{1-\epsilon}(1+\alpha^{-\frac{1}{2}}\xi)+\min\left(1,h + h^{-1}L^{-\frac{1}{2}}\xi \right)  \right).
    \end{split}
\end{equation*}
%By Proposition~\ref{prop:nonzero_condition},
By Remark~\ref{rmk:high_regularity}, we have the non-zero condition \eqref{eqn:nonzero_condition}: 
$$ \sum_{\ell=1}^{L}|\nabla \wl(x)|^2\ge  C_0^2,\quad  \text{for all} ~x\in \Omp.$$ 
%Since $q_h^*\in \mathcal{A}_{q,h}$, $q_h^*=I_h q^\dag$ on $\Omega\setminus\Omp$. 
Hence, we conclude 
\begin{equation*}
    \|  q^\dag-q_h^*\|_{L^2(\Omp)}^2 \le %CC_0^{-2}L^2 h^4+
    CC_0^{-2}L^2(1+\alpha^{-\frac{1}{2}}\xi)\left(h+ h^{1-\epsilon}(1+\alpha^{-\frac{1}{2}}\xi)+\min\left(1,h + h^{-1}L^{-\frac{1}{2}}\xi \right)  \right).
\end{equation*}
This completes the proof.
\end{proof}

\begin{remark}\label{rmk:error_q_delta} 
Theorem~\ref{thm:error_estimate} provides a guideline for the \textit{a priori} choice of the algorithmic parameters $h$ and $\alpha$, in relation to $\delta$. The choice $h^2L^{\frac{1}{2}}\sim \delta $ and $\alpha\sim \delta^2$ yields a convergence rate 
$$  \|q^{\dagger}- q_h^* \|_{L^2(\Om)}\le C L^{\frac{7}{8}} \delta^{\frac{1}{4}-\epsilon},$$
with $\epsilon=0$ for $d=2$, $\epsilon>0$ arbitrary small for $d=3$. This  rate is consistent with the stability in Theorem \ref{thm:stability}, that shows 
\begin{align*}
    \|q^\dagger-q\|_{L^2\II}\le C C_0^{-1}L^{\frac{1}{4}} \Big(\sum_{\ell=1}^L \|\wl(q^\dagger)-\wl(q)\|_{H^1\II}\Big)^{\frac{1}{2}}.
\end{align*}
Thus, the Gagliardo-Nirenberg interpolation inequality \cite{Brezis:2018} 
\begin{equation*}
  \|w\|_{H^1\II}\le C(\Omega)^2\|w\|_{L^2\II}^{\frac{1}{2}} \|w\|_{H^2\II}^{\frac{1}{2}},\qquad w\in H^2\II,
\end{equation*}
and the regularity 
$\|\wl(q^\dagger)\|_{H^2\II}+\| \wl(q)\|_{H^2\II}\le C\|\gl\|_{H^2(\partial\Omega)}\le CL^{\frac{1}{2}}  $ directly yields
\begin{align*}
    \|q^\dagger-q\|_{L^2\II}&\le C C_0^{-1}C(\Omega)L^{\frac{1}{4}} \Big(\sum_{\ell=1}^L (\|\wl(q^\dagger) \|_{H^2\II}+\|\wl(q^\dagger) \|_{H^2\II} )^{\frac{1}{2}} \|\wl(q^\dagger)-\wl(q)\|_{L^2\II}^{\frac{1}{2}}  \Big)^{\frac{1}{2}}\\
    &\le CL^{\frac{7}{8}} \delta^{\frac{1}{4}}.
\end{align*}
{
This  indicates that the rate $O(\delta^{\frac{1}{4}})$ is optimal with respect to the stability estimate in Theorem \ref{thm:stability}. When the problem data possess higher regularity, the convergence rate may be further improved \cite{Falk:1983,Bal:2011}. For instance, if $q^\dagger \in W^{k+1,p}(\Omega)$, $u\in H^{k+1}(\Omega)\cap W^{k,\infty}(\Omega)$  and a $k$-th order finite element method is employed, the argument in the proof of Theorem \ref{thm:error_estimate} can potentially yield an improved rate of $O(\delta^{\frac{k}{2(k+1)}})$ for $k \ge 1$. The present work focuses on deriving stability and error estimates under weakened regularity assumptions.
}

In practice, the noise level $\delta$ is typically unknown. In this case, Theorem \ref{thm:error_estimate} provides an \textit{a priori} estimate for selecting parameters $h$ and $\alpha$  adaptively \cite{Kaltenbacher:2014,Kaltenbacher:2018,Chen:2022}. To achieve this, we approximate $\delta$ by $\frac{1}{L}\sum_{\ell=1}^{L}\|\wl_h(q_h^*) -z_\delta^{(\ell)}\|_{L^2(\Omega)}$ and choose $h$ and $\alpha$ according to 
    \begin{equation*}
        h_j^2L^{\frac{1}{2}}\sim \frac{1}{L}\sum_{\ell=1}^{L}\|\wl_h(q_h^*) -z_\delta^{(\ell)}\|_{L^2(\Omega)} \quad \text{and} \quad \alpha_j^{\frac{1}{2}}\sim \frac{1}{L}\sum_{\ell=1}^{L}\|\wl_h(q_h^*) -z_\delta^{(\ell)}\|_{L^2(\Omega)},
    \end{equation*}
    where $q_h^*$ and $\wl_h$ are solved using parameters $h_{j-1}$ and $\alpha_{j-1}$ starting from {the initial guesses} $h_{0}$ and $\alpha_{0}$.

\end{remark}

\begin{remark}\label{rmk:polygon}
    In two dimensions, the above analysis can be extended to the case where $\Omega$ is a convex polygon. We parameterize $\partial\Omega$ by arc length and generate $H^{\frac{1}{2}}(\partial\Omega)$ orthonormal basis using the eigenvalues and eigenfunctions of Laplace--Beltrami operator on $\partial\Omega$. Indeed, the eigenfunctions are trigonometric functions on each edge which are continuous at each vertex. Therefore, with appropriate normalization, we obtain the $H^{\frac{1}{2}}(\partial\Omega)$ orthonormal basis. With the same argument as in Remark~\ref{rmk:high_regularity}, the following upper bound holds with high probability
    \begin{equation*} 
        \sum_{i=1}^{N}\|\gl\|_{H^2(\Gamma_i)}\le C L^{\frac{1}{2}}, \quad   \ell=1,\dots,L,
    \end{equation*}
     where $\Gamma_i$, $i=1,\dots,N$ are the edges of the polygon $\Omega$. As a consequence, the forward problem \eqref{eqn:cts_constraint} admits $H^2(\Omega)$ solutions \cite[Theorem 5.1.2.4]{Grisvard:1985} and the $L^2(\Omega)$ error estimate $\|\wl (q^\dag)-\wl_h(q^\dag)\|_{L^2(\Omega)}\le c h^2L^{\frac{1}{2}}$ holds as a consequence of  \cite[Corollary 3.29]{Ern:2004}.
\end{remark}

\begin{remark}\label{rmk:compare}
{The discussion so far focuses on the output least-squares formulation \eqref{eqn:dis_functional}--\eqref{eqn:dis_constraint}. An alternative approach is to reformulate the inverse problem \eqref{eqn:IDP} as a transport equation for $q$; see, e.g., \cite{Richter:1981inverse,Richter:1981numerical,Bal:2011}. Since the resulting transport equation is linear for fixed internal data, this approach allows for the construction of efficient direct solvers. Nevertheless, it involves differentiating the measured data $z_\delta^{(\ell)}$, and therefore requires stronger regularity assumptions or additional smoothing to control the instability caused by noise. In contrast, the output least-squares formulation employed in this work can accommodate $L^2$-noise and is the framework under which our stability and error analysis is developed. However, the resulting optimization problem is nonconvex and may admit local minima (even though this may not be the case \cite{alberti-etal-2026}). Consequently, theoretical global convergence is not guaranteed in general. This highlights a trade-off between the regularity requirements of direct transport-based solvers and the optimization challenges of variational approaches.}
\end{remark}

\section{Quantitative Photoacoustic Tomography}\label{sec:QPAT}
In this section, we extend the argument to the numerical inversion scheme for quantitative photoacoustic tomography.
We consider the case where radiation propagation is approximated by a second-order elliptic  equation \eqref{eqn:QPAT}.
Our objective is to numerically reconstruct the true diffusion coefficient $D^\dag $ and absorption coefficient $\sigma^\dag$
from multiple internal observations
$$H^{(\ell)}(x) = \sigma^\dag u^{(\ell)}(x;D^\dag,\sigma^\dag)\quad \text{for all}~~x\in\Omega,$$
where $u^{(\ell)}:=u^{(\ell)}(D^\dag,\sigma^\dag)$ denotes the solution to the elliptic equation \eqref{eqn:QPAT}
with parameters $D^\dag$ and $\sigma^\dag$, and associated with the Dirichlet boundary illuminations $\gl$, $\ell=1,2,\ldots,L+1$:

\begin{equation}\label{eq:u-qpat-3}
    \left\{
        \begin{aligned}
-\nabla \cdot (D^\dag \nabla u^{(\ell)})+\sigma^\dag u^{(\ell)}&=0 && \text { in  $\Omega$,} \\
u^{(\ell)}&=\gl && \text { on $\partial \Omega$.}
        \end{aligned}
    \right.
\end{equation}

We need the following assumptions on the parameters and boundary data. In particular, as in the previous section, we assume the parameters to be known in $\Om\setminus\Omp$.

\begin{assum}\label{assum:QPAT}
We assume that the parameters and boundary data satisfy the following assumptions.
\begin{itemize}
\item[(i)] Let $\Omega\subset \mathbb{R}^d$ ($d=2,3$) be a bounded domain with  $C^{1,1}$ boundary $\partial\Omega$.  The exact diffusion coefficient $D^\dag \in W^{2,p}\II \cap \mathcal{A}_D$ with $p>d$ and the exact absorption coefficient $\sigma^\dag \in  \mathcal{A}_\sigma$,  where
\begin{equation*}
\begin{aligned}
&\mathcal{A}_D = \{D\in W^{1,\infty}\II:~0< \Lambda_D^{-1}\le D\le \Lambda_D \mbox{ in }\Om   \mbox{, }  D=D^\dag \mbox{  in }\Om\setminus\Omp\} ~\text{and}\\
&\mathcal{A}_\sigma = \{\sigma\in L^{\infty}\II: ~0<\Lambda_{\sigma}^{-1} \le \sigma\le \Lambda_{\sigma} \mbox{ a.e.\ in }\Om   \mbox{, }  \sigma=\sigma^\dag \mbox{ a.e.\ in }\Om\setminus\Omp\},
\end{aligned}
\end{equation*}
with some \text{a priori} known positive constants $\Lambda_D$ and $\Lambda_{\sigma}$.
\item[(ii)] Let $g^{(1)} \equiv 1$, and  $g^{(\ell)}$ (with $\ell=2,\dots,L+1$) be independent and identically distributed random boundary data given by the expansion \eqref{eqn:reps_g} satisfying Assumption~\ref{assum:high_regularity}.
\end{itemize}
\end{assum}
\begin{remark}\label{rmk:positive-bdry}
{The representation \eqref{eqn:reps_g} provides a general framework for random illuminations. In practical QPAT, however, the boundary data $g$ represents the photon density and must therefore satisfy $g \ge 0$. Although this paper utilizes the zero-mean representation \eqref{eqn:reps_g} to simplify the mathematical derivation, the gap between this model and physical reality can be bridged by considering the shifted boundary data:
    \begin{equation*}
        \tilde g^{(\ell)} = a +g^{(\ell)}=a+ \sum_{k=1}^{M} a_k^{(\ell)} e_k, \qquad \ell=2,\dots,L+1,
    \end{equation*}
    where $a > 0$ is a sufficiently large constant. In particular, by the upper bound of boundary data \eqref{eqn:bounded_bdry_H2} and the Sobolev embedding theory, it follows 
    \begin{equation*}
        \max_{\ell=2,\dots,L+1}\|\gl\|_{L^\infty(\partial\Omega)}\le c\max_{\ell=2,\dots,L+1}\|\gl\|_{H^2(\partial\Omega)}\le cL^{\frac{1}{2}},
    \end{equation*}
    where $c>0$ depends only on $\Omega$. Taking $a> cL^{\frac{1}{2}}$, the positivity condition is satisfied with probability greater than $1-L^d\exp(-C_1L)-L\exp(-C_2M)$. Furthermore, this modification maintains the required non-zero condition \eqref{eqn:nonzero_condition} for the elliptic problem, since the constant shift of boundary data leads to a constant shift of the corresponding solution to \eqref{eqn:IDP}, which preserves the nonzero condition \eqref{eqn:nonzero_condition}.
    }
\end{remark}

We assume that the empirical observational data, denoted by $Z_\delta^{(\ell)}$ is noisy in the sense that 
\begin{align}\label{eqn:noise_qpat}
\|  Z_\delta^{(\ell)} - H^{(\ell)} \|_{L^2(\Om)} \le \delta, \quad \text{for all} ~~ \ell=1,2,\ldots,L+1.
\end{align} 
Assumption~\ref{assum:QPAT} together with the elliptic maximum principle implies that $ 0<\underline{c}_0 \le H^{(1)} \le 1$ for some positive constant $\underline{c}_0$.  
Without loss of generality, we assume that the empirical observation $Z_\delta^{(1)}$ satisfies the same bound $0< \underline{c}_0 \le  Z_\delta^{(1)} \le 1$. Indeed, otherwise, it is enough to project $Z_\delta^{(1)}$ pointwise onto $[\underline{c}_0,1]$, which preserves \eqref{eqn:noise_qpat}.
 
For $\ell=1,2,\ldots,L$, we define
\begin{equation*}
q^\dag = D^\dag |u^{(1)}|^2,~~ w_{\delta}^{(\ell)} = \frac{Z_{\delta}^{(\ell+1)}}{Z_\delta^{(1)}} ,~~  w^{(\ell)}= \frac{H^{(\ell+1)}}{H^{(1)}} = \frac{u^{(\ell+1)}}{u^{(1)}}  ~\text{in}~ \Omega,
\end{equation*}
and
$$f^{(\ell)} = \frac{g^{(\ell+1)}}{g^{(1)}} = g^{(\ell+1)}  ~\text{on}~\partial\Omega.$$
It is straightforward to observe that
\begin{equation*}
	\begin{split}
		\quad \| w_{\delta}^{(\ell)}  - w^{(\ell)} \|_{L^2(\Om)} &\le \Big\| \frac{Z_{\delta}^{(\ell+1)}H^{(1)} -H^{(1)}H^{(\ell+1)}}{H^{(1)} Z_\delta^{(1)}} \Big\|_{L^2(\Om)}
		+ \Big\| \frac{H^{(1)}H^{(\ell+1)}-Z_\delta^{(1)}H^{(\ell+1)}}{H^{(1)} Z_\delta^{(1)}} \Big\|_{L^2(\Om)}\\
		&\le \frac{1}{\underline{c}_0^2}
		\Big( \| H^{(1)}(Z_{\delta}^{(\ell+1)}-H^{(\ell+1)}) \|_{L^2(\Om)}
		+  \| H^{(\ell+1)}(H^{(1)}-Z_\delta^{(1)})\|_{L^2(\Om)} \Big) \le c\delta.
	\end{split}
\end{equation*} 
A direct calculation (\cite{Bal:2010,Bal:2011}) shows that $w^{(\ell)}$ is the solution of the following elliptic equation
\begin{equation}\label{eqn:ellip-1}
    \left\{
        \begin{aligned}
            -\nabla \cdot (q^\dag \nabla w^{(\ell)}) & =0, && \mbox{in }\Om, \\
            w^{(\ell)} & = f^{(\ell)}, && \mbox{on } \partial\Om.
        \end{aligned}
    \right.
\end{equation}
Thus, the first step of the reconstruction algorithm consists of the recovery of $q^\dag$ from the practical observations $w_{\delta}^{(\ell)}$. This is the inverse diffusivity problem discussed in Section~\ref{sec:IDP}. Indeed, Assumption \ref{assum:QPAT} and elliptic regularity \cite[Theorem 9.15]{Gilbarg:1977} imply $u^{(1)}\in W^{2,p}(\Omega)$ and hence $q^\dagger\in W^{2,p}(\Omega)$. By the bounds on $D^\dag$ and the maximum principle, we may assume that the diffusivity coefficient $q^\dagger$ has positive lower and upper bounds $ 0<\Lambda_q^{-1} \le q^{\dagger}\le \Lambda_q $.
Moreover, since $g^{(1)}\equiv 1$, the boundary data $\fl=g^{(\ell+1)}$  still satisfy Assumption~\ref{assum:high_regularity}  and the non-zero condition given in Proposition~\ref{prop:nonzero_condition} holds for equation \eqref{eqn:ellip-1}. Therefore, as in Section~\ref{sec:fem} we propose to consider the following least-squares formula with $H^1\II$-seminorm penalty:
\begin{equation}\label{eqn:dis_functional_qpat}
	\min_{q_h\in \Uad_{q,h}} J_{\al,h}(q_h)=\frac{1}{2}\sum_{\ell=1}^{L}\| \wl_{h}(q_h)-w_{\delta}^{(\ell)} \|_{L^2(\Om)}^2+\frac{\al L}{2}\|\nabla q_h \|_{L^2\II}^2,
\end{equation}
where the admissible set $\Uad_{q,h}$ is defined in \eqref{eq:Aqh} and $\wl_{h}(q_h)\in V_h$ is the weak solution of
\begin{equation}\label{eqn:dis_constraint_qpat}
	\left\{
	\begin{aligned}	
		(q_h \nabla \wl_{h},\nabla v_h)& =0, && \forall v_h\in \mathring{V}_h, \\
		\wl_{h} & =I_h^{\partial}\fl, && \mbox{on } \partial\Omega.
	\end{aligned}
	\right.
\end{equation}
The following error analysis is a direct consequence of Theorem~\ref{thm:error_estimate}. %\ga{[Shall we try to be more precise with the constants $C$ in the following two results?]}
\begin{proposition}\label{prop:error_estimate_qpat}
Suppose Assumption~\ref{assum:QPAT} holds valid and set $q^\dagger=D^\dagger|u^{(1)}|^2$. 
   Let $q_h^*\in \Uad_{q,h}$ be  a minimizer of problem \eqref{eqn:dis_functional_qpat}-\eqref{eqn:dis_constraint_qpat}. Set $\xi= h^2 L^{\frac{1}{2}}+\delta+\al^{\frac{1}{2}}$. Then, with probability greater than \eqref{eqn:probability_nonzero},
    we have %\ga{** we need to say something about $\epsilon$ **}
    \begin{equation*} 
        \|q^{\dag}- q_h^*\|_{L^2(\Omp)}^2\le C  L^2(1+\alpha^{-\frac{1}{2}}\xi)\left(h+ h^{1-\epsilon}(1+\alpha^{-\frac{1}{2}}\xi)+\min\left(1,h + h^{-1}L^{-\frac{1}{2}}\xi \right)  \right),
    \end{equation*}      
where $\epsilon = 0$ when $d = 2$, and $\epsilon > 0$ is arbitrarily small when $d = 3$. Here, $C$ is a constant independent of $h$, $\delta$, $\alpha$, and $L$, but may depend on $\epsilon$, $\Omega$, and $q^{\dagger}$.    
\end{proposition}

The second step of the inverse algorithm is to recover $u^{(1)}$. The reconstruction of $D^{\dag}$ and $\sigma^\dag$ will follow immediately by using the relations $D^\dag=q^\dag/|u^{(1)}|^2$ and $\sigma^\dag=H^{(1)}/u^{(1)}$.
Since $u^{(1)}|_{\partial\Omega}=g^{(1)}\equiv 1 $, by \eqref{eq:u-qpat-3} we have that $v=1/u^{(1)}-1$ satisfies the following boundary value problem
\begin{equation}\label{eqn:ellip-2}
	\left\{
	\begin{aligned}
		-\nabla \cdot (q^\dag \nabla  v) & =H^{(1)}, && \mbox{in }\Om, \\
		v  & = 0, && \mbox{on } \partial\Om.
	\end{aligned}
	\right.
\end{equation}
We are now ready to show the error bound of the numerically recovered parameters.
\begin{theorem}\label{thm:error_estimate_qpat}
	Suppose that Assumption~\ref{assum:QPAT} holds valid and set $q^\dagger=D^\dagger|u^{(1)}|^2$.  Let  $q_h^*\in \mathcal{A}_{q,h}$ 
    be such that $\|q^\dag-q_h^*\|_{L^2(\Omp)}\le \eta$ for some $\eta\ge 0$ and set the reconstructed coefficient $q^*$ as 
    \begin{equation*}
        q^*=\left\{
            \begin{aligned}
                &q_h^* && \mbox{in }\Omp,\\
                &D^\dagger (Z_\delta^{(1)}/\sigma^\dagger)^2&& \mbox{in }\Omega\setminus\Omp.
            \end{aligned}
        \right.
    \end{equation*}
    Let $v_h\in \mathring{V}_h$ solve
	\begin{equation}\label{eqn:ellip-2_dis}
			 (q^* \nabla  v_h,\nabla \varphi_h)  =(Z^{(1)}_\delta,\varphi_h), \quad \forall \varphi_h\in \mathring{V}_h.
	\end{equation}
	Then there holds
	\begin{equation*}
	 		\|v-v_h\|_{L^2\II}\le C(h+\eta+\delta).
	\end{equation*}
	Moreover, set $D^*=q^*|v_h+1|^2$ and $\sigma^*=Z^{(1)}_\delta (v_h+1)$, we have
	\begin{equation*}
		\|D^\dag-D^*\|_{L^2\II}\le  C(h+\eta+\delta)\quad \text{and} \quad \|\sigma^\dag-\sigma^*\|_{L^2\II}\le C(h+\eta+\delta),
	\end{equation*}
  where $C$ is a constant independent of $h$, $\delta$ and $\eta$.
\end{theorem}
\begin{proof}
With an abuse of notation, several positive constants independent of $h$, $\eta$ and $\delta$ will be denoted by the same letter $C$.    We observe that 
    \begin{equation*}
	\begin{split}
		\quad \| q^\dagger  - q^* \|_{L^2(\Om\setminus\Omp)} &= \Big\| D^\dagger \frac{(H^{(1)})^2}{(\sigma^\dagger)^2} - D^\dagger \frac{(Z_\delta^{(1)})^2}{(\sigma^\dagger)^2}\Big\|_{L^2(\Om\setminus\Omp)}\\
		&\le \Lambda_D\Lambda_\sigma^2
		 \| H^{(1)} +Z_{\delta}^{(1)}  \|_{L^\infty(\Om\setminus\Omp)}\| H^{(1)} -Z_{\delta}^{(1)}  \|_{L^2(\Om\setminus\Omp)}  \le C\delta.
	\end{split}
    \end{equation*} 
	By equation \eqref{eqn:ellip-2} and \eqref{eqn:ellip-2_dis}, we have
	\begin{align*}
		(q^*(\nabla P_hv-\nabla v_h),\nabla \varphi_h)&=(q^*(\nabla P_h v-\nabla v),\nabla \varphi_h)+(q^*(\nabla v-\nabla v_h),\nabla \varphi_h)\\	
		&=(q^*(\nabla P_h v-\nabla v),\nabla \varphi_h)+	((q^*-q^\dag)\nabla v,\nabla \varphi_h)+(H^{(1)}-Z^{(1)}_{\delta},\varphi_h).
	\end{align*}
	Taking  $\varphi_h=P_h v- v_h$, Cauchy-Schwarz inequality and Poinc\'{a}re's inequality yield
	\begin{align*}
	 	\|\nabla \varphi_h\|_{L^2\II}^2
	 	\le& C\|\nabla (P_hv-v)\|_{L^2\II} \|\nabla \varphi_h\|_{L^2\II}+C\|\nabla v\|_{L^\infty\II}\|q^*-q^\dag\|_{L^2\II}\| \nabla \varphi_h\|_{L^2\II}\\
	 	&+C\|H^{(1)}-Z^{(1)}_{\delta}\|_{L^2\II}\|\nabla \varphi_h\|_{L^2\II}.
	\end{align*}
	By elliptic regularity theory, we have  $v\in H^2\II\cap W^{1,\infty}\II$. Hence, by the projection error  \eqref{eqn:error_P_h}, estimate of $q^*$ and the noise level \eqref{eqn:noise_qpat}, we derive that
	\begin{equation*}
		\|\nabla P_h v-\nabla v_h\|_{L^2\II}=	\|\nabla \varphi_h\|_{L^2\II}\le C(h+\eta+\delta).
	\end{equation*}
	Thus, by Poinc\'{a}re's inequality and the error bound \eqref{eqn:error_P_h}, we conclude
	\begin{equation*}
			\|  v-  v_h\|_{L^2\II}\le \|  P_h v-  v_h\|_{L^2\II}+\|  P_h v-  v\|_{L^2\II}\le  C(h+\eta+\delta).
	\end{equation*}
	Moreover,  direct computation yields
	\begin{align*}
			\|D^\dag-D^*\|_{L^2\II}&=	\left\| \frac{q^\dag}{|u^{(1)}|^2}-q^* |v_h+1|^2\right\|_{L^2\II}=\left\|  q^\dag|v+1|^2-q^* |v_h+1|^2\right\|_{L^2\II} \\
			&\le \left\| (q^\dag -q^*)|v+1|^2  \right\|_{L^2\II}+\left\| q^*(| v+1|^2-|v_h+1|^2) \right\|_{L^2\II}\\
			&\le C(\eta+\delta)+ C(h+\eta+\delta)\le  C(h+\eta+\delta),
	\end{align*}	
	and
	\begin{align*}
			\|\sigma^\dag-\sigma^*\|_{L^2\II}&=\left\| \frac{H^{(1)}}{u^{(1)}}-Z^{(1)}_\delta(v_h+1)\right\|_{L^2\II}= \left\| H^{(1)}(v+1)-Z^{(1)}_\delta(v_h+1)\right\|_{L^2\II} \\
			&\le \| (H^{(1)}-Z_\delta^{(1)})(v+1)\|_{L^2\II}+\| Z_\delta^{(1)}(v-v_h)\|_{L^2\II}\\
				&\le C\delta+ C(h+\eta+\delta)\le  C(h+\eta+\delta).
	\end{align*}
\end{proof}
\begin{remark}\label{rmk:conv_rate}
	The error analysis in Proposition~\ref{prop:error_estimate_qpat} and Theorem~\ref{thm:error_estimate_qpat} provide a guideline for choosing the mesh size $h$ and regularization parameter $\al$, see also Remark~\ref{rmk:error_q_delta}. Indeed, by choosing $h^2L^{\frac{1}{2}}\sim \delta $ and $\al \sim \delta^2$, with probability greater than \eqref{eqn:probability_nonzero},
	there holds 
	\begin{equation*}
		\|D^\dag-D^*\|_{L^2\II} + \|\sigma^\dag-\sigma^*\|_{L^2\II}\le CL^{\frac{7}{8}}   \delta^{\frac{1}{4}-\epsilon},
	\end{equation*}
    with $\epsilon=0$ for $d=2$, $\epsilon>0$ arbitrary small for $d=3$.
\end{remark}

\section{Numerical Implementation and Results}\label{sec:numer}
In this section, we provide numerical reconstructions of the diffusion coefficient $D^\dag$ and the absorption coefficient $\sigma^\dag$ based on the two stage algorithm discussed in Section~\ref{sec:QPAT}. We first solve the optimization problem \eqref{eqn:dis_functional_qpat}-\eqref{eqn:dis_constraint_qpat} and then solve the direct problem \eqref{eqn:ellip-2_dis}. We consider the two-dimensional setting ($d=2$).

\subsection{Numerical implementation}\label{subsec:numer_implementation}
In this part, we introduce the numerical implementation for the reconstruction algorithm. We first describe the generation of the boundary illuminations $\gl$, $\ell=1,\dots,L+1$. Recall in Assumption~\ref{assum:QPAT}(iii), $g^{(1)}\equiv 1$ is fixed and $\gl$ (with $\ell=2,\dots,L+1$) are taken as 
\begin{equation*}
	g^{(\ell)}=\sum_{k=1}^{M} a_k^{(\ell)} e_k ,
\end{equation*}
where  $\{e_k\}_{k=1}^{\infty}$ is a fixed orthonormal basis of $H^{\frac{1}{2}}(\partial \Omega )$ generated by the eigenfunctions of the Laplace--Beltrami operator   $-\Delta $ on $\partial \Omega$, see Remark \ref{rmk:high_regularity}.  The coefficients $a_{k}^{(\ell)}\sim N(0,\theta_k^2)$ are independent and identically distributed random variables satisfying Assumption~\ref{assum:high_regularity} with $\theta_k^2\lesssim k^{-6}$. %\ga{[This should be double-checked, taking into accounts my suggestions above]}. % $\theta_k =k^{-2}$ and $s=\frac{5}{2} $.  
In all the examples, we take the first five terms in the series, i.e.\ $M=5$. With the truncated boundary illuminations, we generate noisy measurements as follows:  
\begin{equation*}
	Z_\delta^{(\ell)}(x)=H^{(\ell)}(x)+\delta\sup_{z\in \Omega}|H^{(\ell)}(z) |\xi(x), \quad \ell=1,\dots,L+1,
\end{equation*}
where $\xi$ follows standard Gaussian distribution, while $\delta$ denotes the level of noise. The exact data $H^{(\ell)}=\sigma^\dag u^{(\ell)}(D^\dag,\sigma^\dag ) $ correspond to the precise values of $D^\dag$ and $\sigma^\dag  $, calculated using a highly refined mesh with $h=\frac{1}{500}$.

In the first step of inversion, we solve the least-squares formula \eqref{eqn:dis_functional_qpat}-\eqref{eqn:dis_constraint_qpat} using conjugate gradient method \cite[Section 5.2]{Nocedal:2006}. The derivative of the least-squares functional is computed in the following.
\begin{lemma}\label{lem:dJ}
    Let $J_{\alpha,h}$ be the objective functional given in  \eqref{eqn:dis_functional_qpat}-\eqref{eqn:dis_constraint_qpat}. The directional derivative of $J_{\alpha,h}$ at $q_h$ along the direction $p_h$, denoted as $J^\prime_{\alpha,h}(q_h)[p_h]$, is expressed as 
    \begin{equation*}
        J^\prime_{\alpha,h}(q_h)[p_h]=\frac{1}{2}\sum_{\ell=1}^{L}(p_h\nabla \wl_h,\nabla \phi_h)+\alpha L (\nabla q_h,\nabla p_h),
    \end{equation*}
    where $\phi_h^{(\ell)}\in \mathring{V}_h,\,\ell=1,\dots,L$ satisfying 
    \begin{equation}\label{eqn:dis_adjoint_eq}
        (q_h\nabla \phi^{(\ell)}_h,\nabla \varphi_h)=(\wl_\delta-\wl_h,\varphi_h),\quad \forall\varphi_h\in \mathring{V}_h.
    \end{equation}
\end{lemma}
\begin{proof}
    A direct computation yields
    \begin{equation*}
        J^\prime_{\alpha,h}(q_h)[p_h]= \frac{1}{2}\sum_{\ell=1}^{L}( \wl_h-\wl_\delta, {w^{ (\ell)\prime}_h}(q_h)[p_h])+\alpha L (\nabla q_h,\nabla p_h),
    \end{equation*}
where ${w^{ (\ell)\prime}_h}(q_h)[p_h]$ is the derivative of $ \wl_h(q_h)$ at $q_h$ along the direction $p_h$. Note that ${w^{ (\ell)\prime}_h}(q_h)[p_h] \in \mathring{V}_h$ satisfies 
\begin{equation*}
    ( {w^{ (\ell)\prime}_h}(q_h)[p_h] , \varphi_h)=-(p_h\nabla \wl_h(q_h), \nabla \varphi_h), \forall \varphi_h\in V_h.
\end{equation*} 
This and the weak formulation of $\phi_h$ yields the desired result.
\end{proof}
According to Lemma \ref{lem:dJ}, the derivative $J^\prime_{\alpha,h}(q_h)$ belongs to the dual space of $V_h$, which is not necessarily a FEM space. To overcome this issue, we apply a smoothing map to pull it back to the space $V_h\subset H^1(\Omega)$ and obtain a sufficiently regular direction for updating. The gradient $G_h\in V_h $ of the objective functional $J_{\alpha,h}$ at $q_h$ can be computed as follows
\begin{equation}\label{eqn:grad_J}
    (\nabla G_h,\nabla p_h)= J^\prime_{\alpha,h}(q_h)[p_h]=\frac{1}{2}\sum_{\ell=1}^{L}(p_h\nabla \wl_h,\nabla \phi_h)+\alpha L (\nabla q_h,\nabla p_h),\quad \forall p_h\in V_h.
\end{equation}
We summarize the reconstruction procedure in Algorithm \ref{alg:QPAT}. {Throughout, the initial guess for the optimization problem \eqref{eqn:dis_functional_qpat}-\eqref{eqn:dis_constraint_qpat} is set to a constant background value. The minimization procedure typically converges  within 100 iterations.}
\begin{algorithm}[htbp]
\caption{}\label{alg:QPAT}
\begin{algorithmic}[1]
\Require  Boundary data $\gl,\, \ell=1,\dots,L+1$ and corresponding measurement $Z^{(\ell)}_\delta$.
\State Step 1: inverse diffusivity problem.
\State \hspace{1.12cm} Define $\wl_\delta=Z^{(\ell+1)}_\delta / Z^{(1)}_\delta$ and $f^{(\ell)}=g^{(\ell+1)}/ g^{(1)}$.
\State \hspace{1.12cm} Choose initial guess $q_h^0$ and set $j=0$.
\State \hspace{1.12cm} Solve the forward problem \eqref{eqn:dis_constraint_qpat} with $q_h=q_h^j$.
\State \hspace{1.12cm} Solve the adjoint problem \eqref{eqn:dis_adjoint_eq} and compute the gradient $G_h=G_h(q_h^j)$ by \eqref{eqn:grad_J}.
\State \hspace{1.12cm} Compute the descent direction
\begin{equation*}
    d_j=-G_h(q_h^j)+\beta_{j-1}d_{j-1},
\end{equation*}
\hspace{1.12cm} with the conjugate coefficient $\beta_{j-1}$ given by
\begin{equation*}
    \beta_{j-1}=\frac{\|G_h(q_h^j)\|_{L^2(\Omega)}}{\|G_h(q_h^{j-1})\|_{L^2(\Omega)}}
\end{equation*}
\hspace{1.12cm} and setting $\beta_{-1}=0$
\State \hspace{1.12cm} Calculate the step length $\gamma_j$ by a line search procedure \cite[eq. (5.43)]{Nocedal:2006}.
\State \hspace{1.12cm} Update the diffusivity $  q_h^{j+1}=q_h^j+\gamma_jd_j  $.
%\begin{equation*}
%    q_h^{j+1}=q_h^j+\gamma_jd_j.
%\end{equation*}
\State \hspace{1.12cm} Increase $j$ by one and go to line 4.
\State \hspace{1.12cm} Repeat the above procedure until a stopping criterion is satisfied.
\State Step 2: direct problem 
\State \hspace{1.12cm} Solve forward problem \eqref{eqn:ellip-2_dis} to get $v_h$.
\State \hspace{1.12cm} Set $D^*=q^*|v_h+1|^2$ and $\sigma^*=Z^{(1)}_\delta (v_h+1)$.
\end{algorithmic}
\end{algorithm}

\subsection{Numerical experiments}\label{subsec:numer_experiments}

\begin{figure}[htbp]
	\centering
	\begin{tabular}{ccc}
		\includegraphics[width=0.25\textwidth]{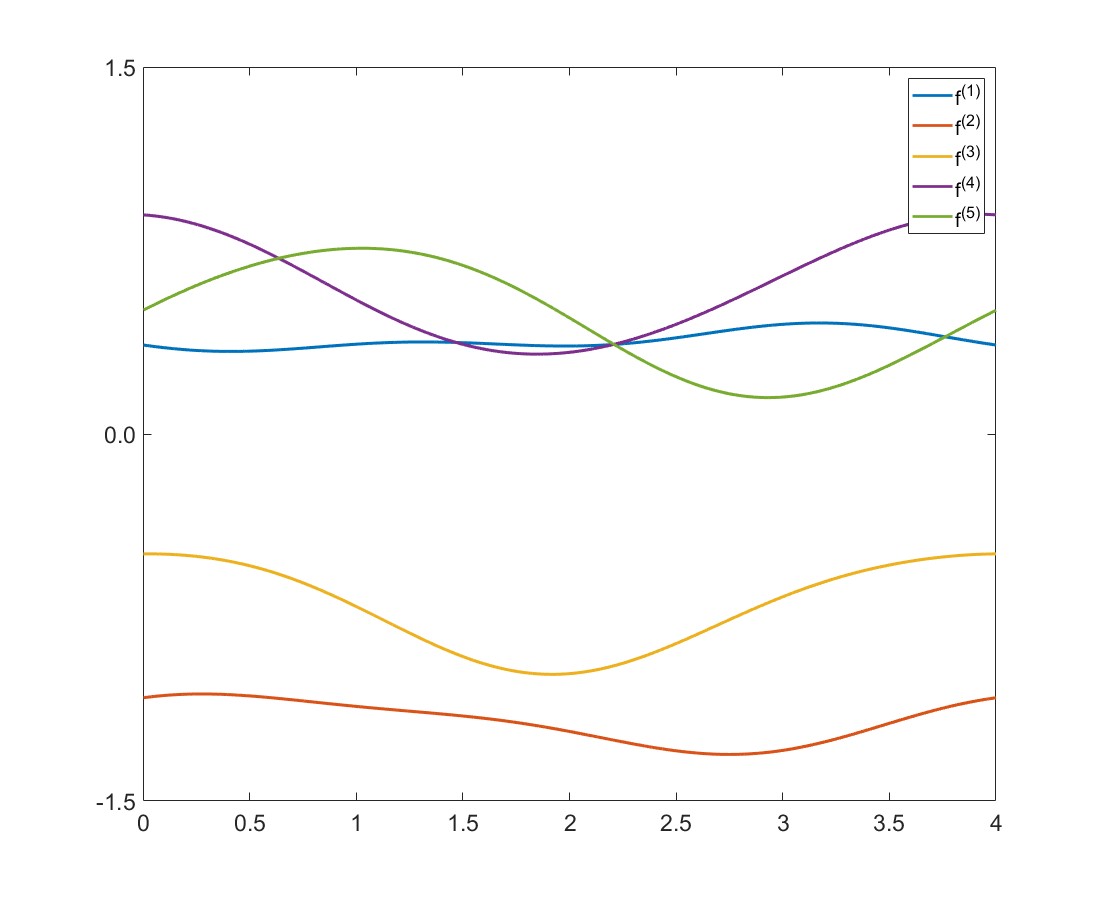}&
		\includegraphics[width=0.25\textwidth]{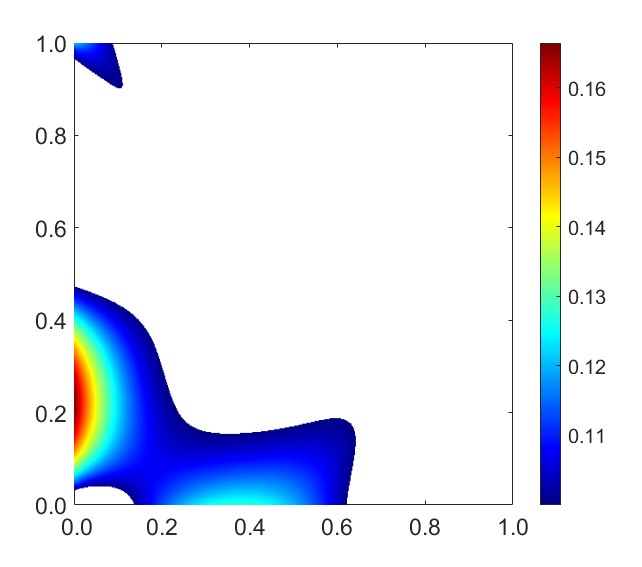}&
		\includegraphics[width=0.25\textwidth]{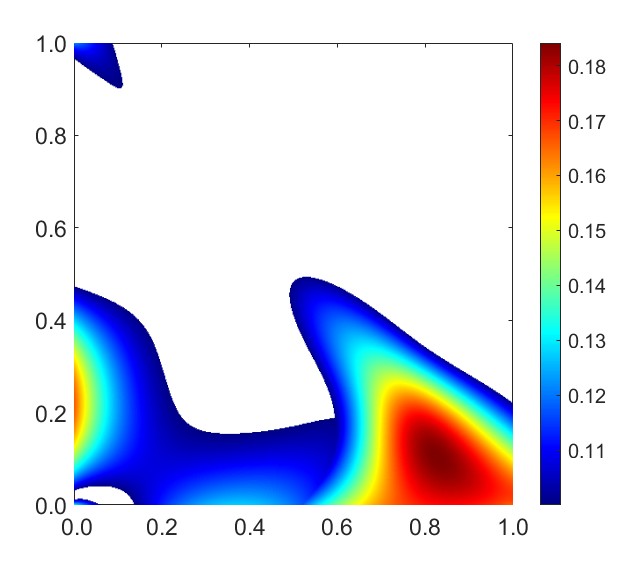}\\
		(a) Boundary $\fl$ & (b) $L=1$ & (c) $L=2$ \\
		\includegraphics[width=0.25\textwidth]{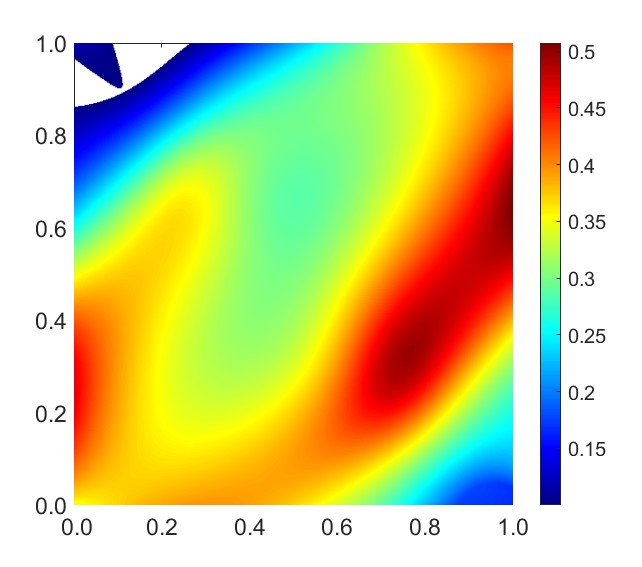}&
		\includegraphics[width=0.25\textwidth]{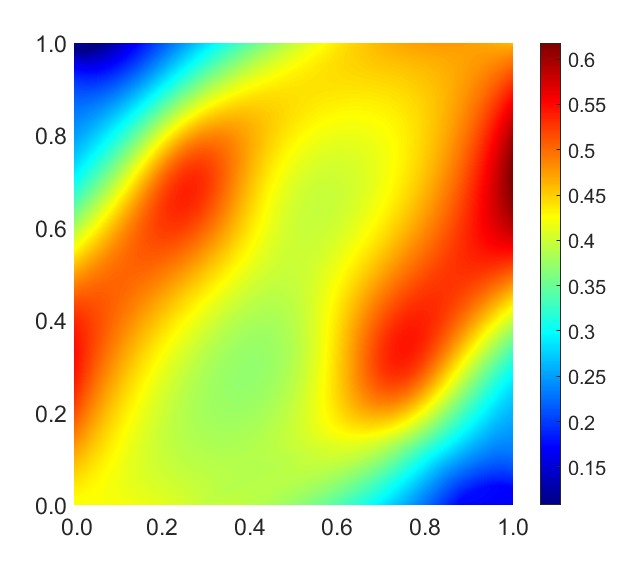}&
		\includegraphics[width=0.25\textwidth]{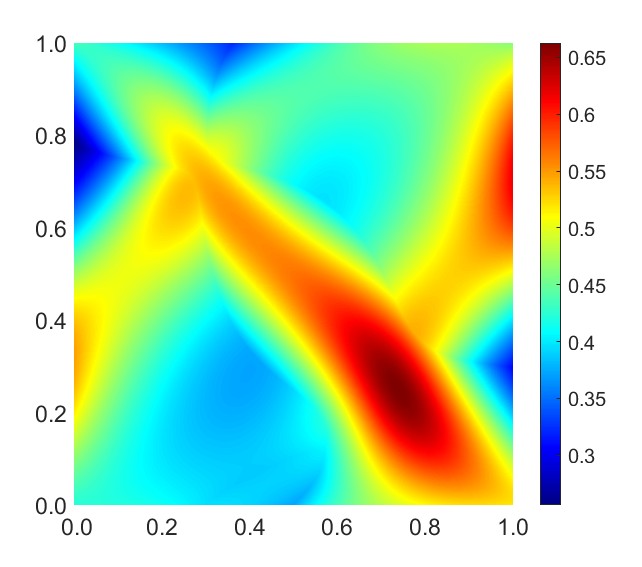}\\
		(d) $L=3$ & (e) $L=4$ & (f) $L=5$ \\
	\end{tabular}
	\caption{Boundary illuminations and the non-zero region of Example~\ref{ex:1}. Top left: plot of boundary data $\fl=g^{(\ell+1)}$. Top middle to bottom right: region where the non-zero condition is satisfied as the number of boundary inputs increases.}
	\label{Fig:ex1_nonzero}
\end{figure}

In this part, we provide numerical verification of the non-zero condition \eqref{eqn:nonzero_condition} and the numerical reconstructions of the diffusion coefficient $D^\dag $ and the absorption coefficient $\sigma^\dag$. { While Proposition \ref{prop:nonzero_condition} establishes that the non-zero condition holds with overwhelming probability, the associated constants $C_0, C_1,$ and $C_2$ are not explicitly determined. Here, we numerically verify this condition.}  To verify the non-zero condition, we plot the region:
\begin{equation*}
	  \{x\in \Omega:\max_{\ell=1,\dots,L}|\nabla w^{(\ell)}(x)\cdot \nu |\ge C_0\}, 
\end{equation*}
where $ w^{(\ell)}(x)=H^{(\ell+1)}/ H^{(1)}$. In the following numerical experiments, we fix the direction $\nu=(1,0)$ and the threshold $C_0=0.1$. To quantify the performance of the numerical reconstruction, we introduce the following relative $L^2(\Omega)$ error:
\begin{equation*}
	e_D=\|D_h^*-D^\dag\|_{L^2\II}/\|D^\dag\|_{L^2\II}
	\quad \mbox{and}\quad e_\sigma=\|\sigma_h^*-\sigma^\dag\|_{L^2\II}/\|\sigma^\dag\|_{L^2\II}.
\end{equation*}
We start with the following examples with smooth coefficients.
\begin{example}\label{ex:1}
	$\Omega=(0,1)^2$,  $D^{\dag}(x,y)=2+\sin(2\pi x)\sin(2\pi y)$ and $\sigma^{\dag}=6+4\sigma_1+4\sigma_2$ with $\sigma_1(x,y)=e^{-20(x-0.3)^2-20(y-0.7)^2}$ and $\sigma_2(x,y)=e^{-20(x-0.7)^2-20(y-0.3)^2}$. %The exact solution $\ul$ are computed  with mesh size $h=1/500$.
\end{example}
In Figure~\ref{Fig:ex1_nonzero}(a), we plot the random boundary data $f^{(\ell)}=g^{(\ell+1)}/g^{(1)}=g^{(\ell+1)}$. We show the region in which the non-zero condition is satisfied with different $L$ in Figures~\ref{Fig:ex1_nonzero}(b)-(f). We observe that the region where the non-zero condition is satisfied expands as the number of random boundary data increases. For $L=1$, the non-zero condition is satisfied only in a small region, while for $L=3$, the non-zero condition is satisfied in most parts of the domain $\Omega$. We also notice that as $L$ increases, the lower bound $C_0$ increases, indicating better stability of the inverse problem.

\begin{figure}[htbp]
	\centering
	\begin{tabular}{ccc}
		\includegraphics[width=0.25\textwidth]{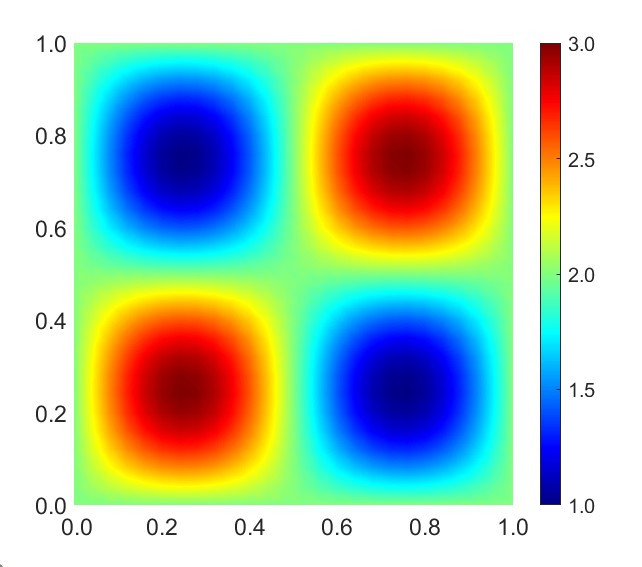}&
		\includegraphics[width=0.25\textwidth]{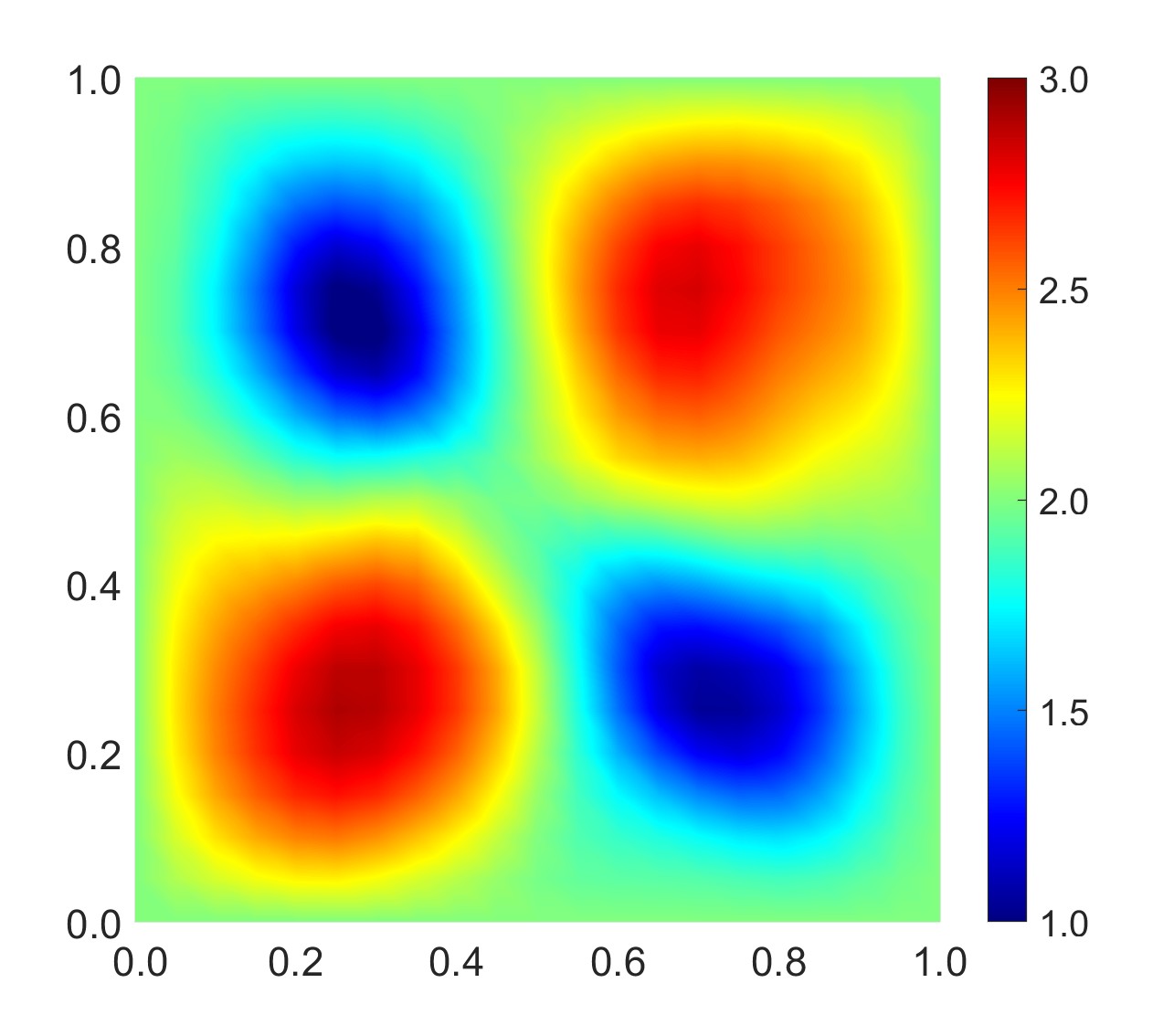}&
		\includegraphics[width=0.25\textwidth]{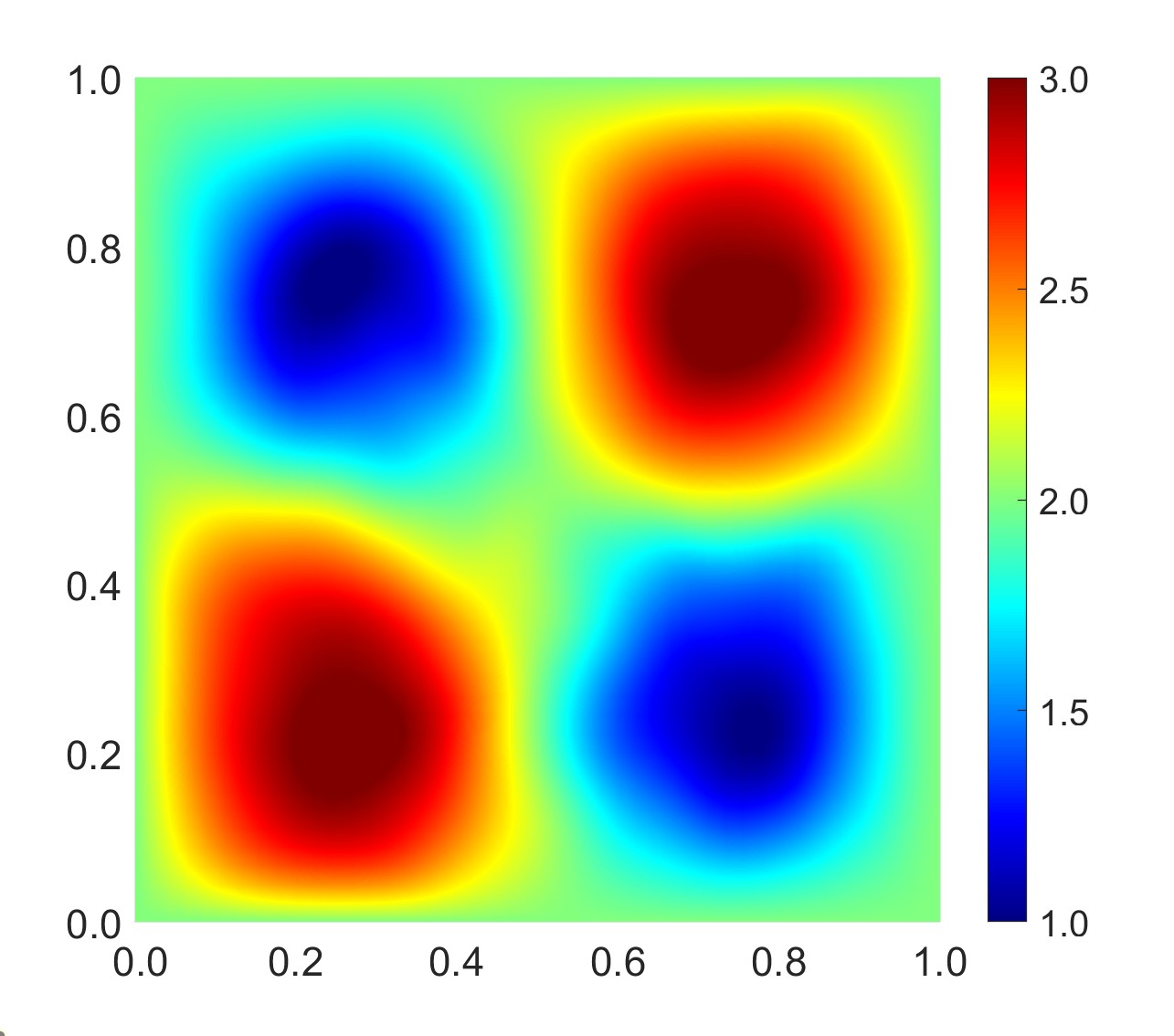}\\
		(a)  $D^\dag$ & (b) $\delta=5e$-$2$ & (c) $\delta=1e$-$2$ \\
		\includegraphics[width=0.25\textwidth]{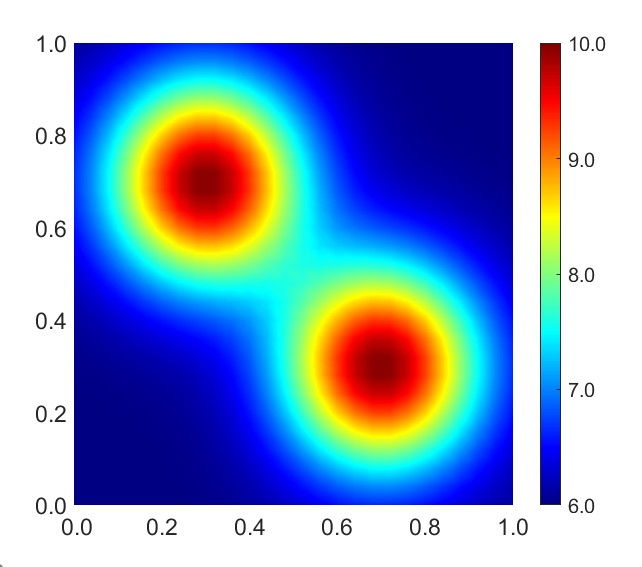}&
		\includegraphics[width=0.25\textwidth]{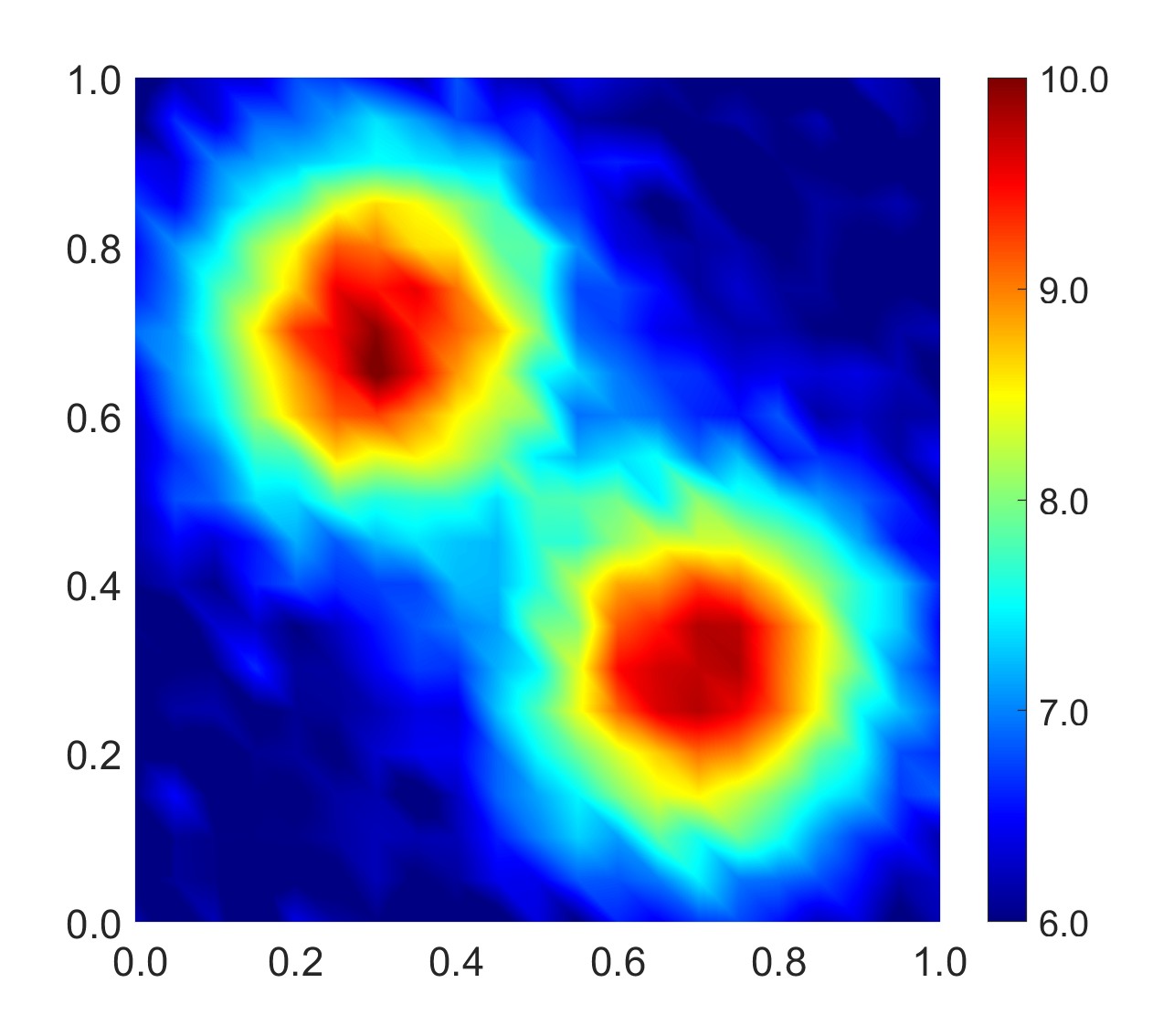}&
		\includegraphics[width=0.25\textwidth]{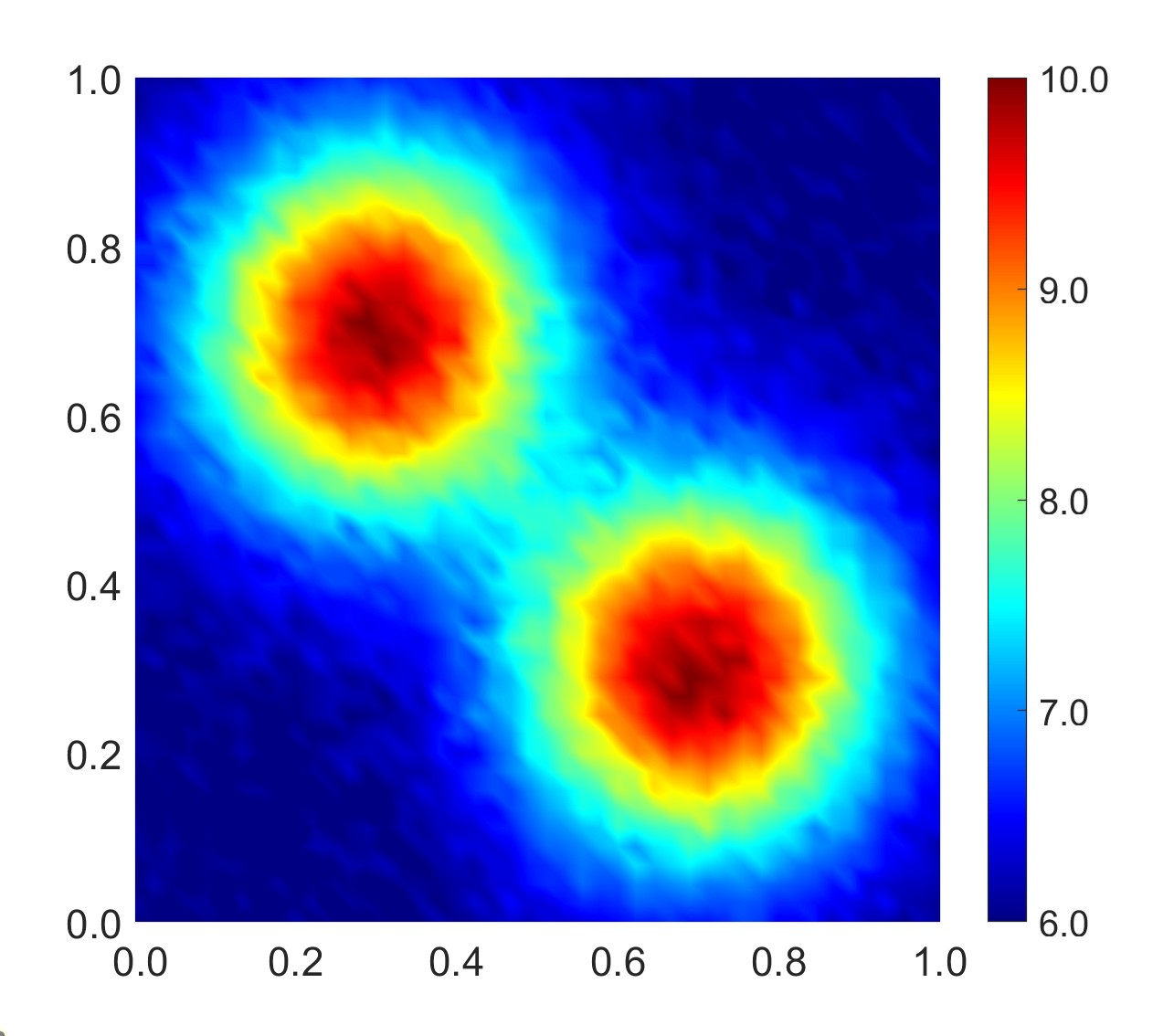}\\
		(d) $\sigma^\dag$ & (e) $\delta=5e$-$2$ & (f) $\delta=1e$-$2$ \\
	\end{tabular}
	\caption{Example~\ref{ex:1}. First row: reconstructions of $D^\dag$.
		Second row: reconstructions of $\sigma^\dag$. }
	\label{Fig:ex1_recon}
\end{figure}

Table~\ref{tab:ex1} displays the convergence rate of the reconstruction errors. The mesh size and the regularization parameter are chosen by following the guidelines in Remark~\ref{rmk:conv_rate} with fixed number of illuminations $L=5$: $h\sim \delta^{\frac{1}{2}}$ and $ \alpha\sim\delta^2$. We initialize the mesh size $h=1/12$ and the regularization parameter $\alpha=3e$-$7$. The numerical results indicate that the error $e_D$ and $e_\sigma$ decay to zero as the noise level tends to zero, with rate $O(\delta^{0.22})$ and $O(\delta^{0.26})$, respectively. These convergence rates are consistent with the rate predicted in Remark~\ref{rmk:conv_rate}, which is $O(\delta^{0.25})$. Figure~\ref{Fig:ex1_recon} shows the recovered diffusion coefficient and absorption coefficient in $5\%$ and $1\%$ noise.  Here we take $h=1/20$, $\alpha=1e$-$6$ for noise level $\delta=5e$-$2$ and  $h=1/45$, $\alpha=4e$-$8$ for  $\delta=1e$-$2$.

\begin{figure}[htbp]
	\centering
	\begin{tabular}{ccc}
		\includegraphics[width=0.25\textwidth]{Figure/ex1_D_exact.jpg}&
		\includegraphics[width=0.25\textwidth]{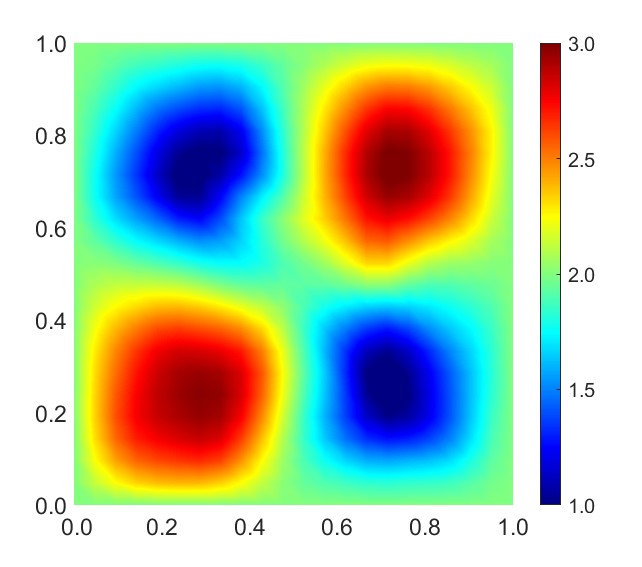}&
		\includegraphics[width=0.25\textwidth]{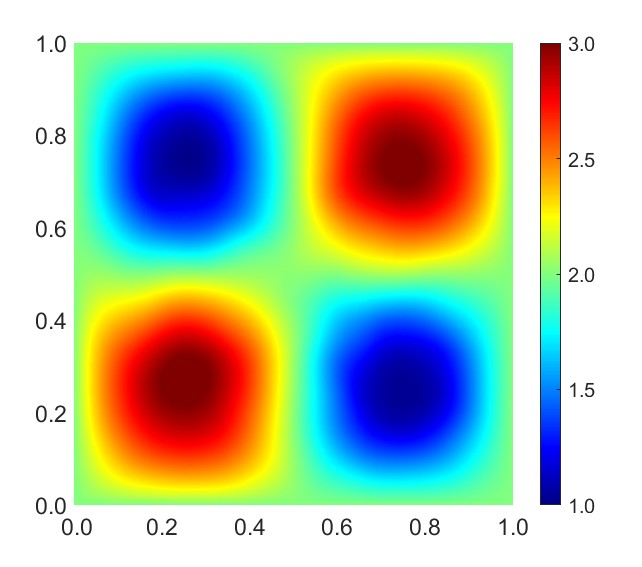}\\
		(a)  $D^\dag$ & (b) $\delta=5e$-$2$ & (c) $\delta=1e$-$2$ \\
		\includegraphics[width=0.25\textwidth]{Figure/ex1_Sigma_exact.jpg}&
		\includegraphics[width=0.25\textwidth]{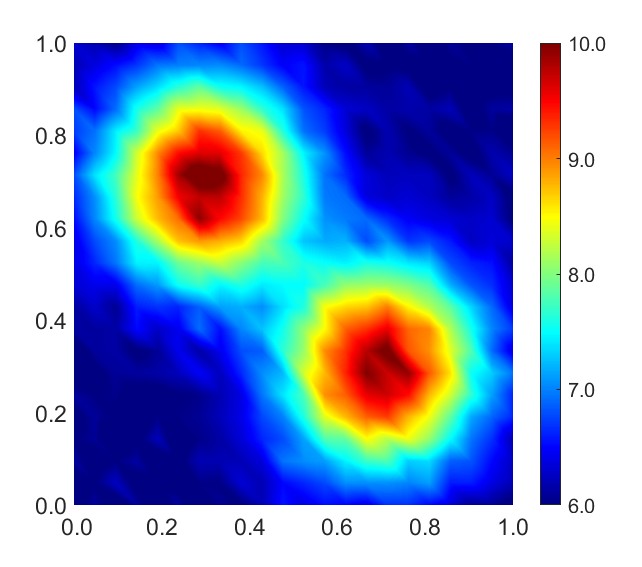}&
		\includegraphics[width=0.25\textwidth]{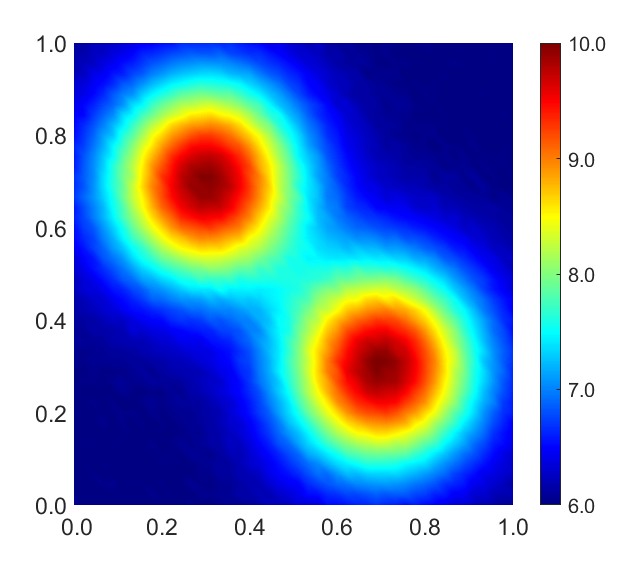}\\
		(d) $\sigma^\dag$ & (e) $\delta=5e$-$2$ & (f) $\delta=1e$-$2$ \\
	\end{tabular}
	\caption{Example~\ref{ex:1} using the adaptive algorithm in Remark \ref{rmk:error_q_delta}. First row: reconstructions of $D^\dag$.
		Second row: reconstructions of $\sigma^\dag$. }
	\label{Fig:ex1_recon_adapt}
\end{figure}

We also investigate the performance of the adaptive algorithm mentioned in Remark \ref{rmk:error_q_delta}. In particular, given initial mesh size $h_0$ and regularization parameter $\alpha_0$, we solve the least-squares problem \eqref{eqn:dis_functional_qpat}-\eqref{eqn:dis_constraint_qpat} with $h_0$ and $\alpha_0$. Then for $j=1,2,\dots$, the mesh size $h_j$ and regularization parameter $\alpha_j$ are chosen from the estimator $\frac{1}{L}\sum_{\ell=1}^{L}\|\wl_h(q_h^*)-\zl_\delta\|_{L^2(\Omega)}$ with $q_h^*$ and $\wl_h$ solved using parameters $h_{j-1} $ and $\alpha_{j-1}$:
    \begin{equation*}
        h_j^2L^{\frac{1}{2}}\sim \frac{1}{L}\sum_{\ell=1}^{L}\|\wl_h(q_h^*) -z_\delta^{(\ell)}\|_{L^2(\Omega)} \quad \text{and} \quad \alpha_j^{\frac{1}{2}}\sim \frac{1}{L}\sum_{\ell=1}^{L}\|\wl_h(q_h^*) -z_\delta^{(\ell)}\|_{L^2(\Omega)}.
    \end{equation*} 
Figure~\ref{Fig:ex1_recon_adapt} shows the recovered diffusion coefficient and absorption coefficient in $5\%$ and $1\%$ noise with the adapted approach. {We observe that the reconstruction quality is comparable or even better than that of Algorithm \ref{alg:QPAT}, in which the mesh size $h$ and regularization parameter $\alpha$ are chosen manually based on prior knowledge of $\delta$.}  Here we take $h_0=1/10$, $\alpha_0=1e$-$5$ for both noise level cases. The resulting parameters are $h=1/21$, $\alpha=4.06e$-7 for $\delta=5e$-2 and   $h=1/45$, $\alpha=1.76e$-8 for $\delta=1e$-2. {These results demonstrate that computational parameters can be selected adaptively without explicit knowledge of the noise level $\delta$, significantly enhancing the robustness and practicality of the reconstruction.}

\begin{table}[htp!]
	\centering
	\caption{The convergence rates for Example~\ref{ex:1} with respect to $\delta$.} \label{tab:ex1}
	\begin{tabular}{c|cccccccc}
		\toprule
		$\delta$   & 1e-2 & 5e-3 & 2e-3 & 1e-3 & 5e-4 & 2e-4  &1e-4   & rate \\
		\midrule
		$e_D$ & 6.53e-2 & 4.17e-2 & 2.98e-2 & 2.61e-2 & 2.33e-2  & 2.30e-2 &  2.12e-2 & {$O(\delta^{0.22})$}  \\
		$e_{\sigma}$ & 1.45e-2  & 7.90e-3 & 5.23e-3 & 4.77e-3 &  4.13e-3  & 4.06e-3  &  3.76e-3 & {$O(\delta^{0.26})$} \\	
		\bottomrule
	\end{tabular}
\end{table}

\begin{figure}[htbp]
	\centering
	\begin{tabular}{ccc}
		\includegraphics[width=0.25\textwidth]{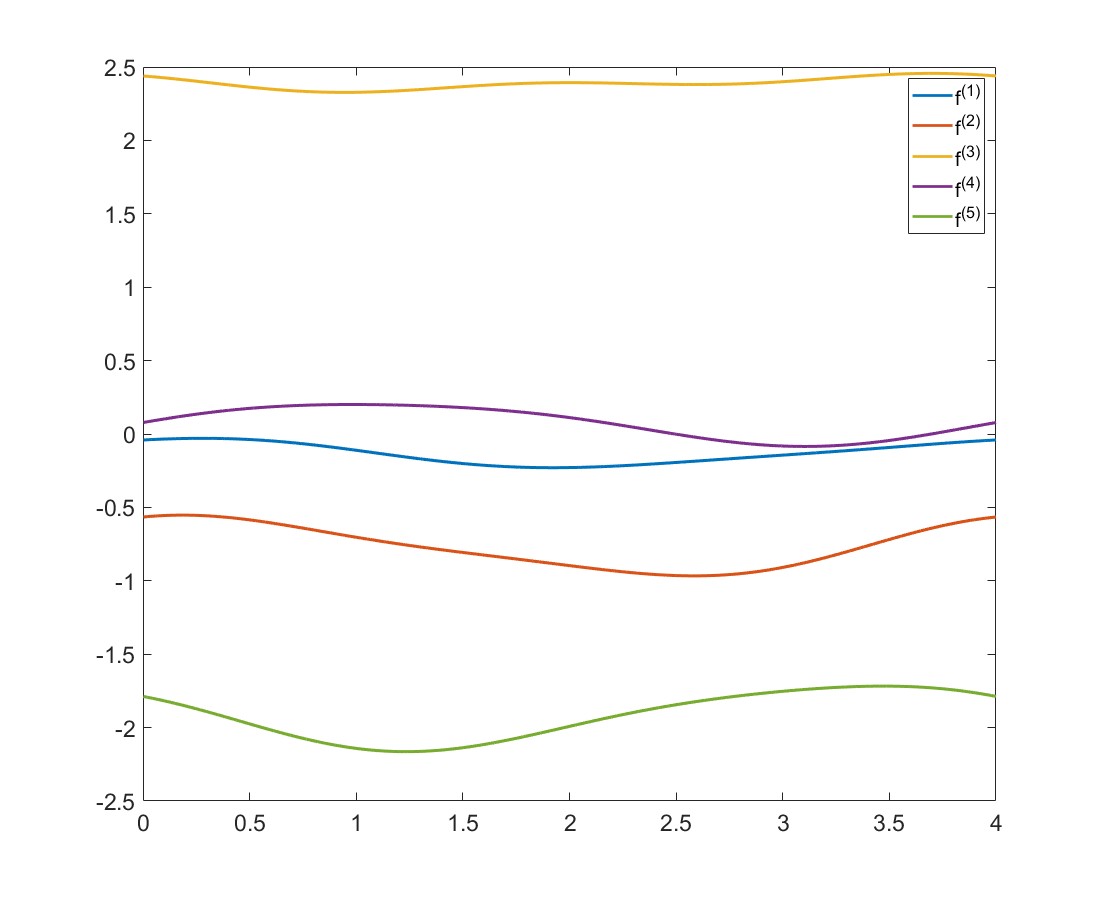}&
		\includegraphics[width=0.25\textwidth]{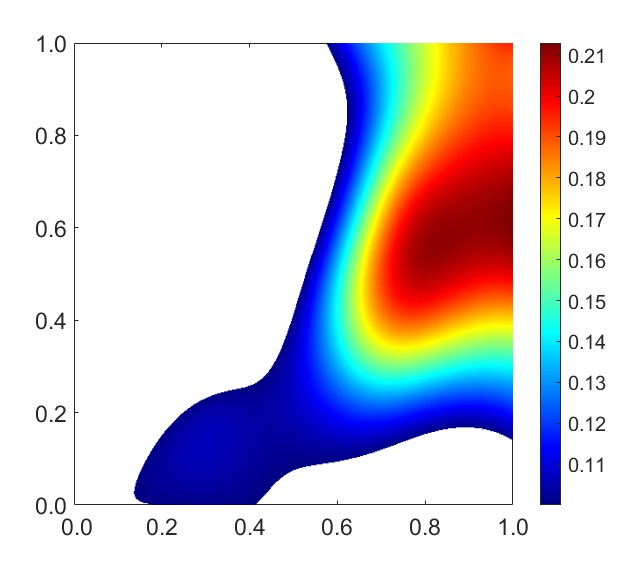}&
		\includegraphics[width=0.25\textwidth]{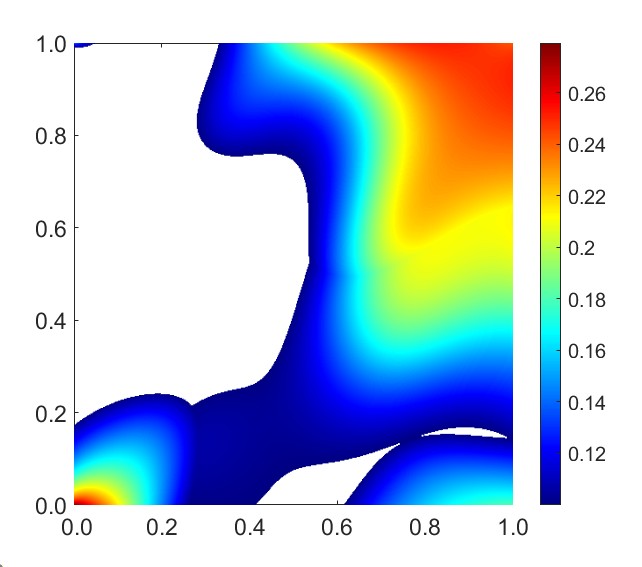}\\
		(a) Boundary $\fl$ & (b) $L=1$ & (c) $L=2$ \\
		\includegraphics[width=0.25\textwidth]{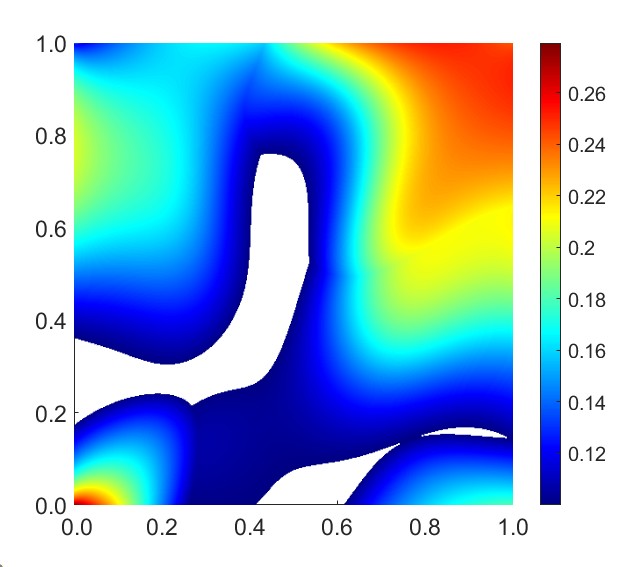}&
		\includegraphics[width=0.25\textwidth]{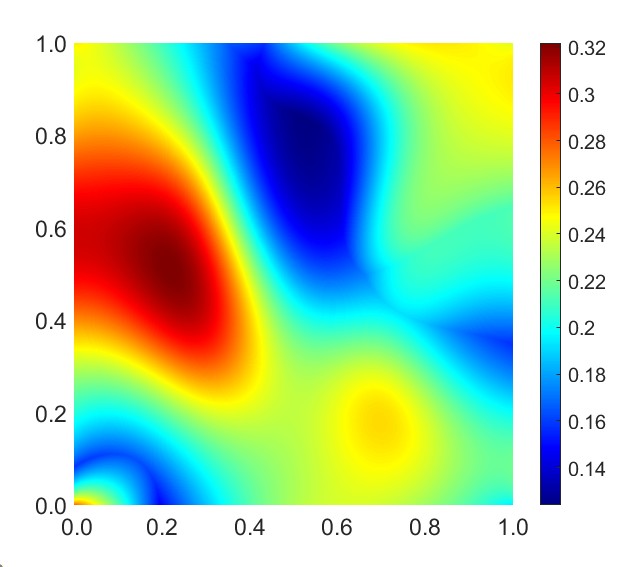}&
		\includegraphics[width=0.25\textwidth]{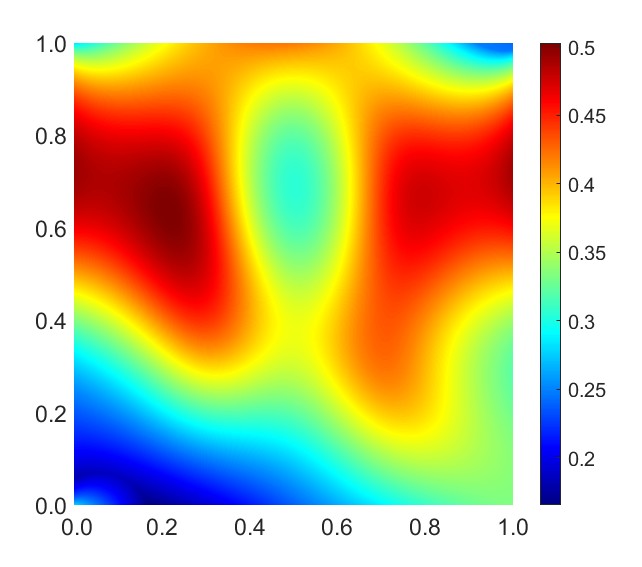}\\
		(d) $L=3$ & (e) $L=4$ & (f) $L=5$ \\
	\end{tabular}
	\caption{Boundary illuminations and the non-zero region of Example~\ref{ex:2}. Top left: plot of boundary data $\fl=g^{(\ell+1)}$. Top middle to bottom right: region which satisfying the non-zero condition as number of boundary input increasing.}
	\label{Fig:ex2_nonzero}
\end{figure}

\begin{figure}[htbp]
	\centering
	\begin{tabular}{ccc}
		\includegraphics[width=0.25\textwidth]{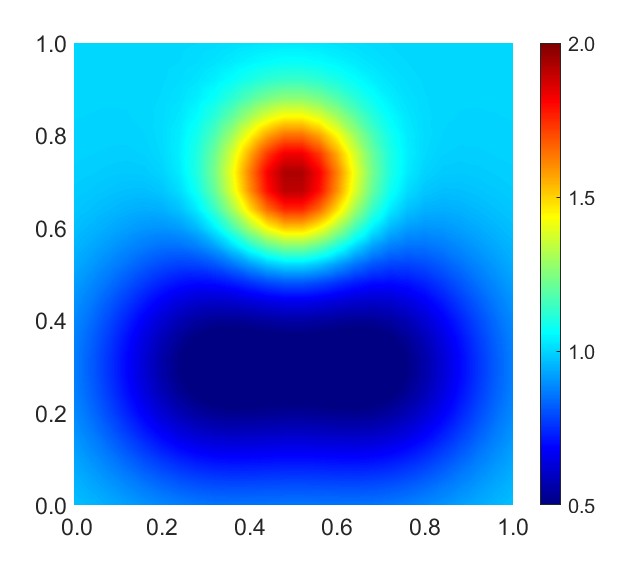}&
		\includegraphics[width=0.25\textwidth]{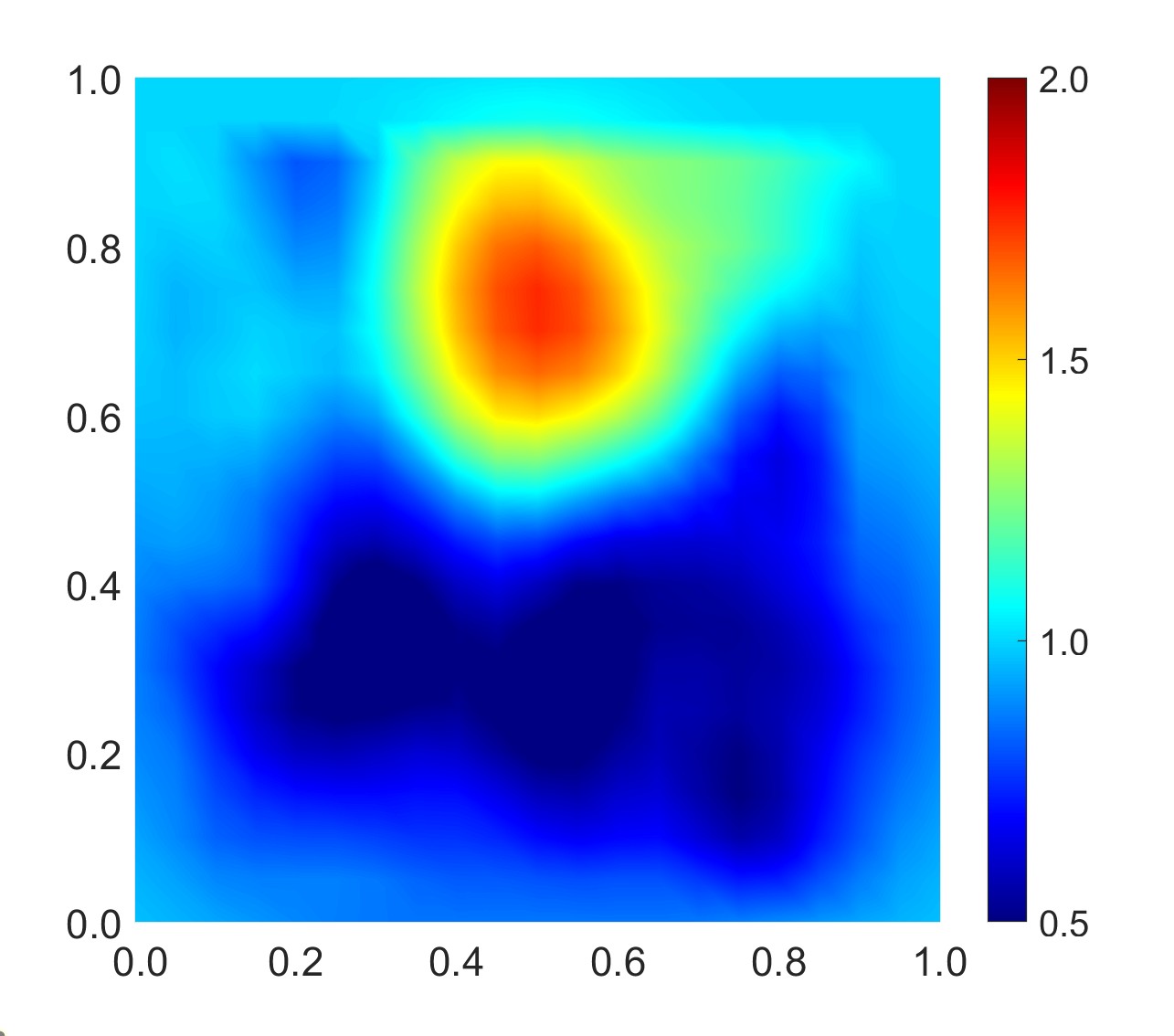}&
		\includegraphics[width=0.25\textwidth]{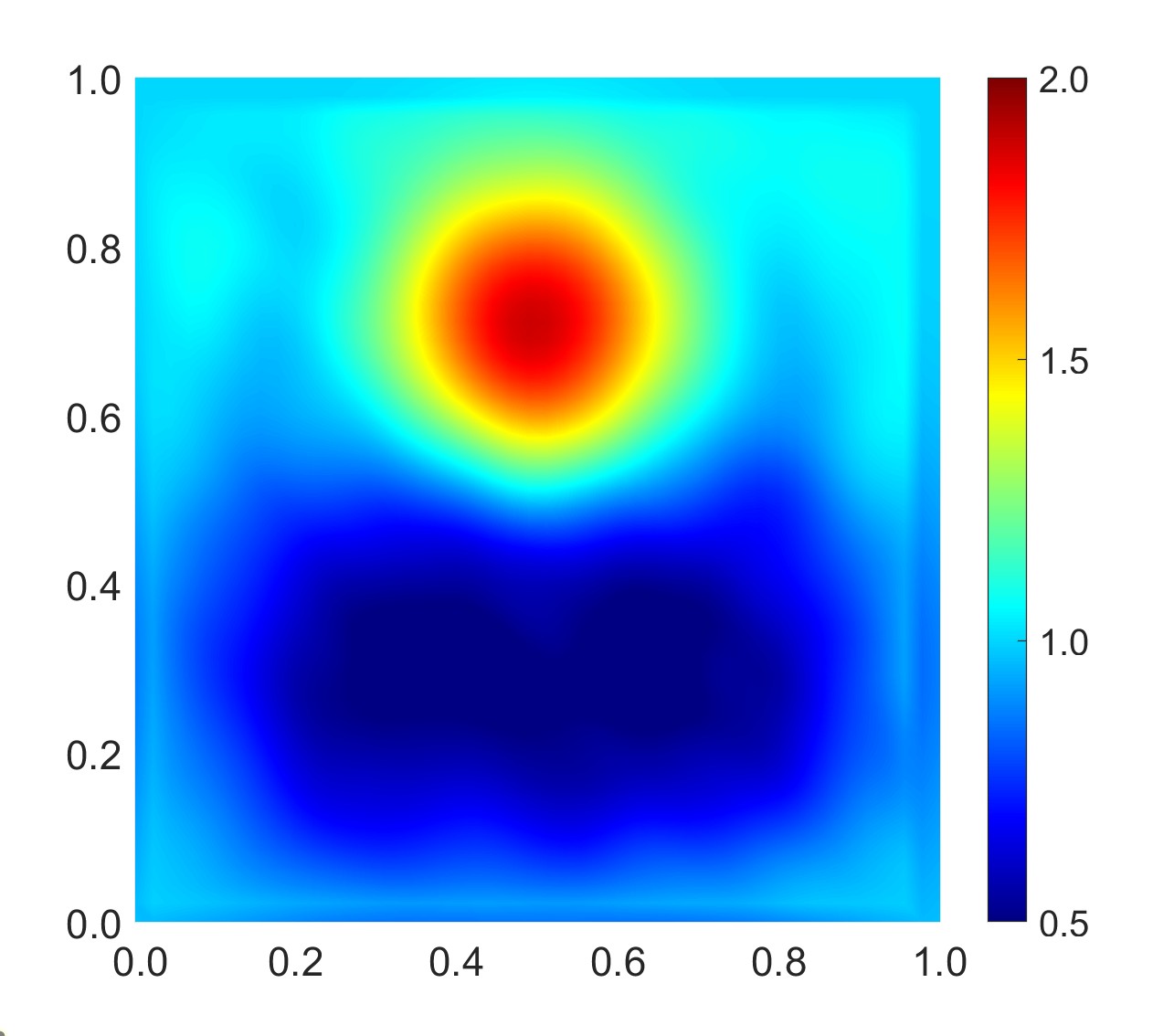}\\
		(a) $D^\dag$ & (b) $\delta=5e$-$2$ & (c) $\delta=1e$-$2$ \\
		\includegraphics[width=0.25\textwidth]{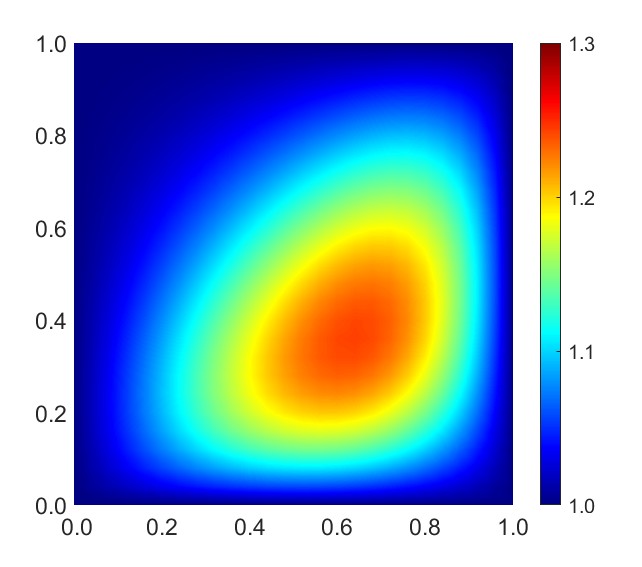}&
		\includegraphics[width=0.25\textwidth]{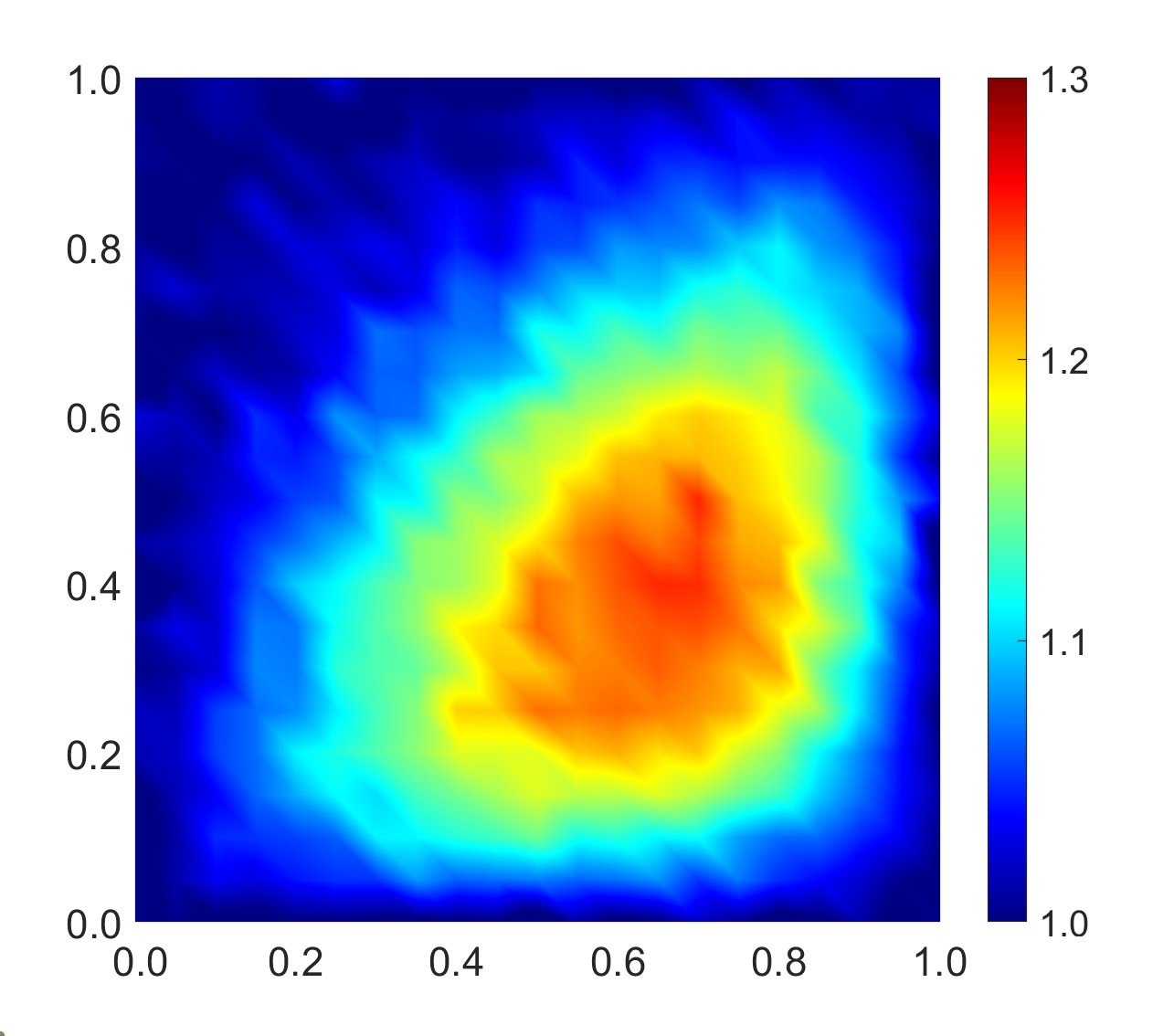}&
		\includegraphics[width=0.25\textwidth]{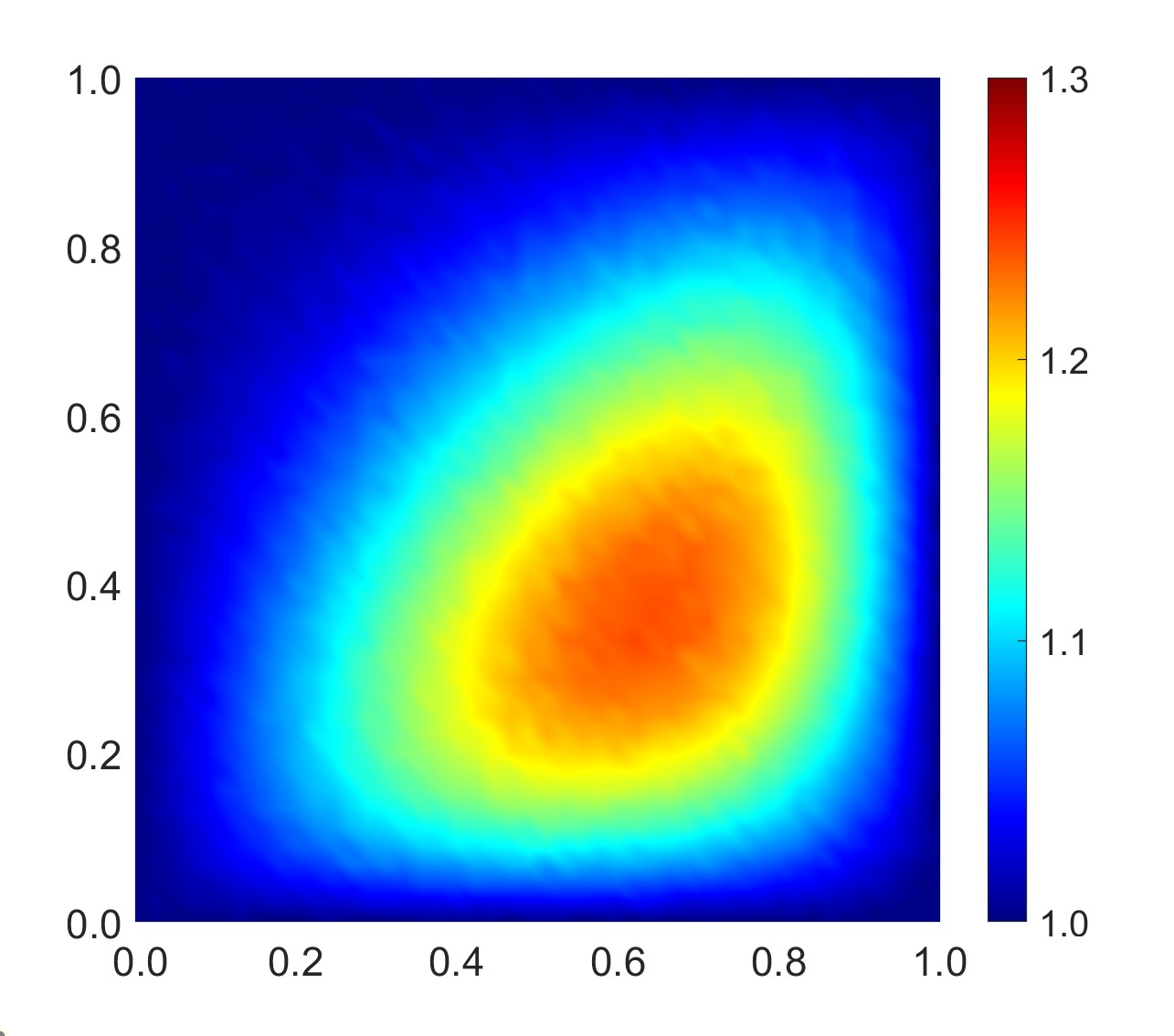}\\
		(d) $\sigma^\dag$ & (e) $\delta=5e$-$2$ & (f) $\delta=1e$-$2$ \\
	\end{tabular}
	\caption{Example~\ref{ex:2}. First row: reconstructions of $D^\dag$.
		Second row: reconstructions of $\sigma^\dag$. }
	\label{Fig:ex2_recon}
\end{figure}

\begin{example}\label{ex:2}
	$\Omega=(0,1)^2$,  $D^\dag=1+D_1-\frac{1}{2}D_2-\frac{1}{2}D_3$ with $D_1(x,y)=e^{-40(x-0.5)^2-40(y-0.7)^2} $, $D_2(x,y)=e^{-15(x-0.3)^2-15(y-0.3)^2}$, $D_3(x,y)=e^{-15(x-0.7)^2-15(y-0.3)^2}$  and the absorption coefficient $\sigma^{\dag}(x,y)=1+0.5\sin(\pi x)\sin(\pi y)e^{-4(1-x)y}  $.% The exact solution $\ul$ are computed  with mesh size $h=1/500$.
\end{example}
The region representing the non-zero condition and the numerical reconstructions of Example~\ref{ex:2} are shown in Figures~\ref{Fig:ex2_nonzero}-\ref{Fig:ex2_recon} and Table~\ref{tab:ex2}. For the non-zero condition region, we observe a similar behavior as in Example~\ref{ex:1}, the region enlarges with the addition of boundary illuminations. For testing the convergence rates of reconstruction errors, we fix $L=5$ and initially choose the mesh size $h=1/12$ and the regularization parameter $\alpha=5e$-$7$. We observe the convergence rate $O(\delta^{0.35})$ for $e_D$ and $O(\delta^{0.42})$ for $e_\sigma$. The experimental convergence rates are slightly higher than the theoretical rate $O(\delta^{0.25})$. {This discrepancy is acceptable, since our theoretical bound is optimal with respect to the stability estimate in Theorem \ref{thm:stability}. Moreover, improved convergence rates may be possible when the parameter and solution possess higher regularity.} Furthermore,  the first step of the reconstruction algorithm requires to solve an optimization problem to get $q_h^*$, the non-convexity of the loss function may lead to local minima, making it challenging to verify the theoretical convergence rates.  
Figure~\ref{Fig:ex2_recon} shows the reconstructions with $5\%$ and $1\%$ noise level, with $h=1/20$, $\alpha=5e$-$7$  and  $h=1/45$, $\alpha=2e$-$8$, respectively.

\begin{table}[htp!]
	\centering
	\caption{The convergence rates for Example~\ref{ex:2} with respect to $\delta$.} \label{tab:ex2}
	\begin{tabular}{c|cccccccc}
		\toprule
		$\delta$   & 1e-2 & 5e-3 & 2e-3 & 1e-3 & 5e-4 & 2e-4  &1e-4   & rate \\
		\midrule
		$e_D $ &6.34e-2 & 5.72e-2 & 3.47e-2 & 2.65e-2 & 2.40e-2  & 1.85e-2 &  1.24e-2 & {$O(\delta^{0.35})$}  \\
		$e_{\sigma} $ & 1.04e-2  & 5.67e-3 & 4.08e-3 & 2.95e-3 &  2.70e-3  & 1.84e-3  &  1.24e-3 & {$O(\delta^{0.42})$} \\		
		\bottomrule
	\end{tabular}
\end{table}

\begin{figure}[htbp]
	\centering
	\begin{tabular}{ccc}
		\includegraphics[width=0.25\textwidth]{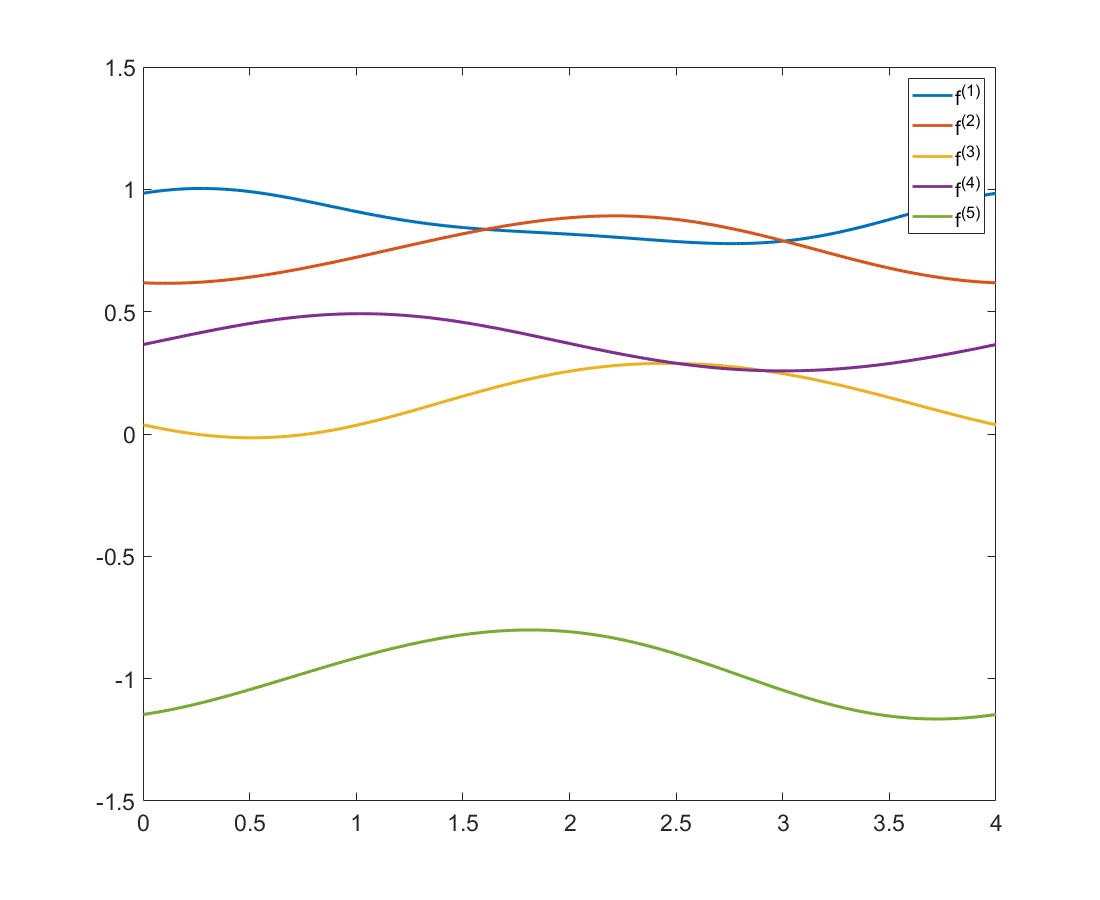}&
		\includegraphics[width=0.25\textwidth]{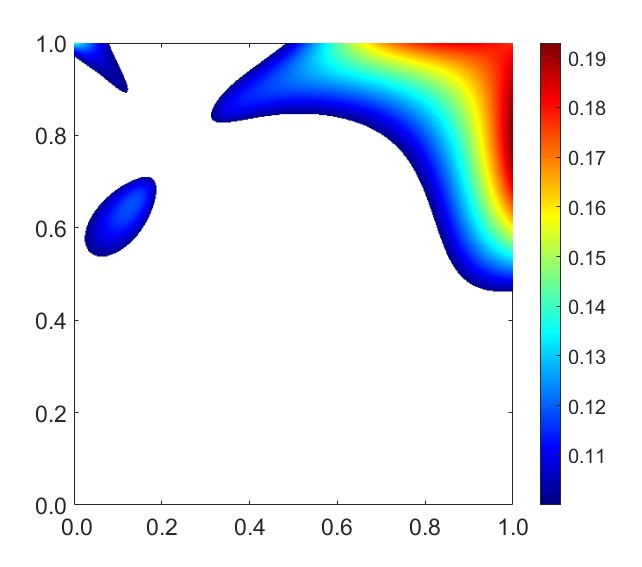}&
		\includegraphics[width=0.25\textwidth]{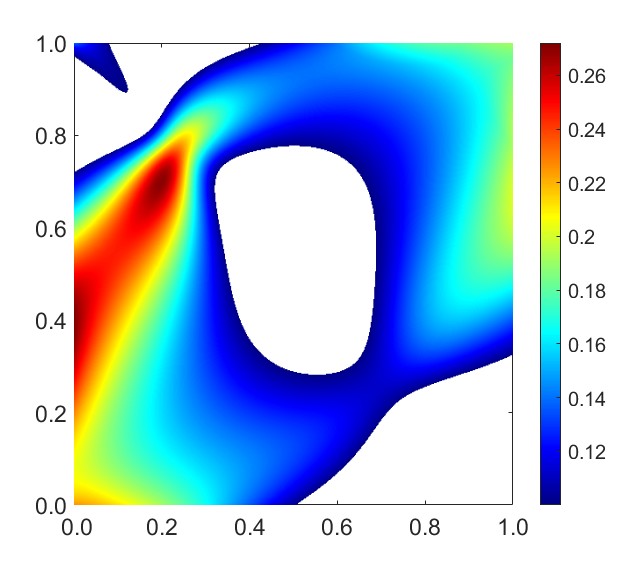}\\
		(a) Boundary $\fl$ & (b) $L=1$ & (c) $L=2$ \\
		\includegraphics[width=0.25\textwidth]{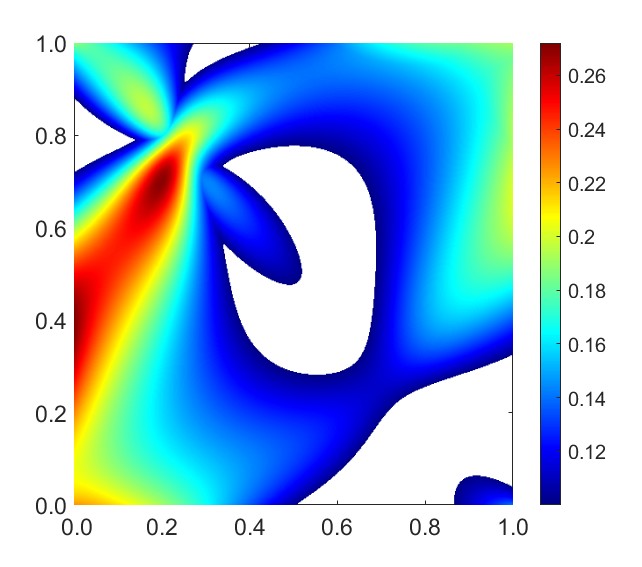}&
		\includegraphics[width=0.25\textwidth]{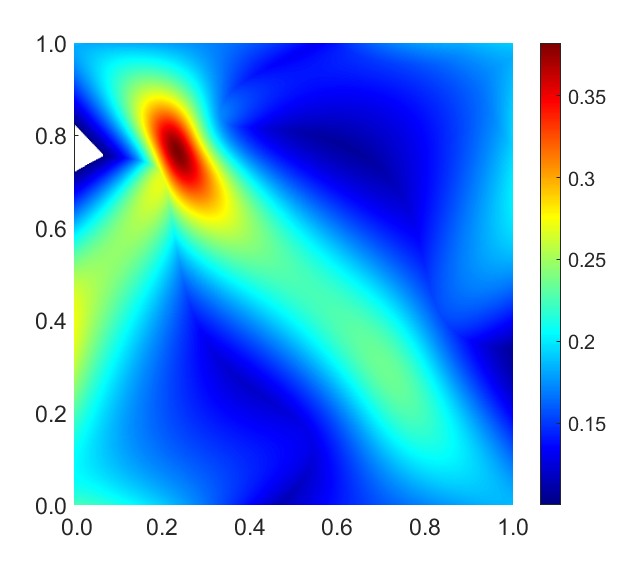}&
		\includegraphics[width=0.25\textwidth]{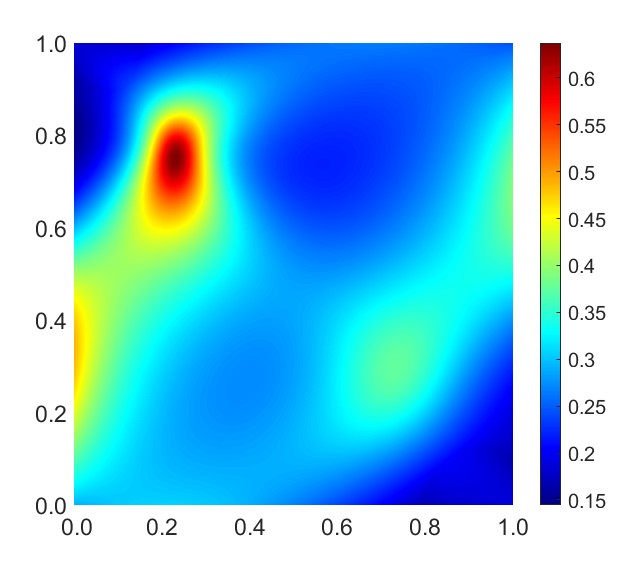}\\
		(d) $L=3$ & (e) $L=4$ & (f) $L=5$ \\
	\end{tabular}
	\caption{Boundary illuminations and the non-zero region of Example~\ref{ex:3}. Top left: plot of boundary data $\fl=g^{(\ell+1)}$. Top middle to bottom right: region which satisfying the non-zero condition as number of boundary input increasing.}
	\label{Fig:ex3_nonzero}
\end{figure}

\begin{figure}[htbp]
	\centering
	\begin{tabular}{ccc}
		\includegraphics[width=0.25\textwidth]{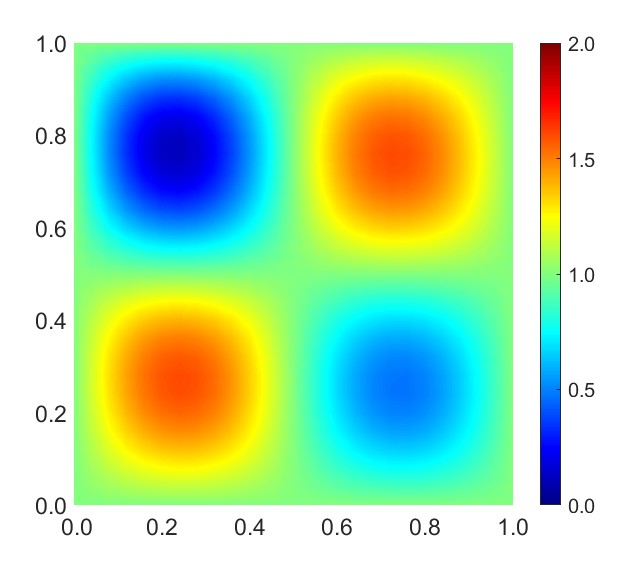}&
		\includegraphics[width=0.25\textwidth]{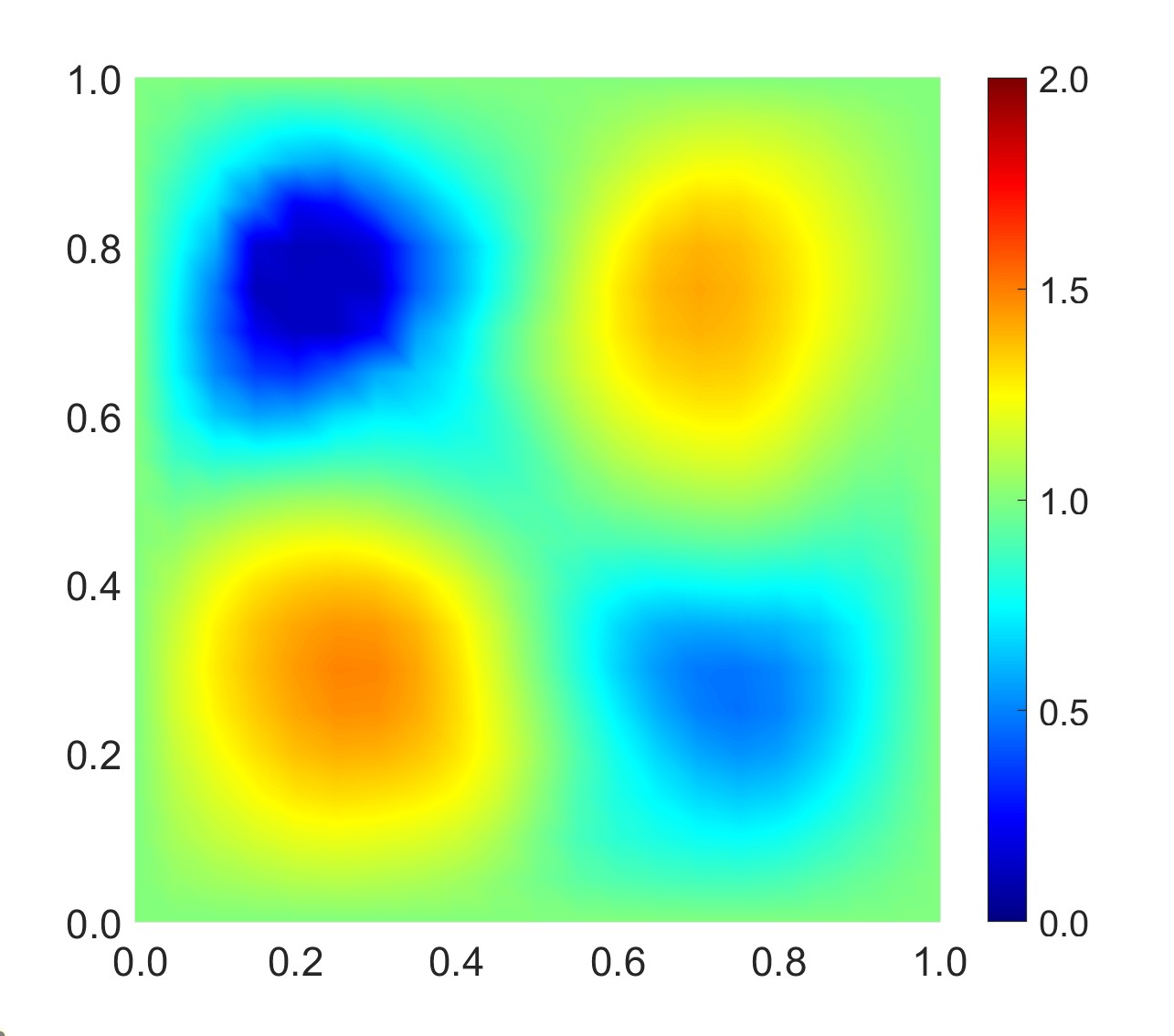}&
		\includegraphics[width=0.25\textwidth]{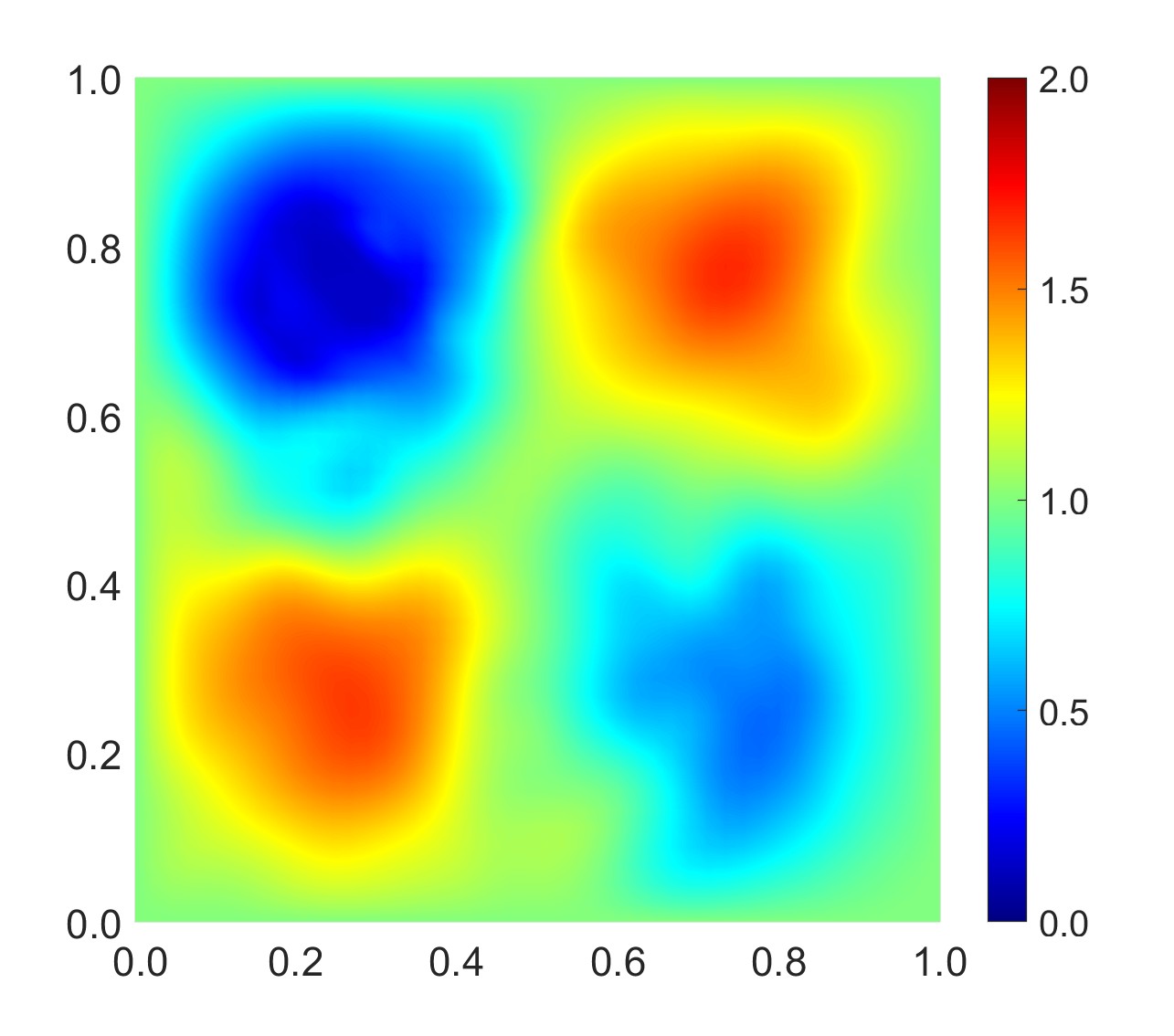}\\
		(a) $D^\dag$ & (b) $\delta=5e$-$2$ & (c) $\delta=1e$-$2$ \\
		\includegraphics[width=0.25\textwidth]{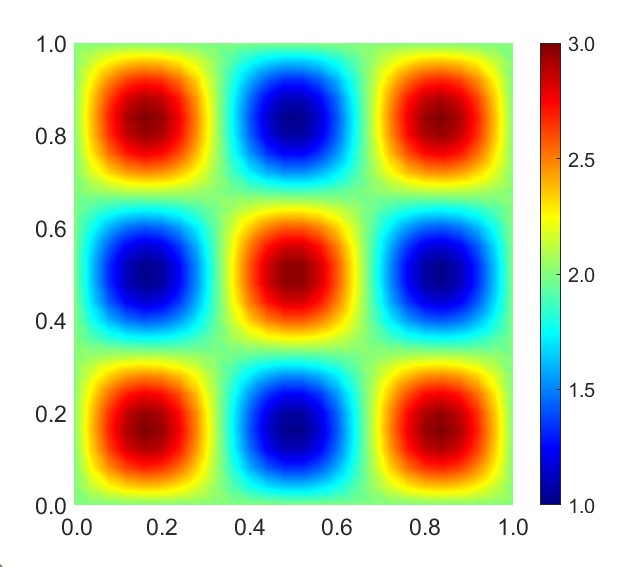}&
		\includegraphics[width=0.25\textwidth]{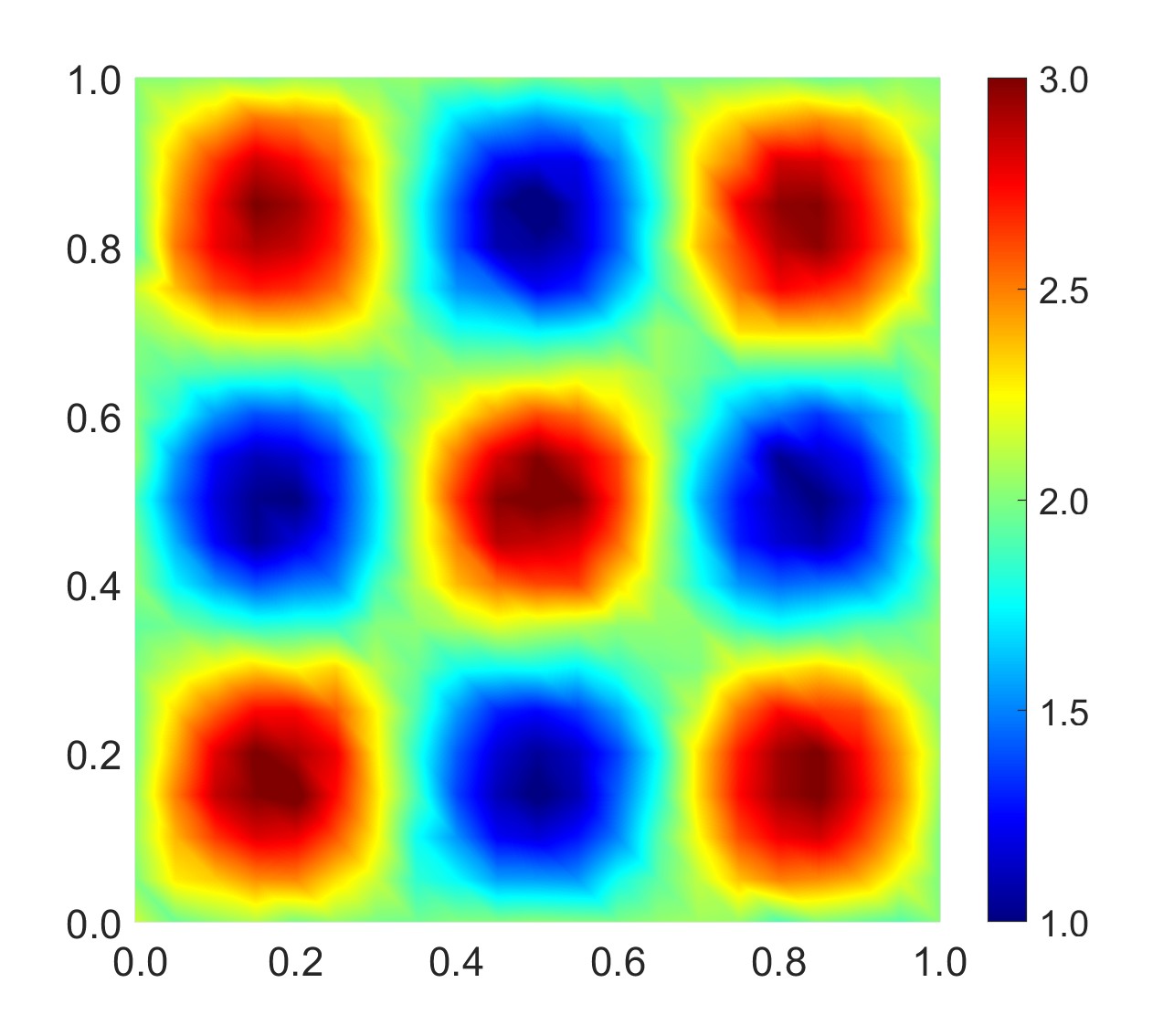}&
		\includegraphics[width=0.25\textwidth]{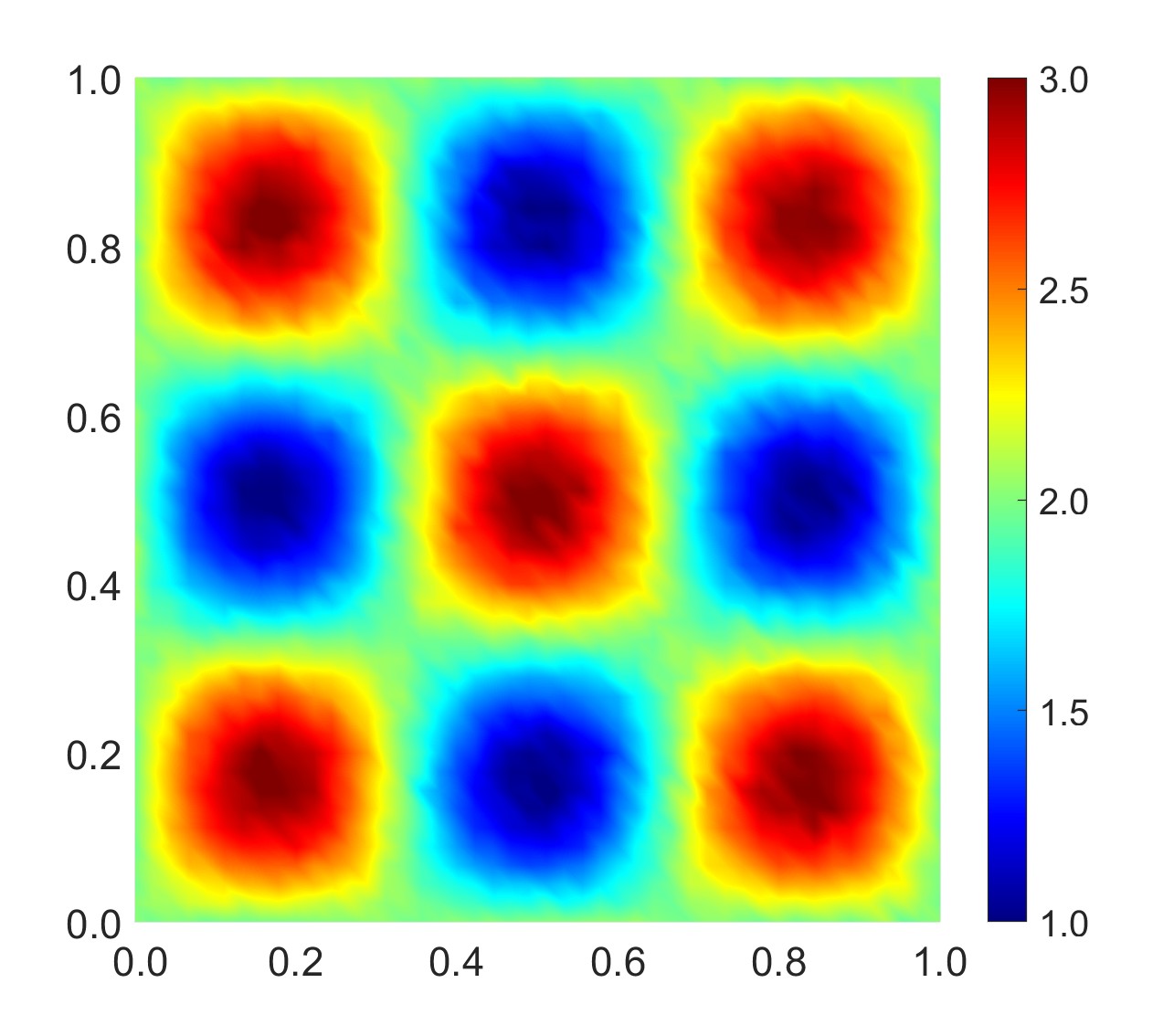}\\
		(d) $\sigma^\dag$ & (e) $\delta=5e$-$2$ & (f) $\delta=1e$-$2$ \\
	\end{tabular}
	\caption{Example~\ref{ex:3}. First row: reconstructions of $D^\dag$.
		Second row: reconstructions of $\sigma^\dag$. }
	\label{Fig:ex3_recon}
\end{figure}

\begin{example}\label{ex:3}
	$\Omega=(0,1)^2$,  $D^{\dag}(x,y)=1+\frac{1}{2}\sin(2\pi x)\sin(2\pi y)e^{xy}$ and $\sigma^{\dag}(x,y)=3+\sin(3\pi x)\sin(3\pi y)$.  
\end{example}
In this example, we consider a  more challenging setting. The  absorption coefficient $\sigma^\dag$ has high oscillations. Figure~\ref{Fig:ex3_nonzero} shows the behavior of the non-zero condition. The non-zero condition is satisfied in the whole domain when sufficiently many random boundary illuminations are used.  Table~\ref{tab:ex3} present the convergence rates. Here, we fix $L=5$ and choose the initial mesh size $h=1/16$ and the regularization parameter $\alpha=2e$-$6$. The convergence rate for $e_D$ is $O(\delta^{0.26})$, which aligns with the predicted rate $O(\delta^{0.25})$. However, we observe a much faster decay for $e_\sigma$, with convergence rate $O(\delta^{0.39})$. Figure~\ref{Fig:ex3_recon} demonstrates that  even for this challenging absorption coefficient, the reconstruction is accurate for high noise levels.  
Here we take $h=1/20$, $\alpha=2e$-$6$ for noise level $\delta=5e$-$2$ and  $h=1/45$, $\alpha=8e$-$8$ for  $\delta=1e$-$2$.

\begin{table}[htp!]
	\centering
	\caption{The convergence rates for Example~\ref{ex:3} with respect to $\delta$.} \label{tab:ex3}
	\begin{tabular}{c|cccccccc}
		\toprule
		$\delta$   & 1e-2 & 5e-3 & 2e-3 & 1e-3 & 5e-4 & 2e-4  &1e-4   & rate \\
		\midrule
		$e_D $ & 7.80e-2 & 5.78e-2 & 3.53e-2 & 3.13e-2 & 2.82e-2  & 2.78e-2 & 2.15e-2 & {$O(\delta^{0.26})$}  \\
		$e_{\sigma} $ & 1.36e-2  & 7.73e-3 & 3.24e-3 & 3.06e-3 &  2.54e-3  & 2.42e-3  &  1.89e-3 & {$O(\delta^{0.39})$} \\		
		\bottomrule
	\end{tabular}
\end{table}

\begin{figure}[htbp]
 	\centering
 	\begin{tabular}{ccc}
 		\includegraphics[width=0.25\textwidth]{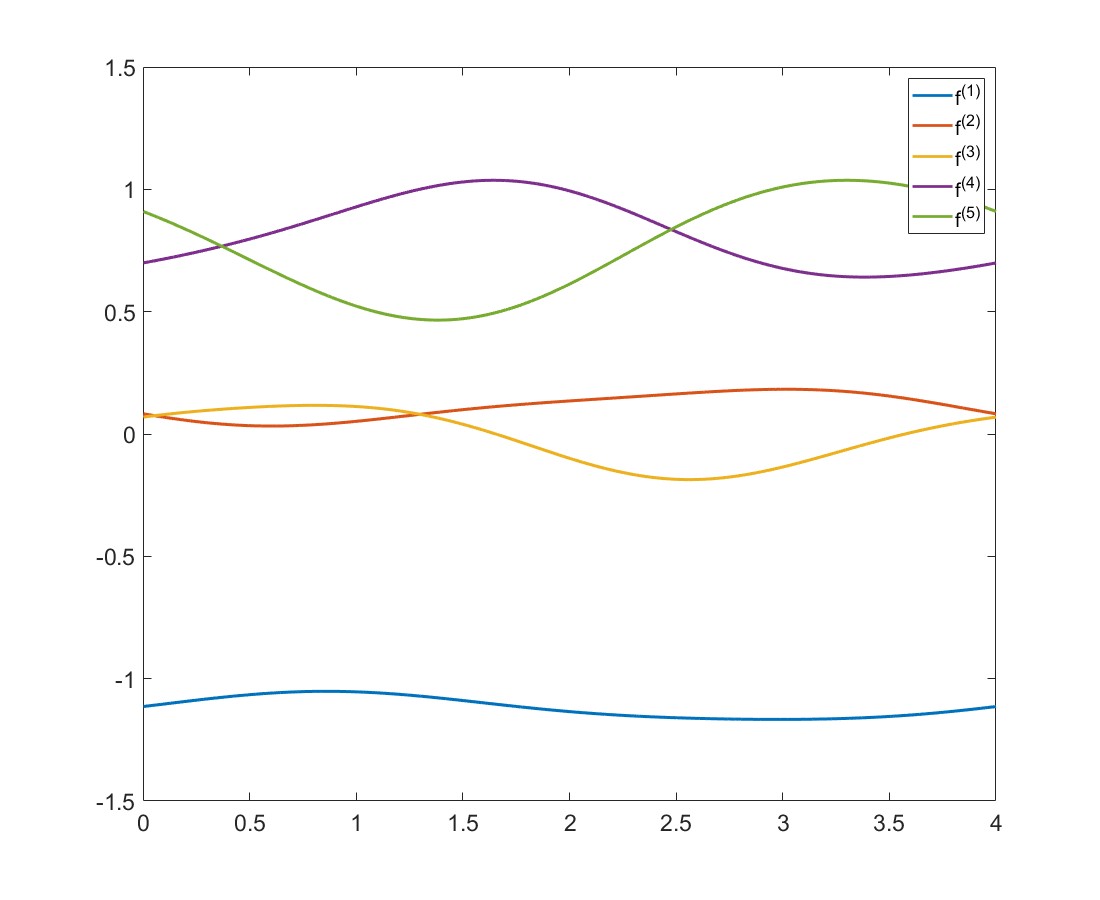}&
 		\includegraphics[width=0.25\textwidth]{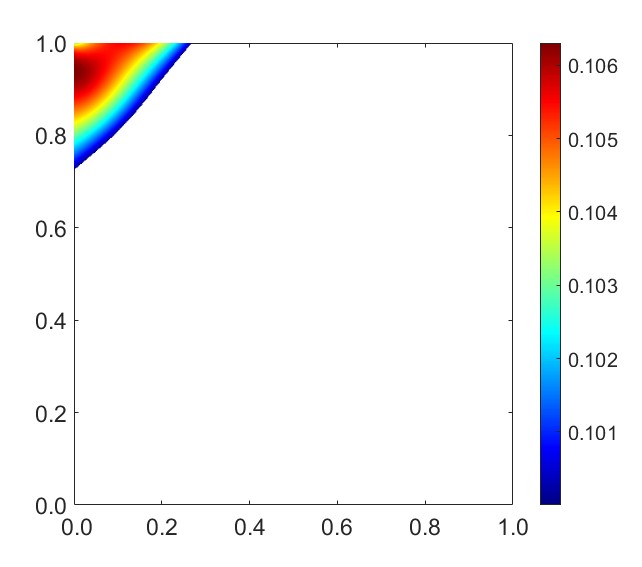}&
 		\includegraphics[width=0.25\textwidth]{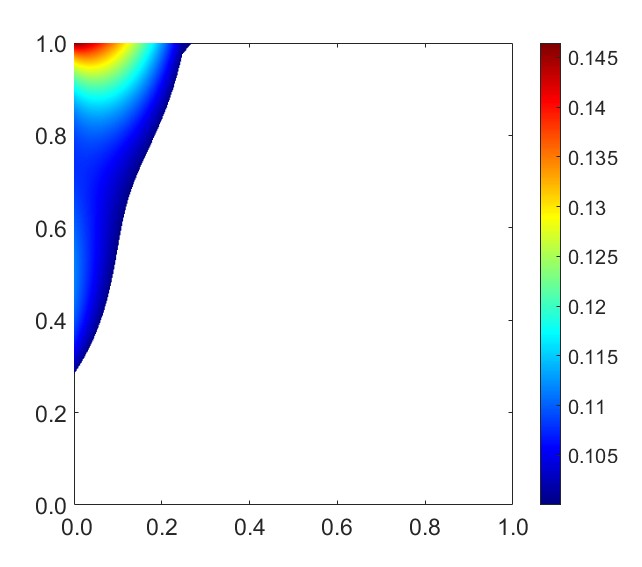}\\
 		(a) Boundary $\fl$ & (b) $L=1$ & (c) $L=2$ \\
 		\includegraphics[width=0.25\textwidth]{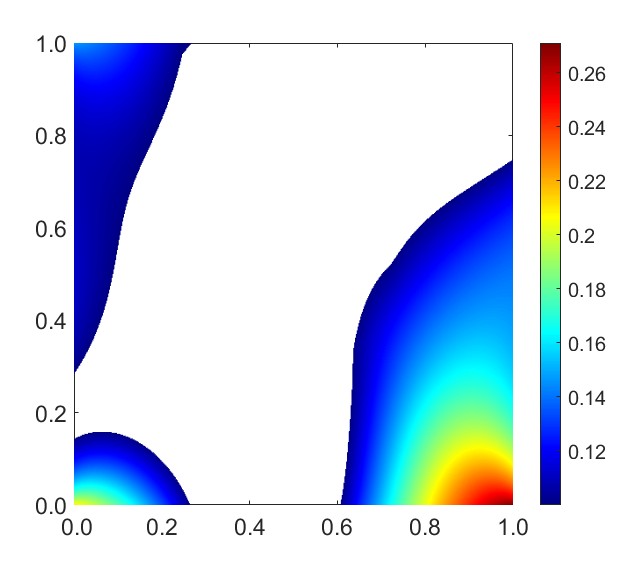}&
 		\includegraphics[width=0.25\textwidth]{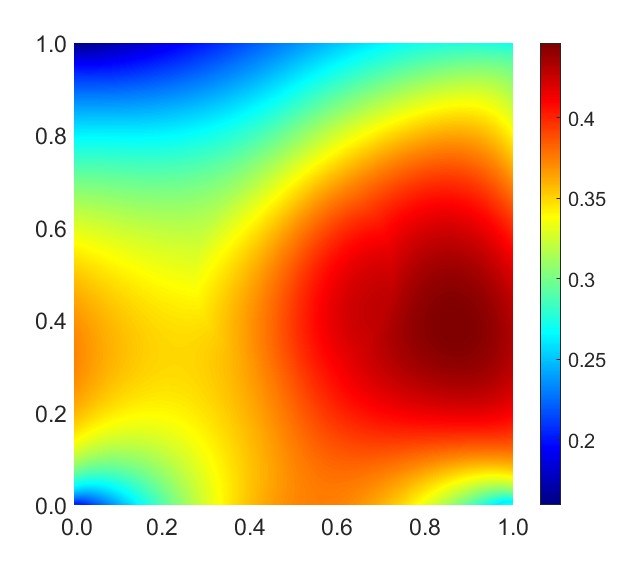}&
 		\includegraphics[width=0.25\textwidth]{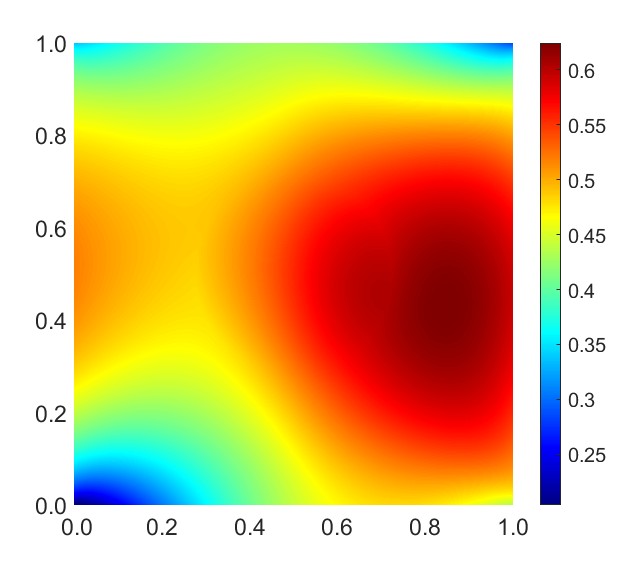}\\
 		(d) $L=3$ & (e) $L=4$ & (f) $L=5$ \\
 	\end{tabular}
 	\caption{Boundary illuminations and the non-zero region of Example~\ref{ex:4}. Top left: plot of boundary data $\fl=g^{(\ell+1)}$. Top middle to bottom right: region which satisfying the non-zero condition as number of boundary input increasing.}
 	\label{Fig:ex4_nonzero}
 \end{figure}

\begin{figure}[htbp]
 	\centering
 	\begin{tabular}{ccc}
 		\includegraphics[width=0.25\textwidth]{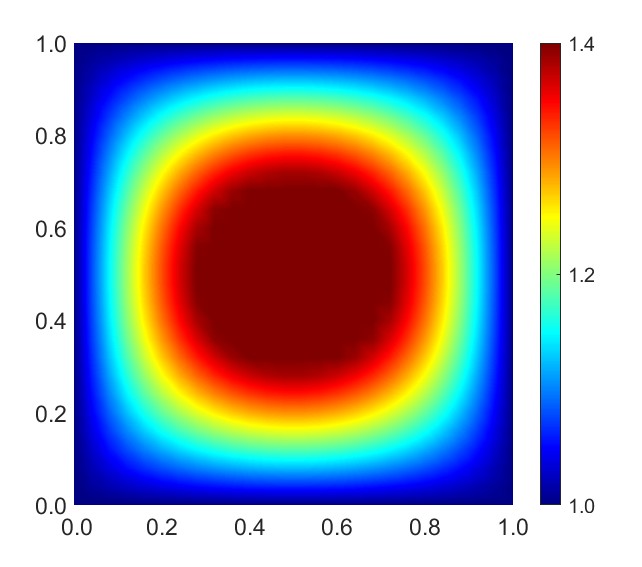}&
 		\includegraphics[width=0.25\textwidth]{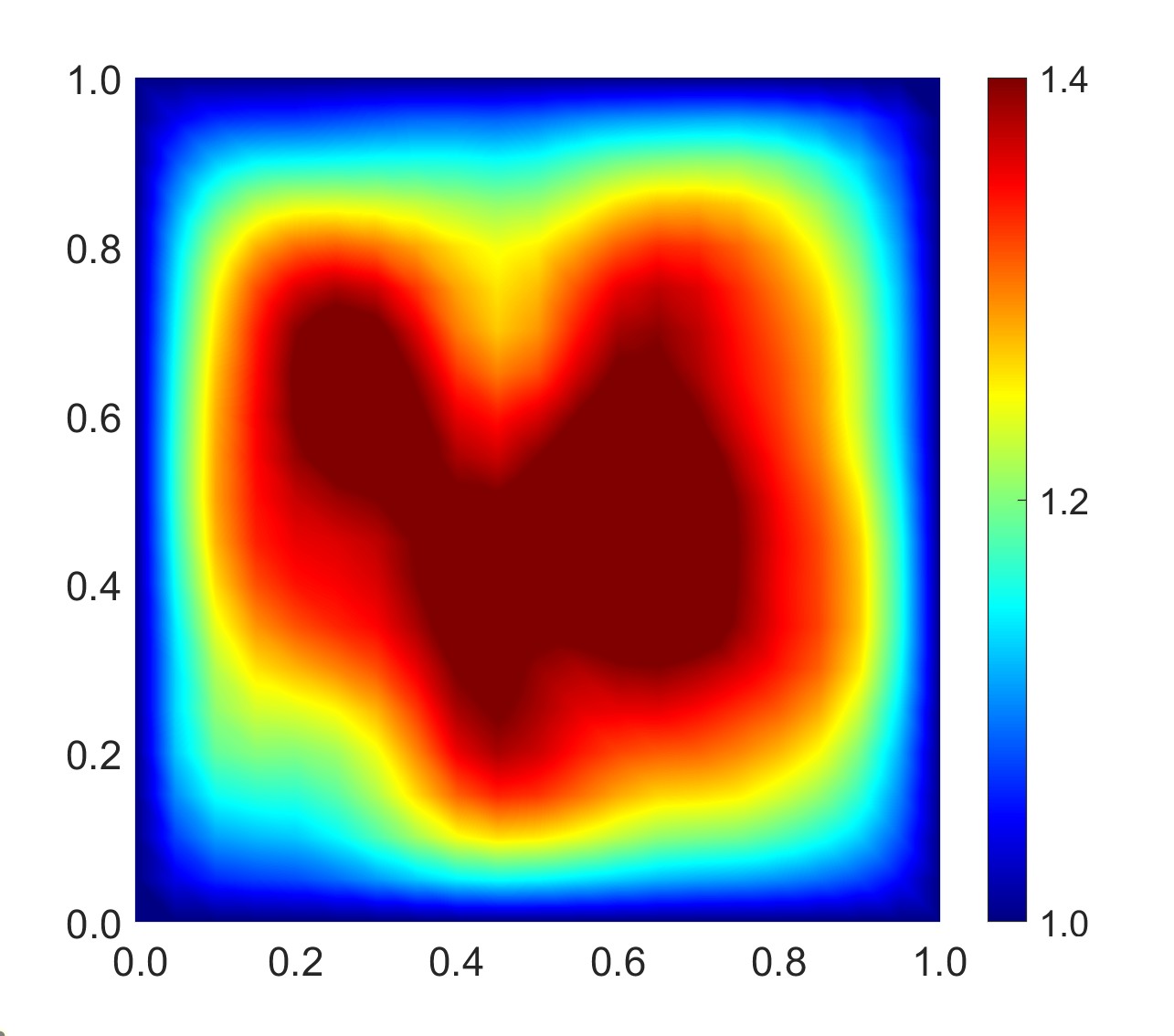}&
 		\includegraphics[width=0.25\textwidth]{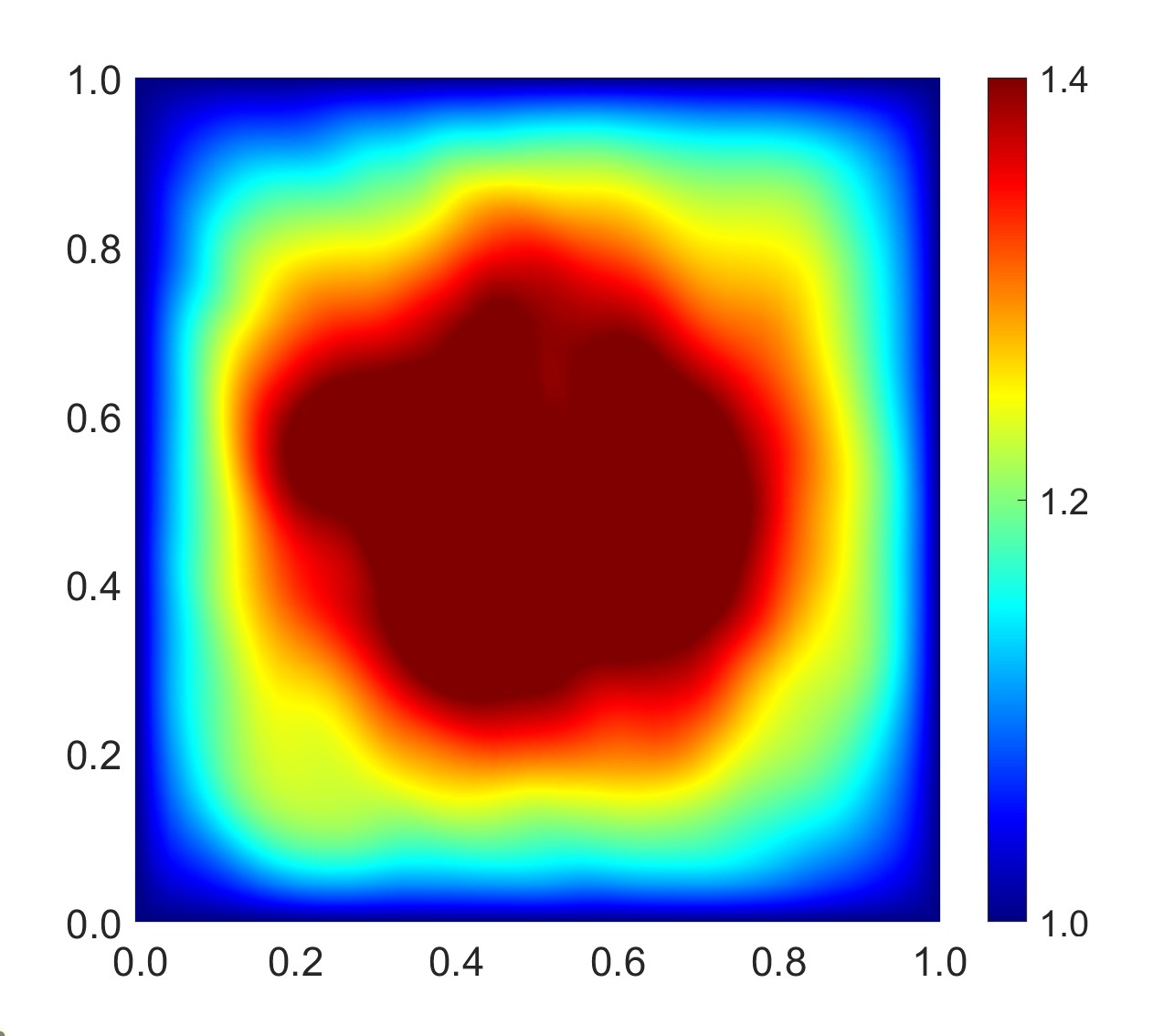}\\
 		(a) $D^\dag$ & (b) $\delta=5e$-$2$ & (c) $\delta=1e$-$2$ \\
 		\includegraphics[width=0.25\textwidth]{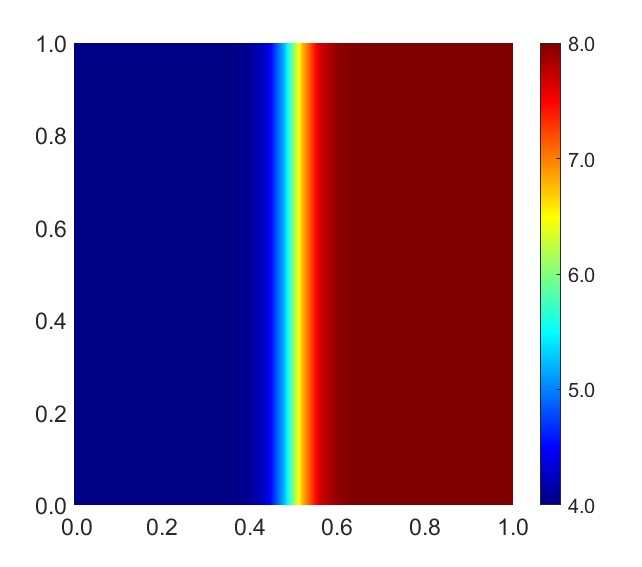}&
 		\includegraphics[width=0.25\textwidth]{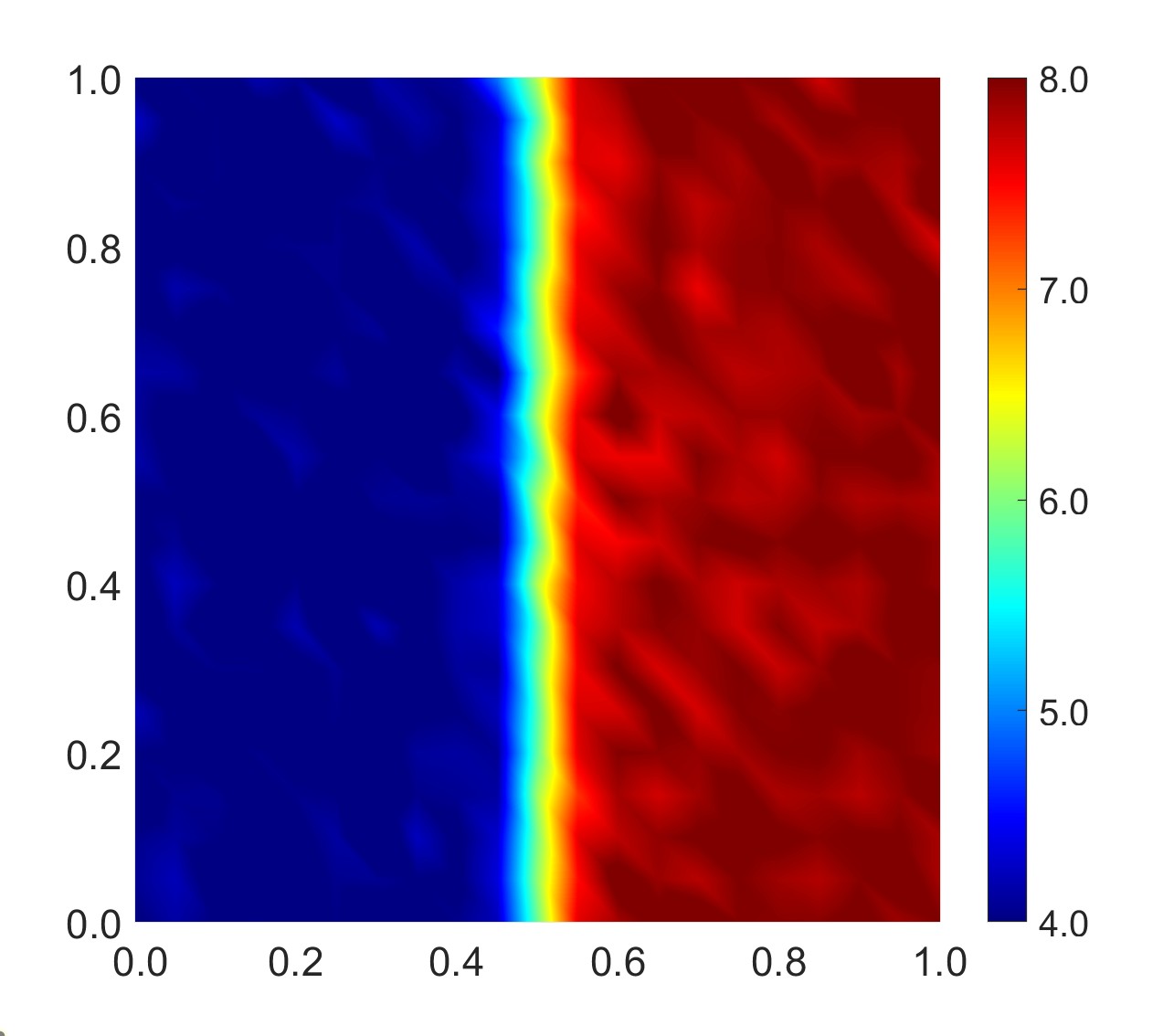}&
 		\includegraphics[width=0.25\textwidth]{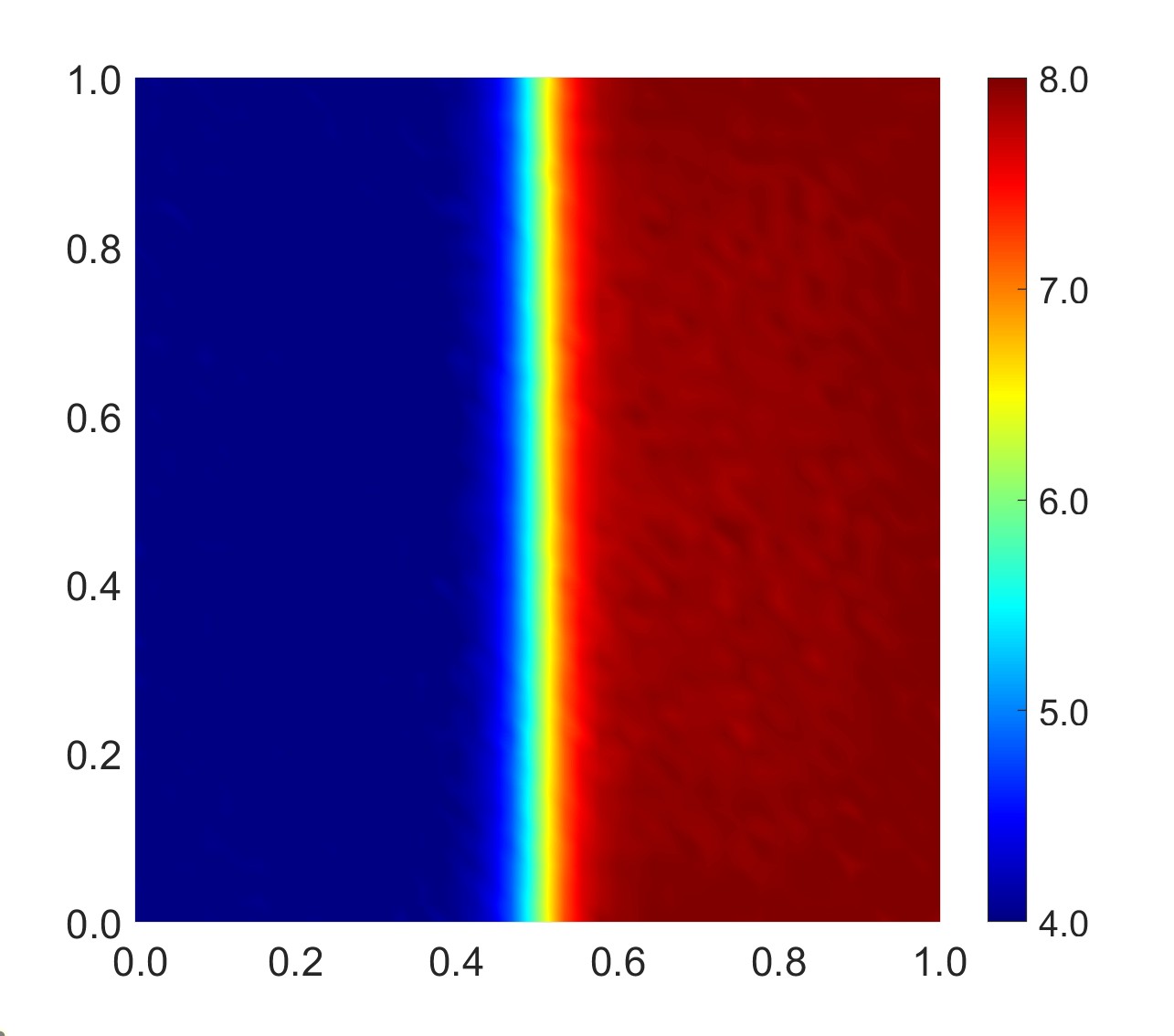}\\
 		(d) $\sigma^\dag$ & (e) $\delta=5e$-$2$ & (f) $\delta=1e$-$2$ \\
 	\end{tabular}
 	\caption{Example~\ref{ex:4}. First row: reconstructions of $D^\dag$.
 		Second row: reconstructions of $\sigma^\dag$. }
 	\label{Fig:ex4_recon}
 \end{figure}

Next, we present numerical results for nonsmooth coefficients.
 \begin{example}\label{ex:4}
 $\Omega=(0,1)^2$,  $D^{\dag}(x,y)=\min(1.4,1+2x(1-x)\sin(\pi y))$ and $\sigma^{\dag}(x,y)=6+2\tanh(20x-10)$.  
 \end{example}
Here, we cut off the diffusion coefficient $D^\dag$ in order to have discontinuous derivatives.  Additionally, the absorption coefficient $\sigma^\dag$ includes a sharp interface where the magnitudes of the derivatives are large. The non-zero condition and the numerical reconstructions are presented in Figures~\ref{Fig:ex4_nonzero}-\ref{Fig:ex4_recon} and Table~\ref{tab:ex4}. The mesh size and the regularization parameter are initialized as $h=1/12$ and $\alpha=1e$-$5$. For this nonsmooth case, we still observe the convergence rates $O(\delta^{0.29})$ and $O(\delta^{0.33})$ for $e_D$ and $e_\sigma$, respectively. The convergence rate for the diffusion coefficient $D^\dag$ matches the predicted rate, whereas the convergence rate for absorption coefficient $\sigma^\dag$ is slightly higher. In the numerical reconstructions Figure~\ref{Fig:ex4_recon}, we take  $h=1/20$, $\alpha=5e$-$6$  and  $h=1/45$, $\alpha=2e$-$7$ for noise level $\delta=$5e-2 and $\delta=$1e-2, respectively.

 \begin{table}[htp!]
 	\centering
 	\caption{The convergence rates for Example~\ref{ex:4} with respect to $\delta$.} \label{tab:ex4}
 	\begin{tabular}{c|cccccccc}
 		\toprule
 		$\delta$   & 1e-2 & 5e-3 & 2e-3 & 1e-3 & 5e-4 & 2e-4  &1e-4   & rate \\
 		\midrule
 		$e_D $ & 4.89e-2 & 4.72e-2 & 3.39e-2 & 2.68e-2 & 2.14e-2  & 1.86e-2 & 1.31e-2 & {$O(\delta^{0.29})$}  \\
 		$e_{\sigma} $ & 1.65e-2  & 1.25e-2 & 7.89e-3 & 6.23e-3 &  5.17e-3  & 4.61e-3  &  3.33e-3 & {$O(\delta^{0.33})$} \\		
 		\bottomrule
 	\end{tabular}
 \end{table}

\begin{figure}[htbp]
	\centering
	\begin{tabular}{ccc}
		\includegraphics[width=0.25\textwidth]{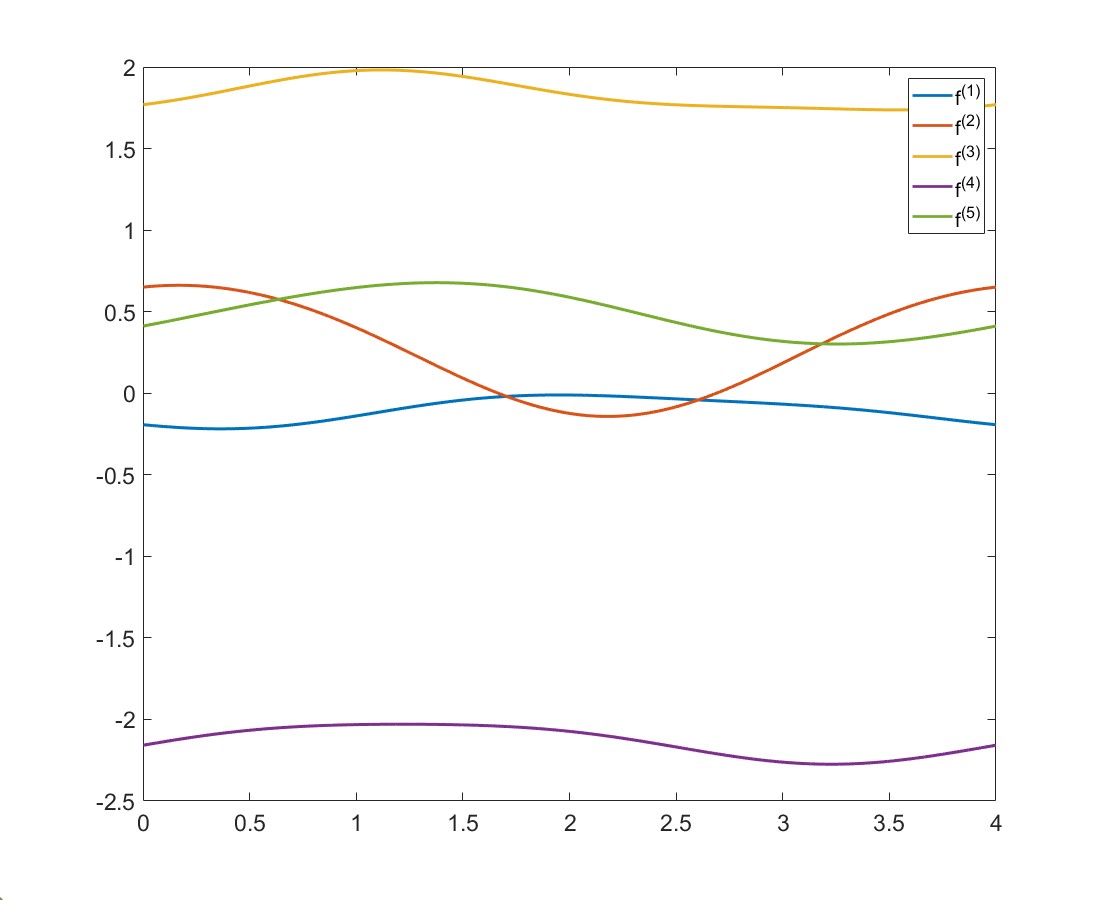}&
		\includegraphics[width=0.25\textwidth]{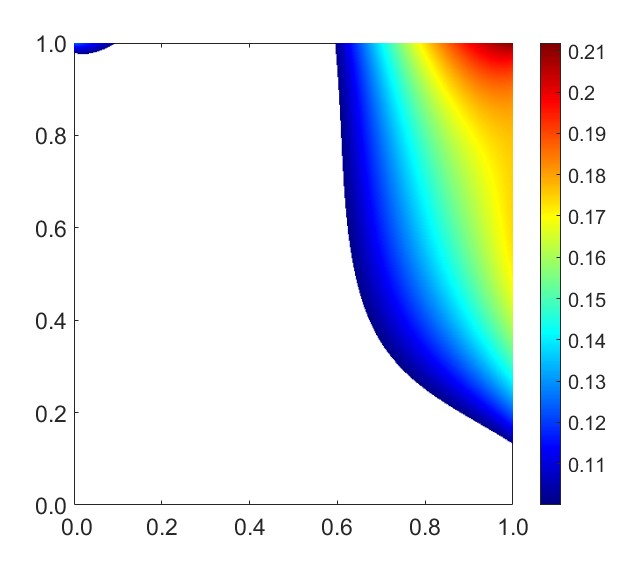}&
		\includegraphics[width=0.25\textwidth]{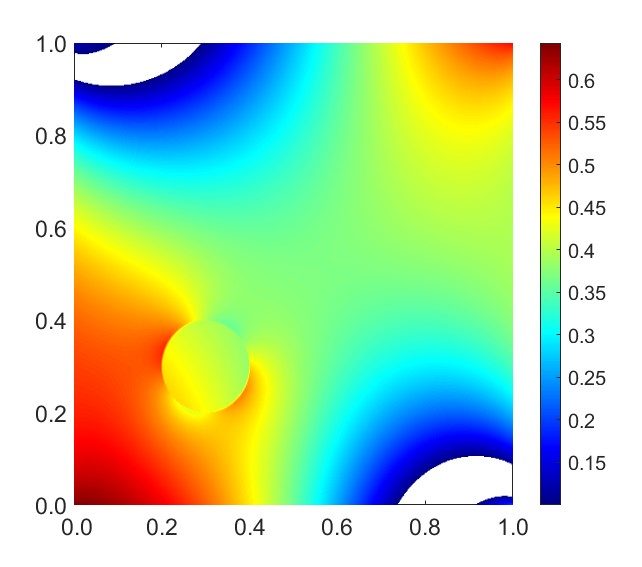}\\
		(a) Boundary $\fl$ & (b) $L=1$ & (c) $L=2$ \\
		\includegraphics[width=0.25\textwidth]{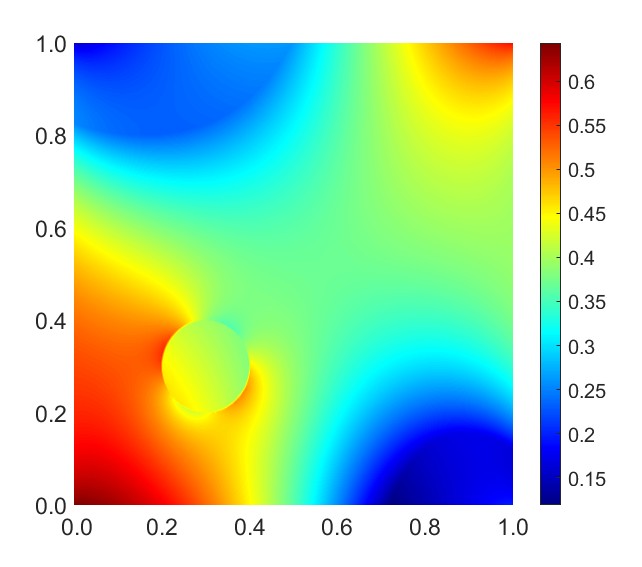}&
		\includegraphics[width=0.25\textwidth]{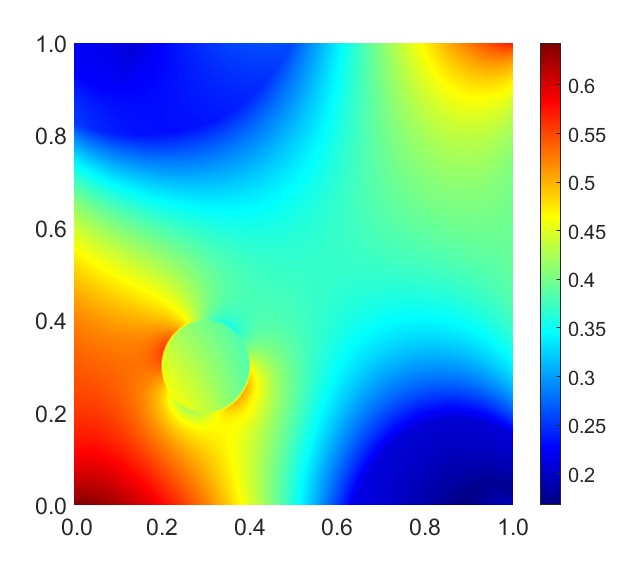}&
		\includegraphics[width=0.25\textwidth]{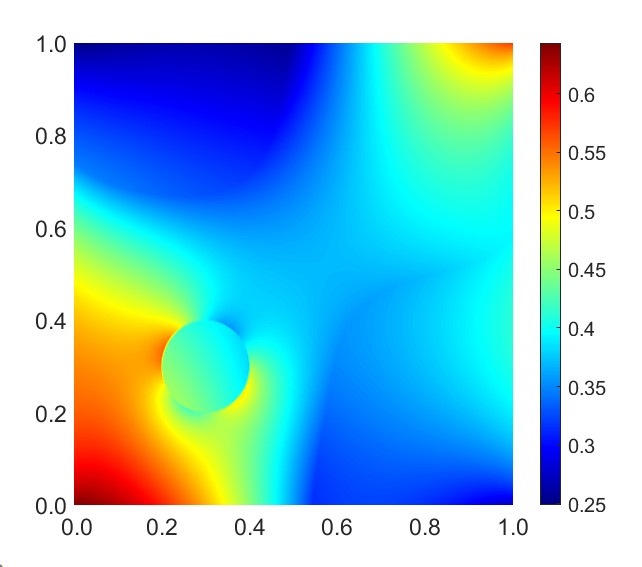}\\
		(d)$L=3$ & (e) $L=4$ & (f) $L=5$ \\
	\end{tabular}
	\caption{Boundary illuminations and the non-zero region of Example~\ref{ex:5}. Top left: plot of boundary data $\fl=g^{(\ell+1)}$. Top middle to bottom right: region which satisfying the non-zero condition as number of boundary input increasing.}
	\label{Fig:ex5_nonzero}
\end{figure}

\begin{figure}[htbp]
	\centering
	\begin{tabular}{ccc}
		\includegraphics[width=0.25\textwidth]{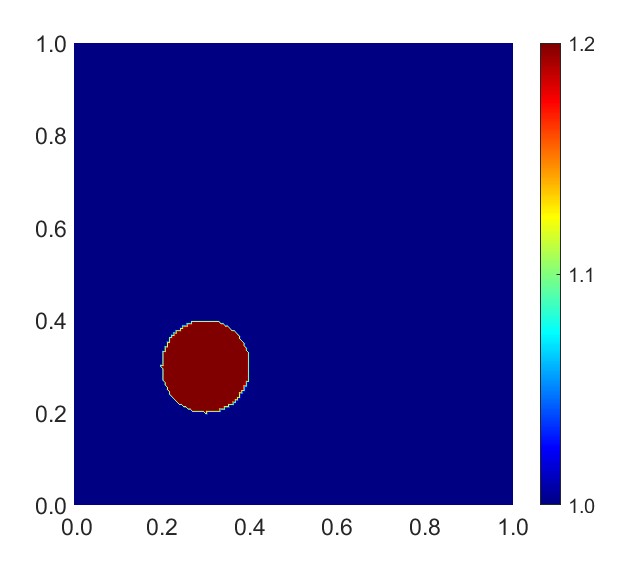}&
		\includegraphics[width=0.25\textwidth]{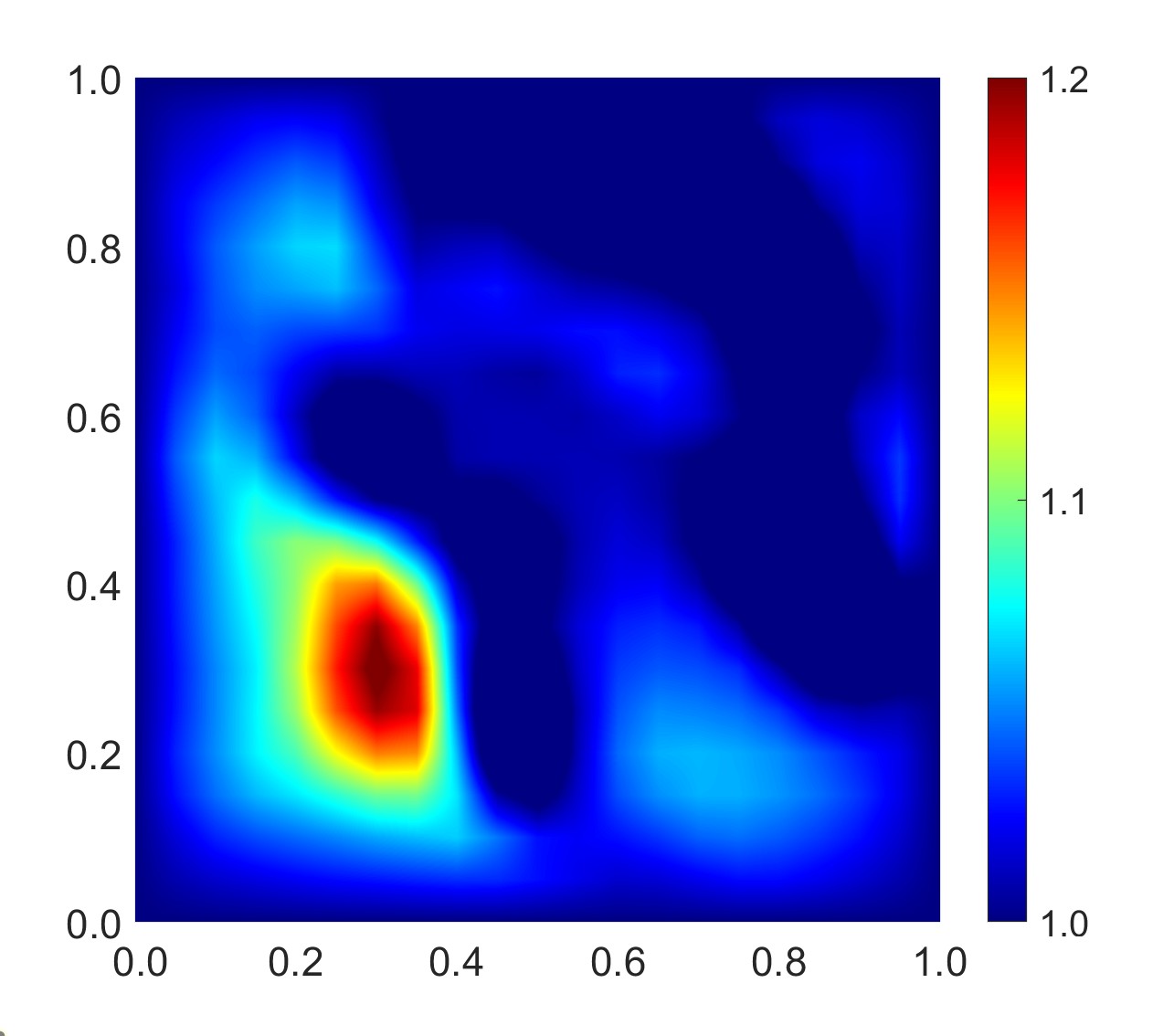}&
		\includegraphics[width=0.25\textwidth]{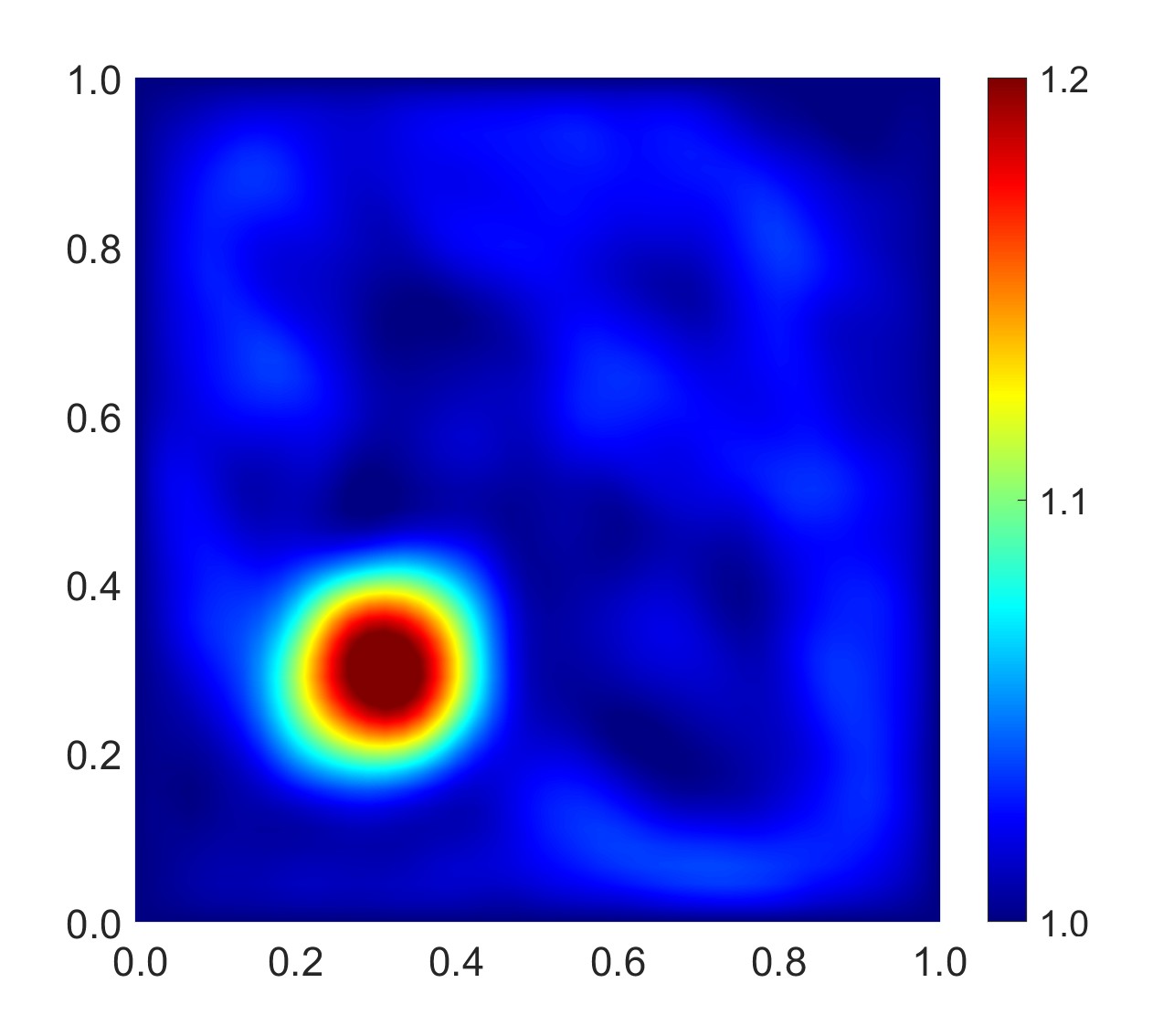}\\
		(a) exact $D$ & (b) $\delta=5e$-$2$ & (c) $\delta=1e$-$2$ \\
		\includegraphics[width=0.25\textwidth]{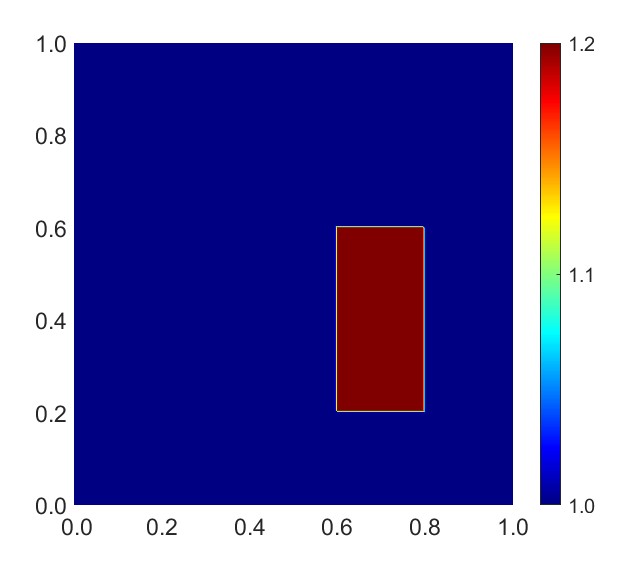}&
		\includegraphics[width=0.25\textwidth]{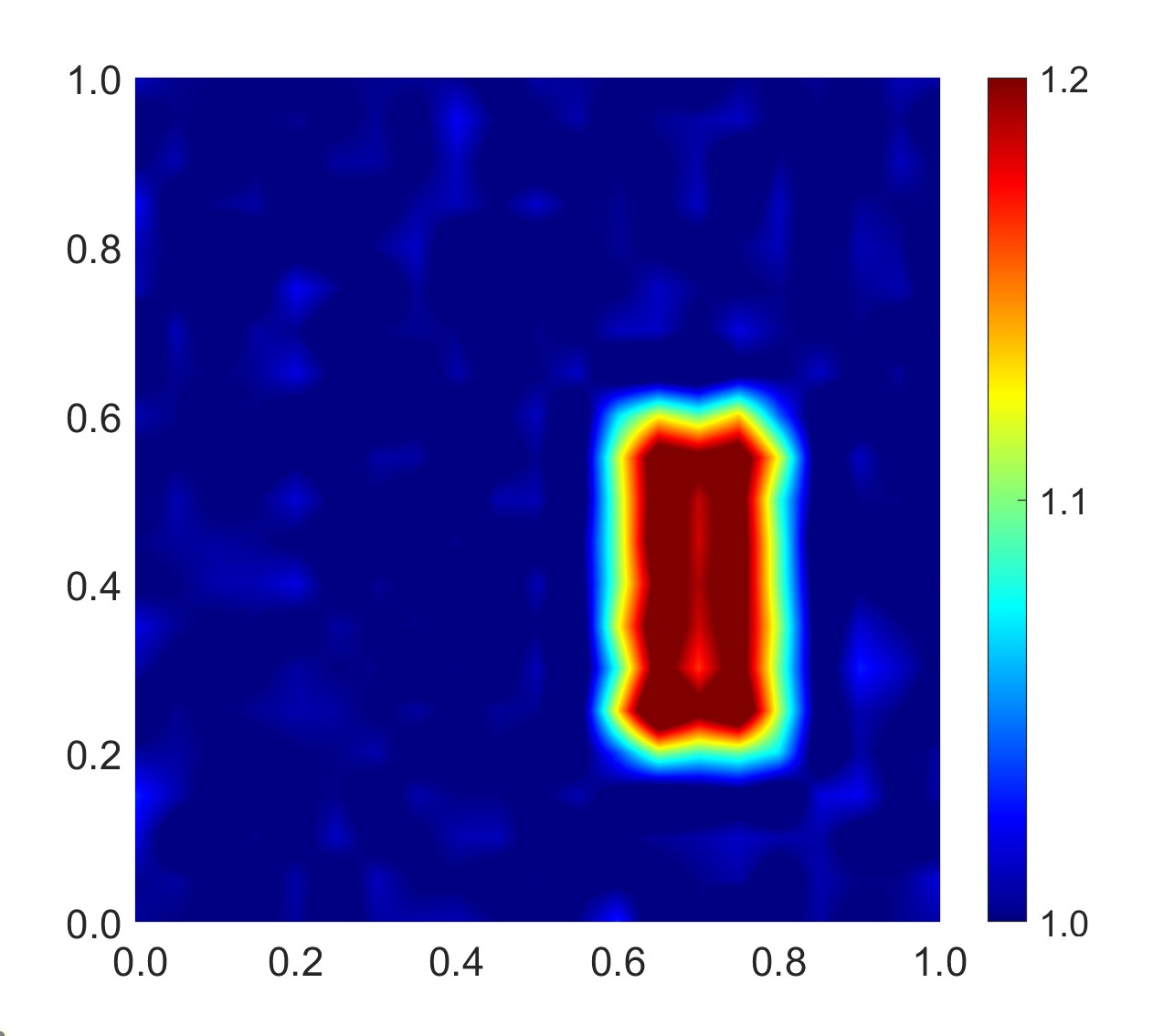}&
		\includegraphics[width=0.25\textwidth]{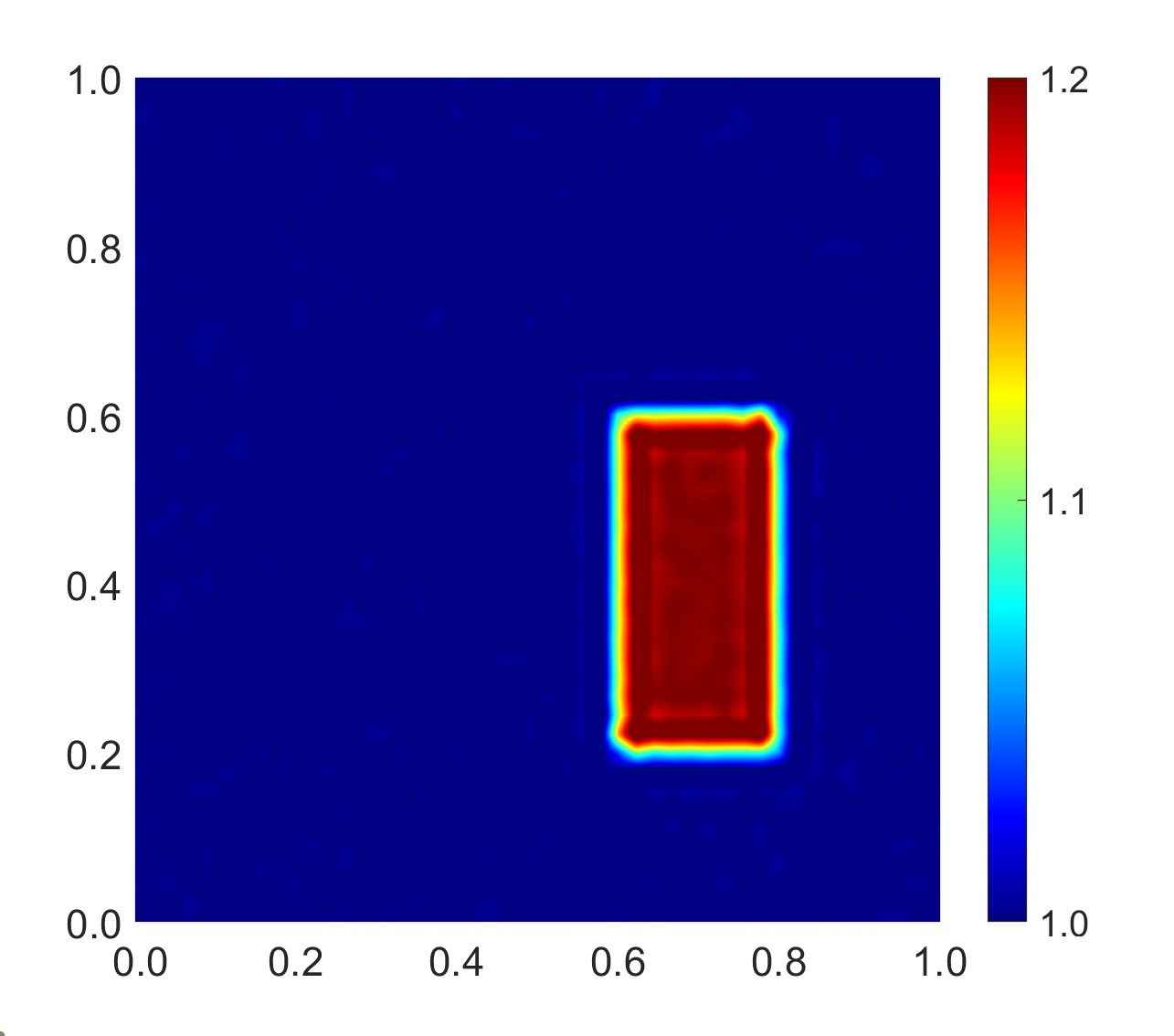}\\
		(d) exact $\sigma$ & (e) $\delta=5e$-$2$ & (f) $\delta=1e$-$2$ \\
	\end{tabular}
	\caption{Example \ref{ex:5}. First row: reconstructions of $D^\dag$.
		Second row: reconstructions of $\sigma^\dag$. }
	\label{Fig:ex5_recon}
\end{figure}

\begin{example}\label{ex:5}
	$\Omega=(0,1)^2$,  $D^{\dag}(x,y)=1+0.2\chi_{\{(x-0.3)^2+(y-0.3)^2<0.1^2 \} } $ and $\sigma^{\dag}(x,y)=1+0.2\chi_{[0.6,0.8]\times[0.2,0.6]} $. % The exact solution $\ul$ are computed  with mesh size $h=1/500$.
\end{example}
In this case, both the diffusion coefficient $D^\dag$ and the absorption coefficient $\sigma^\dag$ are piecewise constant, which is out the scope of our theoretical framework. Figures~\ref{Fig:ex5_nonzero}-\ref{Fig:ex5_recon} show the non-zero condition and the numerical reconstructions. The results indicate that the non-zero condition remains valid numerically even if the coefficients do not satisfy Assumption~\ref{assum:QPAT}. Meanwhile, the reconstructions are satisfactory for these piecewise constant coefficients. Here we take $h=1/20$, $\alpha=1e$-$6$ for noise level $\delta=5e$-$2$ and  $h=1/45$, $\alpha=4e$-$8$ for  $\delta=1e$-$2$. 

\section{Conclusion}\label{sec:conclusion}
In this paper, we investigated the reconstruction of the diffusion and absorption coefficients in QPAT. This is achieved by using multiple internal measurements illuminated by random boundary data. The reconstruction method starts with a straightforward reformulation, leading to an inverse diffusivity problem. A H\"older type stability is established by using energy estimates with special test function as well as the non-zero condition, guaranteed by the use of random boundary illuminations. The diffusivity coefficient is numerically recovered by employing a least-square formulation with a finite element discretization. The stability estimate motivates the approximation error analysis. With appropriate choices of the discretization mesh size and of the regularization parameters in relation with the noise level, the convergence rate of the approximation error is comparable to the stability result. In the subsequent step, we solve a direct problem involving the reconstructed diffusivity and {the} optical energy measurement. The diffusion and absorption coefficients can be recovered by an algebraic relation using the solution of the direct problem and the reconstructed diffusivity in the previous stage.

{Several important questions remain to be resolved in future work. First, while the use of random boundary illuminations is theoretically robust for satisfying the non-zero condition, it would be beneficial to quantify the minimal number of required illuminations $L$ and investigate how the choice of illuminations directly influences the quality of the reconstruction. Second, the current stability and error analyses rely on regularity assumptions that may not hold in practical QPAT applications, where coefficients often exhibit non-smooth or piecewise constant behavior reflecting heterogeneous material properties. Future work will focus on deriving a posteriori error estimates and designing adaptive finite element methods addressing these lower regularity settings. Finally, as our numerical experiments indicated that empirical convergence rates  are higher than the theoretical ones, we aim to investigate the underlying reasons for this discrepancy to develop sharper stability estimates and more representative error analyses.}

\section*{Funding}
Co-funded by the European Union (ERC, SAMPDE, 101041040). Views and opinions expressed are however those of the authors only and do not necessarily reflect those of the European Union or the European Research Council. Neither the European Union nor the granting authority can be held responsible for them. GSA is a member of the ``Gruppo Nazionale per l’Analisi Matematica, la Probabilità e le loro Applicazioni'', of the ``Istituto Nazionale di Alta Matematica''. The research was supported in part by the MIUR Excellence Department Project awarded to Dipartimento di Matematica, Università di Genova, CUP D33C23001110001. Finanziato dall’Unione europea-Next Generation EU, Missione 4 Componente 1 CUP D53D23005770006. The work of Z. Zhou is supported by by National Natural Science Foundation of China (Projects 12422117 and 12426312) and Hong Kong Research Grants Council (15302323).

\bibliographystyle{abbrv}
\bibliography{ref}
\end{document}